\documentclass[8 pt]{amsproc}
\usepackage{amsmath}
\usepackage{amssymb}
\usepackage{graphicx}
\usepackage{blindtext, rotating}
\usepackage[english]{babel}

\oddsidemargin=-1cm\evensidemargin=-1cm
\begin{document}
\numberwithin{equation}{section}

\def\1#1{\overline{#1}}
\def\2#1{\widetilde{#1}}
\def\3#1{\widehat{#1}}
\def\4#1{\mathbb{#1}}
\def\5#1{\frak{#1}}
\def\6#1{{\mathcal{#1}}}

\def\C{{\4C}}
\def\R{{\4R}}
\def\n{{\4n}}
\def\Z{{\4Z}}

 \title{Formal   Embeddings  between $\mathcal{BSD}$-Models}
 \author{Valentin Burcea}
\begin{abstract}  It is studied the Classification Problem for   Formal (Holomorphic)  Embeddings,  between (small open pieces of) Shilov Boundaries of Bounded Symmetric Domains of First  Type    in Complex Spaces of Different Dimensions, using linear changes of coordinates. 
\end{abstract}
\address{V. Burcea: INDEPENDENT}
\email{vdburcea@gmail.com}
\thanks{\emph{Keywords:} Normal Form,  Equivalence Problem,   Bounded Symmetric Domain, Real Submanifold, Shilov Boundary,  Power Series}
\thanks{This
project was supported principally ($80\%$)  by CAPES at  Federal University of Espirito Santo, being initiated ($20\%$) with   Science Foundation Ireland, Grant 10/RFP/MT H2878. Emphasizing that the paper \cite{V1} was fully supported by   Science Foundation Ireland, Grant 06/RFP/MAT 018.}
\maketitle
  
\def\Label#1{\label{#1}{\bf (#1)}~}


\def\cn{{\C^n}}
\def\cnn{{\C^{n'}}}
\def\ocn{\2{\C^n}}
\def\ocnn{\2{\C^{n'}}}


\def\dist{{\rm dist}}
\def\const{{\rm const}}
\def\rk{{\rm rank\,}}
\def\id{{\sf id}}
\def\tr{{\bf tr\,}}
\def\aut{{\sf aut}}
\def\Aut{{\sf Aut}}
\def\CR{{\rm CR}}
\def\GL{{\sf GL}}
\def\Re{{\sf Re}\,}
\def\Im{{\sf Im}\,}
\def\span{\text{\rm span}}
\def\Diff{{\sf Diff}}

\def\codim{{\rm codim}}
\def\crd{\dim_{{\rm CR}}}
\def\crc{{\rm codim_{CR}}}

\def\phi{\varphi}
\def\eps{\varepsilon}
\def\d{\partial}
\def\a{\alpha}
\def\b{\beta}
\def\g{\gamma}
\def\G{\Gamma}
\def\D{\Delta}
\def\Om{\Omega}
\def\k{\kappa}
\def\l{\lambda}
\def\L{\Lambda}
\def\z{{\bar z}}
\def\w{{\bar w}}
\def\Z{{\1Z}}
\def\t{\tau}
\def\th{\theta}

\emergencystretch15pt \frenchspacing

\newtheorem{Thm}{Theorem}[section]
\newtheorem{Cor}[Thm]{Corollary}
\newtheorem{Pro}[Thm]{Proposition}
\newtheorem{Lem}[Thm]{Lemma}

\theoremstyle{definition}\newtheorem{Def}[Thm]{Definition}

\theoremstyle{remark}
\newtheorem{Rem}[Thm]{Remark}
\newtheorem{Exa}[Thm]{Example}
\newtheorem{Exs}[Thm]{Examples}

\def\bl{\begin{Lem}}
\def\el{\end{Lem}}
\def\bp{\begin{Pro}}
\def\ep{\end{Pro}}
\def\bt{\begin{Thm}}
\def\et{\end{Thm}}
\def\bc{\begin{Cor}}
\def\ec{\end{Cor}}
\def\bd{\begin{Def}}
\def\ed{\end{Def}}
\def\br{\begin{Rem}}
\def\er{\end{Rem}}
\def\be{\begin{Exa}}
\def\ee{\end{Exa}}
\def\bpf{\begin{proof}}
\def\epf{\end{proof}}
\def\ben{\begin{enumerate}}
\def\een{\end{enumerate}}
\def\beq{\begin{equation}}
\def\eeq{\end{equation}}
 \section{Introduction and Main Result}
 
 The study of    proper  holomorphic mappings,  between unit balls  in Complex Spaces, goes back to Webster\cite{Webs}. For  $N>n$, two proper holomorphic mappings $f,g :\mathbb{B}^{n}\rightarrow\mathbb{B}^{N}$ are called equivalent if  there exist $\sigma\in \mbox{Aut}\left(\mathbb{B}^{n}\right)$ and  $\tau\in \mbox{Aut}\left(\mathbb{B}^{N}\right)$ such that $g=\tau\circ f\circ\sigma$. 
   
 Huang\cite{huang1} proved that any   proper holomorphic mapping,    between $\mathbb{B}^{n}$ and $\mathbb{B}^{N}$    of class $\mathcal{C}^{2}$ up to the boundary, is equivalent to  $$\mbox{ $\left( z_{1}, z_{2}, \dots,z_{n}\right)\longrightarrow\left( z_{1}, z_{2}, \dots,z_{n},0,\dots,0\right)$, for all $n>1$ and $N<2n-1$.}$$
 
Huang-Ji\cite{HJ} proved that the rational proper holomorphic mappings   between $\mathbb{B}^{n}$ and $\mathbb{B}^{2n-1}$  are equivalent to
 $$\left( z_{1}, z_{2}, \dots,z_{n}\right)\rightarrow \left( z_{1}, z_{2}, \dots,z_{n},0,\dots,0\right),\hspace{0.1 cm} \left(z_{1}, z_{2}, \dots, z_{n-1}, z_{n}z_{1}, z_{n}z_{2}, \dots, z_{n}^{2}\right)
,\quad\mbox{for all
  $n\geq3$.}$$

The Classification Problem, for proper  holomorphic mappings\cite{Far},\cite{huang1},\cite{HJ} between unit balls in Complex Spaces,  is reduced to the study  of   CR mappings between  hyperquadrics (see \cite{BH},\cite{BEH},\cite{LP2}). More generally, the   Classification Problem for C.-R. Mappings, between Shilov Boundaries of Bounded Symmetric Domains of First   Type,   is also very interesting. It has been studied by Kim-Zaitsev\cite{kz}   using the   moving frames  method of Cartan.   Given $q<p$, $q'<p'$ such that $p'-q'<2\left(p-q\right)$ and $p-q>1$, they\cite{kz} proved that,  up to compositions with suitable   automorphisms of   Bounded Symmetric Domains of First Type denoted by $D_{p,q}$ and $D_{p',q'}$, any Smooth  C.-R. Embedding,   between (small open pieces of) their Shilov Boundaries denoted by $S_{p,q}$ and $S_{p',q'}$, is equivalent to the   geodesical emebdding.

 Any Bounded  Symmetric Domain $D_{p,q}$ of First  Type and   its Shilov Boundary (see \cite{kz})  may be defined   by  
\begin{equation}D_{p,q}=\left\{Z\in\mathcal{M}_{p\times q}\left(\mathbb{C}\right);\quad I_{q}-\overline{Z}^{t}Z> 0\right\}\hspace{0.1 cm} \mbox{and} \hspace{0.1 cm}  S_{p,q}=\left\{Z\in\mathcal{M}_{p\times q}\left(\mathbb{C}\right);\quad I_{q}-\overline{Z}^{t}Z= 0\right\},\quad\mbox{for $p>q$.} \label{shi} \end{equation}
 
It is indicated Chirvasitu\cite{Chir} for the standard definition of the Shilov Boundary.  Generalizing classical models as the hyperquadric  (see    \cite{H},\cite{huang1},\cite{huang2},\cite{HJ},\cite{LP2},\cite{JX}), it is fundamental   in the  study  of  Holomorphic Isometries in  Complex Geometry (see \cite{Isom1},\cite{Isom2},\cite{Isom3}). 

In this paper, we use   the language of matrices  in order to establish a normal form  (see also \cite{BEH},\cite{V1},\cite{D1},\cite{D3}) for   Formal (Holomorphic)  Embeddings between Shilov Boundaries of Bounded  Symmetric Domains of First Type (see \cite{kz},\cite{kza}). They are reduced to   the local  study of (local)  Formal Embeddings between   $\mathcal{BSD}$-Models, which are just   quadratic Models derived    from   Shilov Boundaries of Bounded   Symmetric Domains of First Type using   generalized   Cayley Transformations (see \cite{EG}) as 
\begin{equation} \mathcal{C}^{-1}\left(W,Z\right) = \left( \frac{\left[ \sqrt{-1}\left( I_{q}-W\right),\sqrt{-1}Z \right]}{W+ I_{q}}\right)^{t},\quad\mbox{where $W=\left\{w_{ij}\right\}_{1\leq i,j\leq q}$ and   $Z=\left\{z_{ij}\right\}_{1\leq i\leq q\atop{1\leq j\leq N}}$, for $N=p-q$.}\label{Ca}
\end{equation}

In particular,   we use  standard linearization procedure  of the  constructions of normal forms (see  \cite{BEH},\cite{V1},\cite{CM},\cite{D1},\cite{D3}) in order to make computations using formal  power series.     It is an   alternative   to the Method of Cartan  applied by Kim-Zaitsev\cite{kz},\cite{kza}, but it involves complicated methods of linear algebra due to the existence of   non-trivial classes  of mappings derived from the existence of  an analogue   of  the geometrical rank (see \cite{huang1},\cite{huang2}) named    rank in this paper. The geometrical rank may be $0$ or $1$ according to Huang-Ji\cite{HJ} in our case.

The formal computations are reduced to the study of  Formal Embeddings between quadratic manifolds (see also \cite{H}) named $\mathcal{BSD}$-Models. We implement sophisticated linear changes of coordinates    in order to  apply procedures from Baouendi-Huang\cite{BH},   Chern-Moser\cite{CM},   Huang-Ji\cite{HJ} and the moving point trick inspired by Baouendi-Ebenfelt-Huang\cite{BEH} and motivated by Baouendi-Huang\cite{BH} and  Huang\cite{huang1}.  

 The main result is the following: 

 \bt\label{Teorem} Let $p,p',q,q'\in\mathbb{N}^{\star}$   such that $p'-q'=2\left(p-q\right)>2$, $q<q'$  and $p<p'$. It exists $k_{0}\in\left\{0,1,\dots,q\right\}$, and a   set of polynomials of degree $2$, denoted by $\left\{\mathcal{L}_{il}^{k}\left(Z_{i},Z_{k}\right)\right\}_{1\leq i \leq k_{0}\atop{1 \leq k  \leq q \atop{1 \leq l \leq N}}}$, a set    of polynomials of degree $1$, denoted by  $\left\{a_{ku}^{il}\left(Z_{i}\right)\right\}_{1\leq i \leq k_{0}\atop{1 \leq k  \leq q \atop{1 \leq l \leq N}}}$,   such that up to compositions with suitable   automorphisms of the bounded symmetric domains     $D_{p,q}$ and $D_{p',q'}$, a  Formal  Embedding,  between (small open pieces of) their Shilov Boundaries $S_{p,q}$ and $S_{p',q'}$, is equivalent to  
 
\begin{equation}\left(W,Z\right) \rightarrow\begin{pmatrix}  \begin{pmatrix} W   
&  \left(W+ I_{q}\right)\cdot \frac{ \begin{pmatrix}  \tilde{A}_{k_{0}}\left(\left( \frac{\left[ \sqrt{-1}\left( I_{q}-W\right),\sqrt{-1}Z \right]}{W+ I_{q}}\right)^{t}\right)  
&       D_{k_{0}}\left(\left( \frac{\left[ \sqrt{-1}\left( I_{q}-W\right),\sqrt{-1}Z \right]}{W+ I_{q}}\right)^{t}\right)\end{pmatrix}
   }{ \sqrt{-1} }
 \end{pmatrix} & \mbox{O}_{p\times \left(q-p\right)} \\    \mbox{O}_{\left(p'-p\right)\times \left(2q-p\right)}   & \mbox{I}_{\left(q'-q\right)\times \left(q'-q\right)}     \end{pmatrix}, \label{clase} 
\end{equation}
where    we have used the     matrices
 \begin{equation*} 
 \tilde{A}_{k_{0}}\left(W,Z\right)=\begin{pmatrix}  
\displaystyle\sum_{k,u=1\atop{(k,u)\in\mathcal{S}}}^{q}\left(a_{ku}^{il}\left(Z_{i}\right)w_{ku}\right)_{ 1\leq i\leq q\atop{1\leq l\leq 2N}}   \\    \begin{pmatrix}
    z_{k_{0}+1\hspace{0.01 cm}1}   & z_{k_{0}+1\hspace{0.01 cm}2} & \dots & z_{k_{0}+1\hspace{0.01 cm}N} \\ \vdots & \vdots & \ddots & \vdots  \\  z_{q1}    & z_{q2} & \dots & z_{qN}   \end{pmatrix}\end{pmatrix}\hspace{0.1 cm}\mbox{and} \hspace{0.1 cm}  D_{k_{0}}\left(W,Z\right) = \begin{pmatrix}
\left(  \displaystyle\sum_{k=1}^{q}\mathcal{L}_{il}^{k}\left(Z_{i},Z_{k}\right) \right)_{ 1\leq i\leq q\atop{1\leq l\leq 2N}}   
\\ \mbox{O}_{\left(q-k_{0}\right)\times N}
\end{pmatrix} ,
\end{equation*}
such that the following conditions of compatibility hold
\begin{equation} \left(\left\langle\displaystyle\sum_{k,u=1\atop{(k,u)\in\mathcal{S}}}^{q}\left(a_{ku}^{il}\left(Z_{i}\right)w_{ku}\right)_{ 1\leq i\leq k_{0}\atop{1\leq l\leq N}},Z_{j}\right\rangle+\left\langle \left( \displaystyle\sum_{k=1}^{q}\mathcal{L}_{il}^{k}\left(Z_{i},Z_{k}\right)\right)_{1\leq l \leq N},\left(\displaystyle\sum_{k=1}^{q}\mathcal{L}_{jl}^{k}\left(Z_{j},Z_{k}\right) \right)_{1\leq l \leq N}    \right\rangle \right)_{1\leq i,j\leq k_{0} }=\mbox{O}_{k_{0}\times k_{0}} ,\label{compa} 
 \end{equation}
where we have used the set 
 \begin{equation*}\mathcal{S}:=\left\{(k,i)\in \left\{1,\dots,k_{0}\right\}\times \left\{1,\dots,q\right\} \left\vert \left(a_{kk}^{il}\left(Z_{i}\right)\right)_{1 \leq l  \leq N }  \not\equiv (0,0,\dots,0)   \right.\right\},
\end{equation*}  
where $\left\langle \cdot,\cdot\right\rangle$ is the standard hermitian product, and $Z_{1},Z_{2},\dots,Z_{q}$ are the row vectors of the matrix $Z$.
\et 
 
The set of polynomials of degree $1$ and $2$ verify   compatibility conditions defined by elements of the set $\mathcal{S}$.       Generalized Cayley type transformations are used in order to obtain the non-equivalent  classes of equivalence (\ref{clase}).   The    standard linear Embedding described by Kim-Zaitsev\cite{kz} corresponds to the case   $k_{0}=0$.     The generalized Whitney type mapping  introduced by Seo\cite{Seo2} must be equivalent to a certain class described by (\ref{clase}),  but our result is locally proven. An immediate application is the following:  any local proper (formal) holomorphic maps is equivalent to a class defined by (\ref{clase}) in the hypothesis of  Theorem \ref{Teorem}.

 \section{Settings}
\subsection{Identifications}
Throughout this paper, we work in  the following coordinates  
 \begin{equation}\left(w_{11},w_{12},\dots,w_{1q},\dots,w_{q1},w_{q2},\dots,w_{qq};z_{11},z_{12},\dots,z_{1N},\dots,z_{q1},z_{q2},\dots,z_{qN}\right)\in \mathbb{C}^{qN+q^{2}}.\label{coord}\end{equation}   
 \begin{equation}\left({w'}_{11},{w'}_{12},\dots,{w'}_{1{q'}},\dots,{w'}_{{q'}1},{w'}_{{q'}2},\dots,{w'}_{{q'}{q'}};{z'}_{11},{z'}_{12},\dots,{z'}_{1N'},\dots,{z'}_{{q'}1},{z'}_{{q'}2},\dots,{z'}_{{q'}N'}\right)\in \mathbb{C}^{{q'}N'+{q'}^{2}}.\label{coord1}\end{equation} 
 
 In particular, we  consider matrices according to the following   identifications
 \begin{equation}\begin{split}& W=\begin{pmatrix}  w_{11} & w_{12} &\dots &w_{1q} \\ w_{21} & w_{22} &\dots & w_{2q} \\ \vdots &\vdots &\ddots &\vdots  \\ w_{q1} & w_{q2} &\dots & w_{qq}\end{pmatrix}\equiv\left(w_{11},w_{12},  \dots, w_{1q},w_{21},w_{22},  \dots, w_{2q},\dots\dots,w_{q1},w_{q2},  \dots, w_{qq}\right),\\&    Z=\begin{pmatrix} z_{11} &z_{12}& \dots &z_{1N} \\ z_{21} &z_{22}& \dots &z_{2N} \\  \vdots &\vdots &\ddots &\vdots \\ z_{q1} &z_{q2}&\dots &z_{qN}\end{pmatrix} \equiv\left(z_{11}, z_{12}, \dots,z_{1N},z_{21}, z_{22}, \dots,z_{2N},\dots\dots, z_{q1}, z_{q2}, \dots,z_{qN}\right).\end{split}\label{nota}\end{equation}

In consequence, we consider the following    identifications 
\begin{equation}\begin{split}& J:=\begin{pmatrix} j_{11} & j_{12} &\dots &j_{1q} \\ j_{21} & j_{22} &\dots &j_{2q} \\ \vdots & \vdots &\ddots &\vdots  \\ j_{q1} &j_{q2} &\dots &j_{qq}\end{pmatrix}\equiv\left(j_{11},j_{12},  \dots, j_{1q},j_{21},j_{22},  \dots, j_{2q},\dots\dots,j_{q1},j_{q2},  \dots, j_{qq}\right)\in\mathbb{N}^{q^{2}},\\&   I:=\begin{pmatrix} i_{11} &i_{12}& \dots&i_{1N} \\ i_{21} &i_{22}& \dots&i_{2N} \\\vdots &\vdots &\ddots &\vdots \\ i_{m1} &i_{m2}&\dots&i_{qN}\end{pmatrix} \equiv\left(i_{11},i_{12},  \dots,i_{1N},i_{21},i_{22},  \dots,i_{2N},
\dots\dots,i_{q1},i_{q2},  \dots,i_{qN}\right)\in\mathbb{N}^{qN}.\end{split}\label{notaa}\end{equation}

 \begin{equation}  \begin{split} \left\{1,2,3,\dots,q^{2} \right\}\equiv & \left\{ (1,1),(1,2),\dots,\left(1,q\right)\right.,\\&  \hspace{0.2 cm} (2,1),(2,2),\dots,\left(2,q\right),
\\&\quad\hspace{0.2 cm}\vdots\hspace{0.1 cm}\quad\quad\vdots\quad\hspace{0.2 cm}\ddots\quad\hspace{0.1 cm}\vdots\\& \left.\hspace{0.2 cm}(q,1),(q,2),\dots,\left(q,q\right)\right\} \end{split} \quad\mbox{and}  \quad \begin{split}\left\{1,2,3,\dots,qN \right\}\equiv & \left\{ (1,1),(1,2),\dots,\left(1,N\right)\right.,\\&  \hspace{0.2 cm} (2,1),(2,2),\dots,\left(2,N\right),
\\&\quad\hspace{0.2 cm}\vdots\hspace{0.1 cm}\quad\quad\vdots\quad\hspace{0.2 cm}\ddots\quad\hspace{0.1 cm}\vdots\\& \left.\hspace{0.2 cm}(q,1),(q,2),\dots,\left(q,N\right)\right\}.\end{split}\label{GG2}
\end{equation}

\begin{equation}  \begin{split} \left\{1,2,3,\dots,{q'}^{2} \right\}\equiv & \left\{ (1,1),(1,2),\dots,\left(1,q'\right)\right.,\\&  \hspace{0.2 cm} (2,1),(2,2),\dots,\left(2,q'\right),
\\&\quad\hspace{0.2 cm}\vdots\hspace{0.1 cm}\quad\quad\vdots\quad\hspace{0.2 cm}\ddots\quad\hspace{0.1 cm}\vdots\\& \left.\hspace{0.2 cm}\left(q',1\right),\left(q',2\right),\dots,\left(q',q'\right)\right\} \end{split} \quad\mbox{and}  \quad \begin{split}\left\{1,2,3,\dots,q'N'\right\}\equiv & \left\{ (1,1),(1,2),\dots,\left(1,N'\right)\right.,\\&  \hspace{0.2 cm} (2,1),(2,2),\dots,\left(2,N'\right),
\\&\quad\hspace{0.2 cm}\vdots\hspace{0.1 cm}\quad\quad\vdots\quad\hspace{0.2 cm}\ddots\quad\hspace{0.1 cm}\vdots\\& \left.\hspace{0.2 cm}\left(q',1\right),\left(q',2\right),\dots,\left(q',N'\right)\right\}.\end{split}\label{GG222}
\end{equation}

In particular,  we write  
\begin{equation}\begin{split}
W^{J}=&w_{11}^{j_{11}}w_{12}^{j_{12}}\dots w_{1q}^{j_{1q}}\cdot \\&  w_{21}^{j_{21}}w_{22}^{j_{22}}\dots w_{2q}^{j_{2q}}\cdot \\& \quad \vdots  \quad \vdots  \quad \ddots  \quad \vdots\\& w_{q1}^{j_{q1}}w_{q2}^{j_{q2}}\dots w_{qq}^{j_{qq}}\end{split}\quad\mbox{and}  \quad \begin{split} Z^{I}=& z_{11}^{i_{11}}z_{12}^{i_{12}}\dots z_{1N}^{i_{1N}}\cdot \\& z_{21}^{i_{21}}z_{22}^{i_{22}}\dots z_{2N}^{i_{2N}}\cdot \\& \quad \vdots  \quad \vdots  \quad \ddots  \quad \vdots \\&  z_{q1}^{i_{q1}}z_{q2}^{i_{q2}}\dots z_{qN}^{i_{qN}}.\end{split}
\label{IJ}\end{equation}

 We   use similarly    matrices   denoted by $W'$ and  $Z'$. The rows  of the  matrices $Z$ and $Z'$   are denoted by
 \begin{equation} Z_{1},Z_{2},\dots,Z_{q}\hspace{0.1 cm}\mbox{and}\hspace{0.1 cm} Z'_{1},Z'_{2},\dots,Z'_{q'}.\label{row} 
\end{equation}

\subsection{Matrices and Models}Ignoring the considered  numbers, we   generally write
\begin{equation}\left<L, V\right>=L\overline{V}^{t},\quad\mbox{given the matrices $L \in \mathcal{M}_{m,n}\left(\mathbb{C}\right)$ and $V \in \mathcal{M}_{n,p}\left(\mathbb{C}\right)$, for all $m,n,p\in\mathbb{N}^{\star}$.}\label{vb} 
\end{equation}
 
  We replace   the  generalized Cayley transformation defined in (\ref{Ca}) (see \cite{EG}) in (\ref{shi}). We obtain the equation of the  $\mathcal{BSD}$-Model  
  \begin{equation}  \frac{W-\overline{W}^{t}}{2\sqrt{-1}} =Z\overline{Z}^{t}.
\label{tr}\end{equation}

Any (local)  Formal    Embedding, denoted by $\left(\tilde{G},\tilde{F}\right)$, defined between   Shilov Boundaries  of Bounded   Symmetric Domains of First Type, induces naturally by (\ref{Ca})   another (local)  Formal    Embedding, denoted by  $\left(G,F\right)$,   between  the  $\mathcal{BSD}$-Models  defined    by
\begin{equation}\begin{split}&\hspace{0.1 cm}\mathcal{M}:\hspace{0.01 cm}\mbox{Im} W=Z\overline{Z}^{t}\subset \mathbb{C}^{qN+q^{2}},\quad\quad\hspace{0.26 cm}\mbox{for $N=p-q$,}  \\&  \mathcal{M}':\hspace{0.01 cm} \mbox{Im} W'=Z'\overline{{Z'}}^{t}\subset\mathbb{C}^{q'N'+{q'}^{2}},\quad\mbox{for   $N'=p'-q'$,} \end{split} \label{models}
\end{equation}
where  $q<q'$  and $p<p'$. In particular,    we work with  the  commutative diagram

\begin{equation}\begin{array}[c]{ccc}
\mathcal{M}&\stackrel{\left(G,F\right)}{\rightarrow}&\mathcal{M}'\\
 \Updownarrow\scriptstyle{ }&&\Updownarrow\scriptstyle{}\\
S_{p,q}&\stackrel{\left(\tilde{G},\tilde{F}\right)}{\rightarrow}&S_{p',q'},
\end{array}\label{diag}\end{equation}
where each equivalence  is defined using (\ref{Ca}). The   Formal Embedding $\left(G,F\right)$  is written as
\begin{equation}\begin{split}&G\left(W,Z\right):=   \begin{pmatrix} G_{11}\left(W,Z\right)& G_{12}\left(W,Z\right) \\ G_{21}\left(W,Z\right) & G_{22}\left(W,Z\right)\end{pmatrix}\hspace{0.1 cm}\mbox{and}\hspace{0.1 cm} F\left(W,Z\right):=   \begin{pmatrix} F_{1}\left(W,Z\right) \\ F_{2}\left(W,Z\right) \end{pmatrix},  \label{v11}\end{split} \end{equation}
where we have used   the following submatrices:
\begin{itemize}
\item $G_{11}\left(W,Z\right)$ is a $q\times q$ submatrix  having  formal power series in $\left(W,Z\right)$ as   entries,
\item $G_{21}\left(W,Z\right)$ is a $\left(q'-q\right)\times q$ submatrix having  formal power series in $\left(W,Z\right)$ as entries,
\item $G_{12}\left(W,Z\right)$ is a $q\times \left(q'-q\right)$ submatrix having  formal power series in $\left(W,Z\right)$ as entries,
\item $G_{22}\left(W,Z\right)$ is  a $\left(q'-q\right)\times \left(q'-q\right)$ submatrix having  formal power series  in $\left(W,Z\right)$ as entries,
\item  $F_{1}\left(W,Z\right)$ is  a $q\times  N'$ submatrix having  formal power series in $\left(W,Z\right)$  as entries,
\item $F_{2}\left(W,Z\right)$ is  a $\left(q'-q\right)\times N'$ submatrix having  formal power series in $\left(W,Z\right)$ as entries,
\end{itemize}
such that their   linear parts in $Z$ are vanishing. In particular, we obtain
\begin{equation}\begin{split}&\frac{G_{11}\left(W,Z\right)-\overline{G_{11}\left(W,Z\right)}^{t}}{2\sqrt{-1}}= \left<F_{1}\left(W,Z\right) ,F_{1}\left(W,Z\right) \right> ,\quad \frac{G_{12}\left(W,Z\right)-\overline{G_{21}\left(W,Z\right)}^{t}}{2\sqrt{-1}}=  \left<F_{1}\left(W,Z\right) ,F_{2}\left(W,Z\right) \right>, \\&  \frac{G_{21}\left(W,Z\right)-\overline{G_{12}\left(W,Z\right)}^{t}}{2\sqrt{-1}}=  \left<F_{2}\left(W,Z\right) ,F_{1}\left(W,Z\right) \right>,\quad   \frac{G_{22}\left(W,Z\right)-\overline{G_{22}\left(W,Z\right)}^{t}}{2\sqrt{-1}}=  \left<F_{2}\left(W,Z\right) ,F_{2}\left(W,Z\right) \right>    .\end{split}\label{v}
\end{equation}

The equations (\ref{v})    are  the subject of  transformations  preserving $\mathcal{BSD}$-Models using the    product of matrices
\begin{equation} V\otimes Z=\left(\displaystyle\sum_{l=1}^{N}\displaystyle\sum_{k=1}^{q}v_{kl}^{ij}z_{kl}\right)_{1\leq i\leq q\atop 1\leq j\leq N},\quad\mbox{given the matrix $V=\left( v_{\alpha}^{\beta} \right)_{1\leq\alpha\leq qN}^{1\leq \beta\leq qN}\in\mathcal{M}_{qN\times qN}\left(\mathbb{C}\right)$  such that:}\label{edef}
\end{equation}
\begin{equation*}V\equiv  \begin{pmatrix}  {\begin{pmatrix} v^{11}_{11} & v^{11}_{12} & \dots & v^{11}_{1N}\\    v^{12}_{11} & v^{12}_{12} & \dots & v^{12}_{1N} \\  \vdots & \vdots & \ddots & \vdots \\  v^{1N}_{11} & v^{1N}_{12} & \dots & v^{1N}_{1N} \end{pmatrix}} &  {\begin{pmatrix} v^{11}_{21} & v^{11}_{22} & \dots & v^{11}_{2N}\\    v^{12}_{21} & v^{12}_{22} & \dots & v^{12}_{2N} \\  \vdots & \vdots & \ddots & \vdots \\  v^{1N}_{21} & v^{1N}_{22} & \dots & v^{1N}_{2N} \end{pmatrix}} & \begin{pmatrix}
\dots \\ \dots \\ \vdots \\ \dots
\end{pmatrix}  & {\begin{pmatrix} v^{11}_{q1} & v^{11}_{q2} & \dots & v^{11}_{qN}\\    v^{12}_{q1} & v^{12}_{q2} & \dots & v^{12}_{qN} \\  \vdots & \vdots & \ddots & \vdots \\  v^{1N}_{q1} & v^{1N}_{q2} & \dots & v^{1N}_{qN} \end{pmatrix}} \\ {\begin{pmatrix} v^{21}_{11} & v^{21}_{12} & \dots & v^{21}_{1N}\\    v^{22}_{11} & v^{22}_{12} & \dots & v^{22}_{1N} \\  \vdots & \vdots & \vdots & \vdots \\  v^{2N}_{11} & v^{2N}_{12} & \dots & v^{2N}_{1N} \end{pmatrix}} &  {\begin{pmatrix} v^{21}_{21} & v^{21}_{22} & \dots & v^{21}_{2N}\\    v^{22}_{21} & v^{22}_{22} & \dots & v^{22}_{2N} \\  \vdots & \vdots & \ddots & \vdots \\  v^{2N}_{21} & v^{2N}_{22} & \dots & v^{2N}_{2N} \end{pmatrix}} & \begin{pmatrix}
\dots \\ \dots \\ \vdots \\ \dots
\end{pmatrix}  & {\begin{pmatrix} v^{21}_{q1} & v^{21}_{q2} & \dots & v^{21}_{qN}\\    v^{22}_{q1} & v^{22}_{q2} & \dots & v^{22}_{qN} \\  \vdots & \vdots & \ddots & \vdots \\  v^{2N}_{q1} & v^{2N}_{q2} & \dots & v^{2N}_{qN} \end{pmatrix}}\\ {\begin{pmatrix} \vdots & \vdots & \ddots & \vdots
\end{pmatrix}} &{\begin{pmatrix} \vdots & \vdots & \ddots & \vdots
\end{pmatrix}}& \begin{pmatrix}\vdots\end{pmatrix} & {\begin{pmatrix} \vdots & \vdots & \ddots & \vdots
\end{pmatrix}}    \\ {\begin{pmatrix} v^{q1}_{11} & v^{q1}_{12} & \dots & v^{q1}_{1N}\\    v^{q2}_{11} & v_{12}^{q2} & \dots & v^{q2}_{1N} \\   \vdots  & \vdots & \ddots & \vdots \\  v^{qN}_{11} & v^{qN}_{12} & \dots & v^{qN}_{1N} \end{pmatrix}} &  {\begin{pmatrix} v^{q1}_{q1} & v^{q1}_{22} & \dots & v^{q1}_{2N}\\    v^{q2}_{21} & v^{q2}_{22} & \dots & v^{q2}_{2N} \\  \vdots & \vdots & \ddots & \vdots \\  v^{qN}_{21} & v^{qN}_{22} & \dots & v^{qN}_{2N} \end{pmatrix}} & \begin{pmatrix}
\dots \\ \dots \\ \vdots \\ \dots
\end{pmatrix}  & {\begin{pmatrix} v^{q1}_{q1} & v^{q1}_{q2} & \dots & v^{q1}_{qN}\\    v^{q2}_{q1} & v^{q2}_{q2} & \dots & v^{q2}_{qN} \\  \vdots & \vdots & \ddots & \vdots \\  v^{qN}_{q1} & v^{qN}_{q2} & \dots & v^{qN}_{qN} \end{pmatrix}}   \end{pmatrix}.
\end{equation*}   

 The matrices    (\ref{v11}) may be rewritten      as
\begin{equation}  G\left(Z,W\right)= \left(g_{kl}\left(Z,W\right)\right)_{1\leq k,l\leq q'}\hspace{0.1 cm}\mbox{and}\hspace{0.1 cm} F\left(Z,W\right)= \left(f_{kl}\left(Z,W\right)\right)_{ 1\leq k \leq q'  \atop  1\leq  l\leq N'}. \label{10v}  \end{equation} 
  
 In particular, we   have
\begin{equation}\begin{split}& w_{11}=\overline{w}_{11}+2\sqrt{-1}\left<Z_{1},Z_{1}\right>,\quad w_{12} =\overline{w}_{21}+2\sqrt{-1}\left<Z_{1},Z_{2}\right>  , \dots\quad w_{1q} =\overline{w}_{q1}+2\sqrt{-1}\left<Z_{1},Z_{q}\right>, \\&  w_{21} =\overline{w}_{12}+2\sqrt{-1}\left<Z_{2},Z_{1}\right>,\quad w_{22}=\overline{w}_{21}+2\sqrt{-1}\left<Z_{2},Z_{2}\right>  , \dots\quad w_{2q} =\overline{w}_{q2}+2\sqrt{-1}\left<Z_{2},Z_{q}\right>,   \\&  \quad\quad \hspace{0.18 cm}    \vdots \quad\quad\quad\quad\quad\quad \quad\quad\quad \quad \quad\quad\quad\hspace{0.18 cm} \vdots \quad\quad \quad\quad\quad\quad\quad\quad\quad\quad   \ddots   \quad\quad\quad\hspace{0.18 cm}    \vdots  \\& w_{q1} =\overline{w}_{1q}+2\sqrt{-1}\left<Z_{q},Z_{1}\right>, \quad w_{q2} =\overline{w}_{2q}+2\sqrt{-1}\left<Z_{q},Z_{2}\right> ,\dots\quad   w_{qq} =\overline{w}_{qq}+2\sqrt{-1}\left<Z_{q},Z_{q}\right>.\end{split}  \label{bq}
\end{equation}
 
 \begin{equation}\begin{split}&\hspace{0.1 cm} w'_{11} =\overline{w}'_{11}+2\sqrt{-1}\left<Z'_{1},Z'_{1}\right>,\quad\quad w'_{12} =\overline{w}'_{21}+2\sqrt{-1}\left<Z'_{1},Z'_{2}\right>  , \dots\quad\quad\hspace{0.15 cm} w'_{1q'} =\overline{w}_{q'1}+2\sqrt{-1}\left<Z'_{1},Z'_{q'}\right>, \\& \hspace{0.1 cm} w'_{21} =\overline{w}'_{12}+2\sqrt{-1}\left<Z'_{2},Z'_{1}\right>,\quad\quad w'_{22} =\overline{w}'_{21}+2\sqrt{-1}\left<Z'_{2},Z'_{2}\right>  , \dots\quad\quad\hspace{0.15 cm} w'_{2q'} =\overline{w}'_{q'2}+2\sqrt{-1}\left<Z'_{2},Z'_{q'}\right>,   \\&  \quad\quad \hspace{0.18 cm}    \vdots \quad\quad\quad\quad\quad\quad\quad \quad\quad\quad \quad\quad \quad\quad\quad\hspace{0.18 cm} \vdots \quad\quad \quad\quad\quad\quad\quad\quad\quad\quad\quad\quad   \quad \quad\quad\quad\quad   \vdots  \\& w'_{q'1} =\overline{w}'_{1q'}+2\sqrt{-1}\left<Z'_{q'},Z'_{1}\right>, \quad w'_{q'2} =\overline{w}'_{2q'}+2\sqrt{-1}\left<Z'_{q'},Z'_{2}\right> ,\dots\quad   w'_{q'q'} =\overline{w}'_{q'q'}+2\sqrt{-1}\left<Z'_{q'},Z'_{q'}\right>.\end{split}  \label{bqprim}
\end{equation}

We more forward in order    to implement the following:
\section{ Application of the Normalization Procedure  from Baouendi-Huang\cite{BH}}

We study    the  Embedding  (\ref{v11}) using the  commutative  diagram  
   \begin{equation}\begin{array}[c]{ccc}
\mathcal{M}&\stackrel{\left(F,G\right)}{\rightarrow}&\mathcal{M}'\\
 \Updownarrow\scriptstyle{ }&&\Updownarrow\scriptstyle{}\\
\mathcal{M}&\stackrel{\left(F,G\right)}{\rightarrow}&\mathcal{M}',
\end{array}\label{diagg}\end{equation}
where each   Equivalence is defined  using (\ref{edef}). 

Defining $\mbox{wt}\left\{w_{ij}\right\}=2$ and $\mbox{wt}\left\{z_{kl}\right\}=1$, for all $i,j,k=1,\dots,q$ and $l=1,\dots,N$, it follows that
 \bp \label{propo1} Up to compositions with suitable linear   automorphisms of  $\mathcal{BSD}$-Models from  (\ref{models}), we obtain  
 \begin{equation}\begin{pmatrix}G_{11}\left(Z,W\right) &G_{12}\left(Z,W\right) \\ G_{21}\left(Z,W\right)& G_{22}\left(Z,W\right)
\end{pmatrix}=\begin{pmatrix} W+\mbox{O}_{\mbox{wt}}
(2)  &\mbox{O}_{\mbox{wt}}
(2) \\ \mbox{O}_{\mbox{wt}}
(2) & \mbox{O}_{\mbox{wt}}
(2)
\end{pmatrix}\hspace{0.1 cm}\mbox{and}\hspace{0.1 cm}\begin{pmatrix} F_{1}\left(Z,W\right) \\ F_{2}\left(Z,W\right) 
\end{pmatrix}=\begin{pmatrix} Z+\mbox{O}_{\mbox{wt}}
(2) & \mbox{O}_{\mbox{wt}}
(2)\\ \mbox{O}_{\mbox{wt}}
(2) & \mbox{O}_{\mbox{wt}}
(2) 
\end{pmatrix}. \label{lili1}
\end{equation} 
\ep
\begin{proof}   We write  the entries of the first matrix from (\ref{v11})  as follows
\begin{equation}\begin{split}& G_{11}\left(W,Z\right)=A\otimes W+\mbox{O}(2),\quad G_{12}\left(W,Z\right)=B\otimes W+\mbox{O}(2) ,\\& G_{21}\left(W,Z\right)=C\otimes W+\mbox{O}(2)  , \quad G_{22}\left(W,Z\right)=D\otimes W+\mbox{O}(2),\end{split}\label{99911}\end{equation}
where we have used by (\ref{edef}) the following matrices 
\begin{equation} \begin{split}& A=\left(A^{ij}\right)_{1\leq i,j\leq q}=\left(a^{ij}_{kl}\right)_{1\leq k,l\leq q}^{1\leq i,j\leq q},\quad\quad\hspace{0.24 cm}  B=   \left(B^{ij}\right)_{1 \leq i\leq q\atop{1 \leq j \leq q'-q}}=\left(b^{ij}_{kl}\right)_{1\leq k,l\leq q}^{1 \leq i\leq q\atop{1 \leq j \leq q'-q}}, \\&  C=\left(C^{ij}\right)_{1 \leq j \leq q\atop{ 1 \leq i \leq q'-q}}=\left(c^{ij}_{kl}\right)_{1\leq k,l\leq q}^{1 \leq j \leq q\atop{ 1 \leq i \leq q'-q}},  \quad  D=\left(   D^{ij}\right)_{1\leq i,j \leq q'-q}=\left(d^{ij}_{kl}\right)_{1\leq k,l\leq q}^{1\leq i,j \leq q'-q}.
\end{split} \label{bibi3}\end{equation}

  Let $F_{1}^{\left(1\right)}\left(Z\right)$ and  $F_{2}^{\left(1\right)}\left(Z\right)$ be the  linear parts in $Z$ of the matrices  $F_{1}\left(Z,W\right)$ and $F_{2}\left(Z,W\right)$ from (\ref{v11}). We replace (\ref{99911}) in    (\ref{v}) in order to identify  the terms of bidegree $2$ in $\left(Z,\overline{Z}\right)$.  We obtain

\begin{equation} \begin{split}&\frac{ A\otimes W-\overline{A\otimes W}^{t}}{2\sqrt{-1}}= \left<F_{1}^{\left(1\right)}\left(Z\right),F_{1}^{\left(1\right)}\left(Z\right)\right>,\quad \hspace{0.08 cm} \frac{B\otimes W-\overline{C\otimes W}^{t}}{2\sqrt{-1}}= \left<F_{1}^{\left(1\right)}\left(Z\right),F_{2}^{\left(1\right)}\left(Z\right)\right>,\\& \frac{C\otimes W-\overline{B\otimes W}^{t}}{2\sqrt{-1}}= \left<F_{2}^{\left(1\right)}\left(Z\right),F_{1}^{\left(1\right)}\left(Z\right)\right>,\quad  \frac{D\otimes W-\overline{D\otimes W}^{t}}{2\sqrt{-1}}= \left<F_{2}^{\left(1\right)}\left(Z\right),F_{2}^{\left(1\right)}\left(Z\right)\right>.\end{split} \label{8}
\end{equation}
 
We rewrite  the diagonal entries separately from the non-diagonal entries in  (\ref{models}). We obtain
\begin{equation}
 \begin{split}& \frac{w_{kl}-\overline{w_{lk}}}{2\sqrt{-1}}=\left<Z_{k},Z_{l}\right> ,\quad\hspace{0.03 cm}\mbox{for all $k\neq l$ and $k,l=1,\dots,q$,}\\&\quad\hspace{0.2 cm}\mbox{Im} w_{kk}=\left<Z_{k},Z_{k}\right> ,\quad \mbox{for all $k=1,\dots,q$.}\end{split} \label{yu1} 
\end{equation}

It remains to use   (\ref{bibi3}) and (\ref{yu1}) in order to  study   (\ref{8}). We obtain
\begin{equation}  \begin{split}& b^{ij}_{ll}\left(\mbox{Re}  w_{ll}+\sqrt{-1}\left<Z_{l},Z_{l}  \right>\right)-\overline{c^{ji}_{ll}}\left(\mbox{Re} w_{ll}-\sqrt{-1}\left<Z_{l},Z_{l}  \right>\right)=T_{ijll}\left(Z,\overline{Z}\right),\quad\hspace{0.06 cm} \mbox{for all   corresponding $i,j$,} \\&\quad\quad\quad\quad\quad\quad\quad\quad\hspace{0.02 cm}  b^{ij}_{kl}\left( \overline{w_{lk}}+2\sqrt{-1}\left<Z_{k},Z_{l}  \right>\right)-\overline{c^{ji}_{lk}}\left(  \overline{w_{lk}}\right)=T_{ijkl}\left(Z,\overline{Z}\right),\quad \mbox{for all     corresponding $i,j$ and $k\neq l$.} \end{split} \label{lak}
 \end{equation} 
 
Their right-hand sides  depend  only on $Z$ and $\overline{Z}$. We obtain  
$$\mbox{$\left(b^{ij}_{kl}\right)_{1 \leq k,l \leq q}=\left(\overline{c^{ji}_{lk}}\right)_{1 \leq k,l \leq q}$,}\quad\mbox{for all    corresponding $i,j$.}$$

More generally, we obtain
\begin{equation*} \begin{split}& A^{ij}=\overline{A^{ji}}^{t}, \quad\mbox{for all $i,j=1\dots,q$,} \quad  \quad\quad\hspace{0.1 cm} B^{ij}=\overline{C^{ji}}^{t}, \quad\mbox{for all corresponding $i,j$,}  \\& C^{ij}=\overline{B^{ji}}^{t}, \quad\mbox{for all corresponding $i,j$,} \quad  D^{ij}=\overline{D^{ji}}^{t},\quad\mbox{for all  $i,j=1\dots,q'-q$.} \end{split} 
 \end{equation*} 

 Reformulating (\ref{edef}) for the $\mathcal{BSD}$-Model $\mathcal{M}'$ from (\ref{models}),  we consider the Jacobian-matrix of   $G\left(0,W\right)$. It is denoted by
\begin{equation} \begin{split}& A'=\left({a'}^{ij}_{kl}\right)_{1\leq k,l\leq q'}^{1\leq i,j\leq q'}   \in\mathcal{M}_{{q'}^{2}\times {q'}^{2}}\left(\mathbb{C}\right),\quad\mbox{where}\hspace{0.1 cm}
 \left({a'}^{ij}_{kl}\right)_{1\leq k,l\leq q'}^{1\leq i,j\leq q'} = \left(\overline{{a'}^{ji}_{lk}}\right)_{1\leq k,l\leq q'}^{1\leq i,j\leq q'} .  \end{split} 
\label{6BB}
\end{equation} 
 
Given $k',l'\in 1,\dots, q'$ such that $k'<l'$, we consider  the transformation    $  Z'_{k'}\mapsto  Z'_{l'}$, $Z'_{l'}\mapsto Z'_{k'}$, $Z'_{k}\mapsto Z'_{k}$, for all  $k\in\left\{1,\dots,q'\right\}-\left\{k',l'\right\}$.  Therefore $  w'_{mk'}\mapsto w'_{ml'}$,    $w'_{k'm}\mapsto w'_{l'm}$,    $w'_{ml'}\mapsto w'_{mk'}$,    $w'_{l'm}\mapsto w'_{k'm}$,  for all $m,n=k',\dots,l'$. Other   entries of the matrix $W'$ remain unchanged.  In particular, we consider the number
\begin{equation*}q_{1}:=\displaystyle\max_{I =\{1,2,\dots,q_{1}\}\atop \det\left(\frac{g_{ij}}{w_{\tilde{i}\tilde{j}}}\right)_{i,j\in I\atop \tilde{i},\tilde{j}\in I}(0)\neq 0 } \mbox{Card}I. 
\end{equation*}

  It order to prove that $q=q_{1}$, we consider the matrix
\begin{equation} \begin{split}& A_{1}=\left({a}^{ij}_{kl}\right)_{1\leq k,l\leq q_{1}}^{1\leq i,j\leq q_{1}}   \in\mathcal{M}_{{q_{1}}^{2}\times {q_{1}}^{2}}\left(\mathbb{C}\right),\quad\mbox{where}\hspace{0.1 cm}
\left({a}^{ij}_{kl}\right)_{1\leq k,l\leq q_{1}}^{1\leq i,j\leq q_{1}} =\left(\overline{{a}^{ji}_{lk}}\right)_{1\leq k,l\leq q_{1}}^{1\leq i,j\leq q_{1}}. \end{split} 
\label{6VVa}
\end{equation}

We write the linear part of the $F$-component of the Embedding (\ref{v11}) as
  \begin{equation}V'\equiv \begin{pmatrix} V'_{11}  & V'_{12}  & \dots &    V'_{1q'}  \\ V'_{21}  & V'_{22}  & \dots &    V'_{2q'}  \\ \vdots & \vdots & \ddots & \vdots  \\ V'_{q'1}  & V'_{q'2}  & \dots &   V'_{q'q'} \end{pmatrix}\in\mathcal{M}_{q'N'\times q'N'}\left(\mathbb{C}\right). \label{9V1}\end{equation}

 We replace (\ref{bibi3}) in (\ref{8}) in order to extract the terms of degree $2$ in $\left(Z,\overline{Z}\right)$.  We obtain
\begin{equation} \begin{split}&{a'}_{11}^{11}\left<Z_{1}, Z_{1}\right>
 +  {a'}_{12}^{11}\left<Z_{1}, Z_{2}\right>
+ \overline{{a'}_{12}^{11}}\left<Z_{2}, Z_{1}\right>
+{a'}_{22}^{11}\left<Z_{2}, Z_{2}\right>+\dots
= \left(\displaystyle\sum_{k=1 }^{q'}\left(V'_{1k}\cdot Z_{k}^{t}\right)^{t}  \right)\cdot\overline{ \left(\displaystyle\sum_{l=1 }^{q'} V'_{1l}\cdot Z_{l}^{t} \right)}   ,\\& {a'}_{11}^{12}\left<Z_{1}, Z_{1}\right>
+ {a'}_{12}^{12}\left<Z_{1}, Z_{2}\right>
+ \overline{{a'}_{12}^{12}}\left<Z_{2}, Z_{1}\right>
+{a'}_{22}^{12}\left<Z_{2}, Z_{2}\right>+\dots
=\left(\displaystyle\sum_{k=1 }^{q'}\left(V'_{1k}\cdot Z_{k}^{t}\right)^{t}  \right)\cdot\overline{ \left(\displaystyle\sum_{l=1 }^{q'} V'_{2l}\cdot Z_{l}^{t} \right)} ,\\& {a'}_{11}^{22}\left<Z_{1}, Z_{1}\right>
+ {a'}_{12}^{22}\left<Z_{1}, Z_{2}\right>
+\overline{{a'}_{12}^{22}}\left<Z_{2}, Z_{1}\right>
+{a'}_{22}^{22}\left<Z_{2}, Z_{2}\right>+\dots
= \left(\displaystyle\sum_{k=1}^{q'}\left({V'}_{2k}\cdot Z_{k}^{t} \right)^{t}\right)\cdot
\overline{ \left(\displaystyle\sum_{l=1}^{q'} V'_{2l}\cdot Z_{l}^{t} \right)}   ,\\& \quad\quad\quad\quad\quad\quad \vdots\quad\quad\quad\quad\quad\quad\quad\quad\quad\quad\vdots\quad\quad\quad\quad\quad\quad\quad\quad\quad\quad  \vdots \quad\quad\quad\quad\quad\quad\quad\quad \quad\quad \vdots\quad \quad\quad\quad\quad\quad\quad    \quad\quad   \end{split} \label{777V}
\end{equation}
because the  $\left(i,j\right)$-entry of   matrix  from the left-hand side from (\ref{8})  is defined by
\begin{equation*}  \frac{ \displaystyle\sum_{k,l=1}^{q'}{a'}_{kl}^{ij}{w}_{kl}-\overline{\displaystyle\sum_{k,l=1}^{q'}{a'}_{kl}^{ij}{w}_{kl} }}{2\sqrt{-1}}  =   \displaystyle\sum_{k,l=1}^{q'}{a'}_{kl}^{ij} \left<Z_{k}, Z_{l}\right>,\quad\mbox{for all $i,j=1,\dots,q'$.}
\end{equation*}

We   collect     terms  in $\left(Z,\overline{Z}\right)$ from (\ref{777V}). We obtain
 \begin{equation}{a'}_{kl}^{ij}\left<Z_{k}, Z_{l}\right>=\left(\left(  V'_{ik}\cdot Z_{k}^{t}\right)^{t}  \right)\cdot\overline{ \left(  V'_{jl}\cdot Z_{l}^{t} \right)} , \quad  \mbox{for all $i,j,k,l=1,\dots,q'$.} \label{auxiliarV}
  \end{equation}
   
In particular for $i=j=k=l=1$ in (\ref{auxiliarV}), we obtain   
\begin{equation*}{a'}_{11}^{11}\cdot I_{N}=V'_{11}\cdot \overline{V'_{11}}^{t}.
\end{equation*}   
   
The first column   of the matrix $A_{1}$  vanishes in the light of   (\ref{777V}) if   $\det \left(V'_{11}\right)=0$ . Contradiction. In consequence, we obtain
\begin{equation*}\alpha_{11}=\sqrt{{a'}_{11}^{11}}\neq 0\hspace{0.15 cm}\mbox{and}\hspace{0.15 cm}V'_{11}=\alpha_{11}\cdot U_{11},\quad\mbox{where $U_{11}\in \mathcal{M}_{N\times N}\left(\mathbb{C}\right)$ is a unitary matrix.}
\end{equation*} 

In particular for $i=j=k=1$ and $l=2$ in (\ref{auxiliarV}), we obtain   
\begin{equation*}{a'}_{12}^{11}\cdot I_{N}=V'_{11}\cdot \overline{V'_{12}}^{t},\quad\mbox{or equivalently}\hspace{0.1 cm}{a'}_{12}^{11}\cdot I_{N}=\alpha_{11}\cdot U_{11}\cdot \overline{V'_{12}}^{t}.
\end{equation*}   
   
The second column  of the matrix $A_{1}$  vanishes in the light of (\ref{777V}) if   $\det \left(V'_{12}\right)=0$. Contradiction. In consequence, we obtain
\begin{equation*}\alpha_{12}=\overline{\frac{{a'}_{12}^{11}}{\sqrt{{a'}_{11}^{11}}}}\neq 0\hspace{0.15 cm}\mbox{and}\hspace{0.15 cm}V'_{12}=\alpha_{12}\cdot U_{11}^{-1}.
\end{equation*}

More generally, it exists a matrix $\left(\alpha_{ij}\right)_{1\leq i,j\leq q_{1}}$ such that $$ A_{1}= \left({a}^{ij}_{kl}\right)_{1\leq k,l\leq q_{1}}^{1\leq i,j\leq q_{1}} =\left(\alpha_{ik}\overline{\alpha}_{jl}\right)_{1\leq k,l\leq q_{1}}^{1\leq i,j\leq q_{1}}.$$  

The matrix $\left(\alpha_{ij}\right)_{1\leq i,j\leq q_{1}}$ must be invertible. 
Contrary, we can assume that    there exist $\beta_{1},\beta_{2},\dots,\beta_{q_{1}}\in\mathbb{C}$, not all vanishing, such that
\begin{equation*}\begin{split}&\hspace{0.11 cm} \alpha_{11}=\beta_{2}\alpha_{12}+\dots+\beta_{q_{1}}\alpha_{1q_{1}},\\& \hspace{0.11 cm} \alpha_{21}=\beta_{2}\alpha_{22}+\dots+\beta_{q_{1}}\alpha_{2q_{1}},\\&\quad\quad\quad\hspace{0.05 cm}  \quad\quad\quad\quad \ddots \\&\alpha_{q_{1}1}=\beta_{2}\alpha_{q_{1}2}+\dots+\beta_{q_{1}}\alpha_{q_{1}q_{1}}.\end{split} \end{equation*}
 
Contradiction, because the matrix $A_{1}$ is invertible. Otherwise, it exists a linear non-trivial combination of  columns for the matrix
\begin{equation*}  \begin{pmatrix}  {\begin{pmatrix} \alpha_{11}\overline{\alpha}_{11} & \alpha_{11}\overline{\alpha}_{12} & \dots & \alpha_{11}\overline{\alpha}_{1q_{1}}\\    \alpha_{11}\overline{\alpha}_{21} & \alpha_{11}\overline{\alpha}_{22} & \dots & \alpha_{11}\overline{\alpha}_{2q_{1}} \\  \vdots & \vdots & \ddots & \vdots \\  \alpha_{11}\overline{\alpha}_{q_{1}1} & \alpha_{11}\overline{\alpha}_{q_{1}2} & \dots & \alpha_{11}\overline{\alpha}_{q_{1}q_{1}} \end{pmatrix}} &  {\begin{pmatrix} \alpha_{12}\overline{\alpha}_{11} & \alpha_{12}\overline{\alpha}_{12} & \dots & \alpha_{12}\overline{\alpha}_{1q_{1}}\\    \alpha_{12}\overline{\alpha}_{21} & \alpha_{12}\overline{\alpha}_{22} & \dots & \alpha_{12}\overline{\alpha}_{2q_{1}} \\  \vdots & \vdots & \ddots & \vdots \\  \alpha_{12}\overline{\alpha}_{q_{1}1} & \alpha_{12}\overline{\alpha}_{q_{1}2} & \dots & \alpha_{12}\overline{\alpha}_{q_{1}q_{1}} \end{pmatrix}} & \begin{pmatrix}
\dots \\ \dots \\ \vdots \\ \dots
\end{pmatrix}  & {\begin{pmatrix} \alpha_{1q_{1}}\overline{\alpha}_{11}& \alpha_{1q_{1}}\overline{\alpha}_{12} & \dots & \alpha_{1q_{1}}\overline{\alpha}_{1q_{1}}\\    \alpha_{1q_{1}}\overline{\alpha}_{21} & \alpha_{1q_{1}}\overline{\alpha}_{22} & \dots & \alpha_{1q_{1}}\overline{\alpha}_{2q_{1}} \\  \vdots & \vdots & \ddots & \vdots \\  \alpha_{1q_{1}}\overline{\alpha}_{q_{1}1} & \alpha_{1q_{1}}\overline{\alpha}_{q_{1}2} & \dots &\alpha_{1q_{1}}\overline{\alpha}_{q_{1}q_{1}}\end{pmatrix}} 
\\ {\begin{pmatrix}  \alpha_{21}\overline{\alpha}_{11} & \alpha_{21}\overline{\alpha}_{12} & \dots & \alpha_{21}\overline{\alpha}_{1N}\\    \alpha_{21}\overline{\alpha}_{21} & \alpha_{21}\overline{\alpha}_{22} & \dots & \alpha_{21}\overline{\alpha}_{2q_{1}} \\  \vdots & \vdots & \vdots & \vdots \\  \alpha_{21}\overline{\alpha}_{q_{1}1} 
& \alpha_{21}\overline{\alpha}_{q_{1}2}& \dots & \alpha_{21}\overline{\alpha}_{q_{1}q_{1}}\end{pmatrix}} &  {\begin{pmatrix} \alpha_{22}\overline{\alpha}_{11} & \alpha_{22}\overline{\alpha}_{12} & \dots & \alpha_{22}\overline{\alpha}_{1q_{1}}\\    \alpha_{22}\overline{\alpha}_{21}& \alpha_{22}\overline{\alpha}_{22} & \dots &
 \alpha_{22}\overline{\alpha}_{2q_{1}} \\  \vdots & \vdots & \ddots & \vdots \\  \alpha_{22}\overline{\alpha}_{q_{1}1} & \alpha_{22}\overline{\alpha}_{q_{1}2} & \dots & \alpha_{22}\overline{\alpha}_{q_{1}q_{1}} \end{pmatrix}} & \begin{pmatrix}
\dots \\ \dots \\ \vdots \\ \dots
\end{pmatrix}  & {\begin{pmatrix} \alpha_{2q_{1}}\overline{\alpha}_{11} & \alpha_{2q_{1}}\overline{\alpha}_{12} & \dots & 
\alpha_{2q_{1}}\overline{\alpha}_{1q_{1}}\\    \alpha_{2q_{1}}\overline{\alpha}_{21} & \alpha_{2q_{1}}\overline{\alpha}_{22} & \dots & \alpha_{2q_{1}}\overline{\alpha}_{2q_{1}} \\  \vdots & \vdots & \ddots & \vdots \\  \alpha_{2q_{1}}\overline{\alpha}_{q_{1}1} & \alpha_{2q_{1}}\overline{\alpha}_{q_{1}2} & \dots & \alpha_{2q_{1}}\overline{\alpha}_{q_{1}q_{1}} \end{pmatrix}}\\ {\begin{pmatrix} \vdots & \vdots & \ddots & \vdots
\end{pmatrix}} &{\begin{pmatrix} \vdots & \vdots & \ddots & \vdots
\end{pmatrix}}& \begin{pmatrix}\vdots\end{pmatrix} & {\begin{pmatrix} \vdots & \vdots & \ddots & \vdots
\end{pmatrix}}    \\ {\begin{pmatrix} \alpha_{q_{1}1}\overline{\alpha}_{11}
& \alpha_{q_{1}1}\overline{\alpha}_{12}
 & \dots & \alpha_{q_{1}1}\overline{\alpha}_{1q_{1}}
\\    \alpha_{q_{1}1}\overline{\alpha}_{21} & \alpha_{q_{1}1}\overline{\alpha}_{22} & \dots & \alpha_{q_{1}1}\overline{\alpha}_{2q_{1}} \\   \vdots  & \vdots & \ddots & \vdots \\  \alpha_{q_{1}1}\overline{\alpha}_{q_{1}1} & \alpha_{q_{1}1}\overline{\alpha}_{q_{1}2} & \dots & \alpha_{q_{1}1}\overline{\alpha}_{q_{1}q_{1}} \end{pmatrix}} &  {\begin{pmatrix} \alpha_{q_{1}2}\overline{\alpha}_{11}
 & \alpha_{q_{1}2}\overline{\alpha}_{12}
 & \dots & \alpha_{q_{1}2}\overline{\alpha}_{1q_{1}}
\\    \alpha_{q_{1}2}\overline{\alpha}_{21} & \alpha_{q_{1}2}\overline{\alpha}_{22} & \dots & \alpha_{q_{1}2}\overline{\alpha}_{2q_{1}} \\  \vdots & \vdots & \ddots & \vdots \\  \alpha_{q_{1}2}\overline{\alpha}_{q_{1}1} & \alpha_{q_{1}k}\overline{\alpha}_{N2} & \dots & \alpha_{q_{1}k}\overline{\alpha}_{q_{1}}\end{pmatrix}} & \begin{pmatrix}
\dots \\ \dots \\ \vdots \\ \dots
\end{pmatrix}  & {\begin{pmatrix} \alpha_{q_{1}q_{1}}\overline{\alpha}_{11}
 & \alpha_{q_{1}q_{1}}\overline{\alpha}_{12}
 & \dots & \alpha_{q_{1}q_{1}}\overline{\alpha}_{1q_{1}}
\\    \alpha_{q_{1}q_{1}}\overline{\alpha}_{21} & \alpha_{q_{1}q_{1}}\overline{\alpha}_{22} & \dots & \alpha_{q_{1}q_{1}}\overline{\alpha}_{2q_{1}} \\  \vdots & \vdots & \ddots & \vdots \\  \alpha_{q_{1}q_{1}}\overline{\alpha}_{q_{1}1} & \alpha_{q_{1}q_{1}}\overline{\alpha}_{q_{1}2} & \dots & \alpha_{q_{1}q_{1}}\overline{\alpha}_{q_{1}q_{1}}\end{pmatrix}}   \end{pmatrix}.
\end{equation*}

 It suffices to take a unitary matrix $\tilde{V}\in\mathcal{M}_{N\times N}\left(\mathbb{C}\right)$  in order to consider the matrix 
\begin{equation*}\tilde{V}'\equiv \begin{pmatrix} \alpha_{11}\tilde{V}  &  \alpha_{12}\tilde{V} & \dots &   \alpha_{1q_{1}} \tilde{V}  \\ \alpha_{21}\tilde{V}   & \alpha_{22}\tilde{V}  & \dots &    \alpha_{2q_{1}}\tilde{V}   \\ \vdots & \vdots & \ddots & \vdots  \\ \alpha_{q_{1}1}\tilde{V}   & \alpha_{q_{1}2}\tilde{V}   & \dots &   \alpha_{q_{1}q_{1}}\tilde{V}  \end{pmatrix}\in\mathcal{M}_{q_{1}N\times q_{1}N}\left(\mathbb{C}\right). \end{equation*}
 
The  matrix $\tilde{V}'$ must be   invertible, because the matrix $\tilde{V}'$ is the Kronecker Product defined by two invertible matrices. The change of coordinates  $\left(\tilde{W}, \tilde{Z}\right)=\left(A_{1}\otimes W,\tilde{V}'\otimes Z\right)$ is further used in order to study the equations (\ref{8}), because
\begin{equation*}  \frac{A_{1}\otimes W-\left(\overline{A_{1}\otimes W}\right)^{t}}{2\sqrt{-1}}= \left(\tilde{V}' \otimes Z\right)\cdot\left(\overline{ \tilde{V}' \otimes Z}\right)^{t}. \end{equation*}

In particular, we obtain
\begin{equation} \begin{split}&\quad\quad\quad\hspace{0.16 cm}\frac{\tilde{W}-\overline{\tilde{W}}^{t}}{2\sqrt{-1}}= \left<F_{1}^{\left(1\right)}\left(\tilde{V}'^{-1}\otimes \tilde{Z}\right),F_{1}^{\left(1\right)}\left(\tilde{V}'^{-1}\otimes \tilde{Z}\right)\right>,\quad  \frac{\tilde{B}\otimes \left(\tilde{W}-\overline{  \tilde{W}}^{t}\right)}{2\sqrt{-1}}= \left<F_{1}^{\left(1\right)}\left(\tilde{V}'^{-1}\otimes \tilde{Z}\right),F_{2}^{\left(1\right)}\left(\tilde{V}'^{-1}\otimes \tilde{Z}\right)\right>, \\&  \frac{\tilde{C}\otimes \left(\tilde{W}-\overline{  \tilde{W}}^{t}\right)}{2\sqrt{-1}}= \left<F_{2}^{\left(1\right)}\left(\tilde{V}'^{-1}\otimes \tilde{Z}\right),F_{1}^{\left(1\right)}\left(\tilde{V}'^{-1}\otimes \tilde{Z}\right)\right>,\quad  \frac{\tilde{D}\otimes\left( \tilde{W}-\overline{\tilde{W}}^{t}\right)}{2\sqrt{-1}}= \left<F_{2}^{\left(1\right)}\left(\tilde{V}'^{-1}\otimes \tilde{Z}\right),
F_{2}^{\left(1\right)}\left(\tilde{V}'^{-1}\otimes \tilde{Z}\right)\right>,\end{split} 
\label{550q}\end{equation}
where we have used by (\ref{bibi3}) the following matrices
\begin{equation}
\begin{split}&  
\quad\quad \quad\hspace{0.1 cm} \mbox{$\tilde{B}=B\diamond A_{1}^{-1}\in \mathcal{M}_{q_{1}\left(q'-q_{1}\right) \times  q_{1}^{2}}\left(\mathbb{C}\right)$ of  entries  $\left<\left(c^{ij}_{kl}\right)_{1\leq k,l\leq q_{1}\atop{1\leq k,l\leq q}}^{1 \leq i \leq q'-q_{1}\atop{1 \leq j\leq q_{1}}},\left(\overline{A}_{1}^{-1}\right)^{t}\right>$,} \\&\quad\quad\quad\hspace{0.1 cm} \mbox{$\tilde{C}= C\diamond A_{1}^{-1}\in \mathcal{M}_{q_{1}^{2}\times q_{1}\left(q'-q_{1}\right)   }\left(\mathbb{C}\right)$ of entries $\left<\left(b^{ij}_{kl}\right)_{1\leq k,l\leq q_{1}\atop{1\leq k,l\leq q_{1}}}^{1 \leq i\leq  q_{1}\atop{1 \leq  j\leq q'-q_{1}}},\left(\overline{A}_{1}^{-1}\right)^{t}\right>$,}\\& \mbox{$ \tilde{D}=D\diamond A_{1}^{-1}\in\mathcal{M}_{ q_{1}\left(q'-q_{1}\right) \times q_{1}\left(q'-q_{1}\right) }\left(\mathbb{C}\right)$ of entries  $\left<\left(d^{ij}_{kl}\right)_{1\leq k,l\leq q_{1}}^{1 \leq i,j\leq q'-q_{1}},\left(\overline{A}_{1}^{-1}\right)^{t}\right>$, }\end{split}\label{chil}
\end{equation}
because $\diamond$ defines an obvious rule of   multiplication of matrices, becase the matrices $\left(\overline{A}_{1}^{-1}\right)^{t}$, $B$, $C$ and $D$  may be seen as   vectors according to (\ref{nota}). In particular, the matrices $\tilde{B},\tilde{C}$ and $\tilde{D}$ may be assumed  null according to following  change of coordinates: given  $\tilde{Z}_{1},\tilde{Z}_{2}, \dots,\tilde{Z}_{q_{1}}$   the rows of the matrix $\tilde{Z}$ and $ R_{1}\left(\tilde{Z}\right) ,R_{2}\left(\tilde{Z}\right) ,\dots,R_{q_{1}}\left(\tilde{Z}\right)$    the  rows   of the matrix $F_{1}^{(1)}\left(\tilde{V}'^{-1}\otimes \tilde{Z}\right)$, we write
  $R_{k}\left(\tilde{Z}\right)=\displaystyle\sum_{i=1}^{q_{1}} R_{k}\left(\tilde{Z}_{i}\right)$, for all $k=1,\dots,q_{1}$. In particular, we focus on all terms depending on $\left(Z,\tilde{Z}\right)$ in (\ref{550q}). We obtain
$$\left( \left<\tilde{Z}_{i} ,\tilde{Z}_{j}\right>\right)_{1\leq i,j \leq q_{1}}=\left(
\displaystyle\sum_{k,l=1}^{q_{1}}\left<R_{k}\left(\tilde{Z}_{j}\right),\displaystyle R_{l}\left(\tilde{Z}_{j}\right)\right>\right)_{1\leq i,j \leq q_{1}}.$$

Therefore
 $R_{j}\left(\tilde{Z}\right)=R_{j}\left(\tilde{Z}_{j}\right)$ for all $j=1,\dots,q_{1}$, and
 $$\left(\left< R_{i}\left(\tilde{Z}_{i}\right),\displaystyle R_{j}\left(\tilde{Z}_{j}\right)\right>\right)_{1\leq i,j \leq q_{1}}=\left(\delta_{i}^{j} \left<\tilde{Z}_{i} ,\tilde{Z}_{j}\right>\right)_{1\leq i,j \leq q_{1}}.$$

Given $\mathcal{L}_{1}\left(\tilde{Z}\right),\mathcal{L}_{2}\left(\tilde{Z}\right),\dots,\mathcal{L}_{q'-q_{1}}\left(\tilde{Z}\right)$   the row vectors of the matrix $F_{2}^{\left(1\right)}\left(\tilde{V}'^{-1}\otimes \tilde{Z}\right)$, we write    $\mathcal{L}_{k}\left(\tilde{Z}\right)=\displaystyle\sum_{i=1}^{q_{1}} \mathcal{L}_{k}^{(i)}\left(\tilde{Z}_{i}\right)$, for all $k=1,\dots,q'-q_{1}$.  Returning to (\ref{550q}), we obtain
\begin{equation*}\left(\displaystyle\sum_{m,n=1}^{q_{1}}\tilde{b}^{jk}_{mn}\left< \tilde{Z}_{m} ,    \tilde{Z}_{n} \right>\right)_{1\leq j \leq q_{1}\atop{1 \leq k \leq q'-q_{1}}}=\left(\left<R_{j}\left(\tilde{Z}_{j}\right),  \displaystyle\sum_{i=1}^{q_{1}} \mathcal{L}_{k}^{(i)}\left(\tilde{Z}_{i}\right)\right>\right)_{1\leq j \leq q_{1}\atop{1 \leq k \leq q'-q_{1}}}. 
\end{equation*}

  In consequence, we obtain  \begin{equation}\begin{split}&  \tilde{b}^{jk}_{mn} =0 ,\quad\quad\quad\quad\quad\quad\quad\quad\quad\quad \quad\quad\quad\quad  \mbox{for all $m,n,j=1,\dots,q_{1}$ with $m\neq j$ and $k=1,\dots, q'-q_{1}$,}\\&\tilde{b}^{jk}_{ji}\left< \tilde{Z}_{j} ,    \tilde{Z}_{i} \right>=\left<R_{j}\left(\tilde{Z}_{j}\right),   \mathcal{L}_{k}^{(i)}\left(\tilde{Z}_{i}\right)\right>,\quad  \mbox{for all $i,j=1,\dots,q_{1}$  and $k=1,\dots, q'-q_{1}$.} \end{split}\label{rcH}\end{equation}

In particular  for $k=i=1$  and $j=1$ in (\ref{rcH}), we obtain $\tilde{b}^{11}_{11}\left< \tilde{Z}_{1} ,    \tilde{Z}_{1} \right>=\left<R_{1}\left(\tilde{Z}_{1}\right),   \mathcal{L}_{1}^{(1)}\left(\tilde{Z}_{1}\right)\right>$, therefore
\begin{equation*} \left<\tilde{b}^{11}_{11}\cdot \tilde{Z}_{1} ,    \tilde{Z}_{1} \right>=\left<\left(\mathcal{L}_{1}^{(1)}\right)^{\star}\left(   R_{1}\left(\tilde{Z}_{1}\right)\right),    \tilde{Z}_{1} \right> \hspace{0.1 cm}\mbox{and}\hspace{0.1 cm}\mathcal{L}_{1}^{(1)}=\left(\overline{\tilde{b}^{11}_{11}\cdot R_{1}^{-1}}\right)^{t}.  
 \end{equation*} 
 
  In particular  for $k=i=1$  and $j=2$ in (\ref{rcH}), we obtain $\tilde{b}^{11}_{12}\left< \tilde{Z}_{1} ,    \tilde{Z}_{1} \right>=\left<R_{1}\left(\tilde{Z}_{1}\right),   \mathcal{L}_{2}^{(1)}\left(\tilde{Z}_{1}\right)\right>$, therefore
\begin{equation*} \left<\tilde{b}^{11}_{12}\cdot \tilde{Z}_{1} ,    \tilde{Z}_{1} \right>=\left<\left(\mathcal{L}_{2}^{(1)}\right)^{\star}\left(   R_{1}\left(\tilde{Z}_{1}\right)\right),    \tilde{Z}_{1} \right> \hspace{0.1 cm}\mbox{and}\hspace{0.1 cm}\mathcal{L}_{2}^{(1)}=\left(\overline{\tilde{b}^{11}_{12}\cdot R_{1}^{-1}}\right)^{t}.  
 \end{equation*}
 
$$\vdots\quad \vdots\quad\quad\quad\quad \vdots\quad\quad\quad\quad \vdots\quad\quad\quad\quad \vdots\quad\quad\quad\quad\vdots\quad\quad\quad\quad\vdots\quad\quad\quad\quad \vdots $$ 
 
 In particular  for $k=i=1$  and $j=q_{1}$ in (\ref{rcH}), we obtain $\tilde{b}^{11}_{1q_{1}}\left< \tilde{Z}_{1} ,    \tilde{Z}_{1} \right>=\left<R_{1}\left(\tilde{Z}_{1}\right),   \mathcal{L}_{q_{1}}^{(1)}\left(\tilde{Z}_{1}\right)\right>$, therefore
\begin{equation*} \left<\tilde{b}^{11}_{1q_{1}}\cdot \tilde{Z}_{1} ,    \tilde{Z}_{1} \right>=\left<\left(\mathcal{L}_{q_{1}}^{(1)}\right)^{\star}\left(   R_{1}\left(\tilde{Z}_{1}\right)\right),    \tilde{Z}_{1} \right> \hspace{0.1 cm}\mbox{and}\hspace{0.1 cm}\mathcal{L}_{q_{1}}^{(1)}=\left(\overline{\tilde{b}^{11}_{1q_{1}}\cdot R_{1}^{-1}}\right)^{t}.  
 \end{equation*}

The  computations can be similarly repeated for all  $k=1,\dots, q'-q_{1}$. In particular, it suffices to consider the matrix
\begin{equation*} \mathcal{V}= \begin{pmatrix}\left(\overline{\tilde{b}^{11}_{11}\cdot R_{1}^{-1}}\right)^{t} 
   &  \left(\overline{\tilde{b}^{11}_{12}\cdot R_{1}^{-1}}\right)^{t}  &  \dots &   \left(\overline{\tilde{b}^{11}_{1\hspace{0.05 cm}q_{1}}\cdot R_{1}^{-1}}\right)^{t}          \\  \left(\overline{\tilde{b}^{12}_{11}\cdot R_{1}^{-1}}\right)^{t} 
   &  \left(\overline{\tilde{b}^{12}_{12}\cdot R_{1}^{-1}}\right)^{t}  &  \dots &   \left(\overline{\tilde{b}^{12}_{1\hspace{0.05 cm}q_{1}}\cdot R_{1}^{-1}}\right)^{t}     \\   \vdots
  &  \vdots & \ddots & \vdots  \\ \left(\overline{\tilde{b}^{1\hspace{0.05 cm}q'-q_{1}}_{11}\cdot R_{1}^{-1}}\right)^{t} 
   &  \left(\overline{\tilde{b}^{1\hspace{0.05 cm}q'-q_{1}}_{12}\cdot R_{1}^{-1}}\right)^{t}  &  \dots &   \left(\overline{\tilde{b}^{1\hspace{0.05 cm}q'-q_{1}}_{1\hspace{0.05 cm}q_{1}}\cdot R_{1}^{-1}}\right)^{t} 
\end{pmatrix}.
\end{equation*}

In order to leave invariant the equations (\ref{550q}), we implement the coordinates  
\begin{equation*}  \begin{pmatrix}
\tilde{Z}_{1} \\ \tilde{Z}_{2} 
\end{pmatrix}:=\begin{pmatrix}
\tilde{Z}_{1} \\ \tilde{Z}_{2}- \mathcal{V}\otimes \tilde{Z}_{1} 
\end{pmatrix},  \end{equation*}

\begin{equation*}  \begin{pmatrix}\left(\tilde{w}_{ij}\right)_{1 \leq i\leq  q_{1}\atop{1 \leq  j\leq q_{1}}} & \left(\tilde{w}_{ij}\right)_{1 \leq i\leq q'- q_{1}\atop{1 \leq  j\leq q_{1}}}  \\ \left(\tilde{w}_{ij}\right)_{1 \leq j\leq  q'-q_{1}\atop{1 \leq  i\leq q_{1}}} &\left(\tilde{w}_{ij}\right)_{1 \leq i\leq q'- q_{1}\atop{1 \leq  j\leq q'- q_{1}}}
\end{pmatrix}:=\begin{pmatrix}\left(\tilde{w}_{ij}\right)_{1 \leq i\leq  q_{1}\atop{1 \leq  j\leq q_{1}}} & \left(\tilde{w}_{ij}\right)_{1 \leq i\leq q'- q_{1}\atop{1 \leq  j\leq q_{1}}}- \tilde{B} \otimes\tilde{W} \\ \left(\tilde{w}_{ij}\right)_{1 \leq i\leq  q'-q_{1}\atop{1 \leq  j\leq q_{1}}}- \tilde{C} \otimes\tilde{W} & \left(\tilde{w}_{ij}\right)_{1 \leq i\leq q'- q_{1}\atop{1 \leq  j\leq q'- q_{1}}}-\tilde{D}_{1}  \otimes\tilde{W}-\overline{\tilde{D}}^{t}_{1}  \otimes\tilde{W}    +\tilde{D}_{2}  \otimes\tilde{W}
\end{pmatrix},  \end{equation*}
where we have used the following matrices
 \begin{equation*} \begin{split}& \mathcal{D}_{1}= \begin{pmatrix}
\begin{pmatrix}
b^{11}_{11} & b^{11}_{12} & \dots & b^{11}_{1q_{1}} \\ b^{11}_{11} & b^{11}_{12} & \dots & b^{11}_{1q_{1}} \\ \vdots & \vdots & \ddots & \vdots  \\ b^{11}_{11} & b^{11}_{12} & \dots & b^{11}_{1q_{1}}   
\end{pmatrix} & \begin{pmatrix}
0 & 0 & \dots & 0 \\ 0 & 0 & \dots & 0 \\ \vdots & \vdots & \ddots & \vdots  \\ 0 & 0 & \dots & 0 \\ 
\end{pmatrix} & \begin{pmatrix}
  \dots  \\   \dots   \\ \vdots    \\   \dots   \\ 
\end{pmatrix} & \begin{pmatrix}
0 & 0 & \dots & 0 \\ 0 & 0 & \dots & 0 \\ \vdots & \vdots & \ddots & \vdots  \\ 0 & 0 & \dots & 0 \\ 
\end{pmatrix} \\ \begin{pmatrix}
0 & 0 & \dots & 0 \\ 0 & 0 & \dots & 0 \\ \vdots & \vdots & \ddots & \vdots  \\ 0 & 0 & \dots & 0 \\ 
\end{pmatrix} & \begin{pmatrix}
b^{12}_{11} & b^{12}_{12} & \dots & b^{12}_{1q_{1}} \\ b^{12}_{11} & b^{12}_{12} & \dots & b^{12}_{1q_{1}} \\ \vdots & \vdots & \ddots & \vdots  \\ b^{12}_{11} & b^{12}_{12} & \dots & b^{12}_{1q_{1}}   
\end{pmatrix} & \begin{pmatrix}
  \dots  \\   \dots   \\ \vdots    \\   \dots   \\ 
\end{pmatrix} & \begin{pmatrix}
0 & 0 & \dots & 0 \\ 0 & 0 & \dots & 0 \\ \vdots & \vdots & \ddots & \vdots  \\ 0 & 0 & \dots & 0 \\ 
\end{pmatrix} \\ \begin{pmatrix}
\vdots & \vdots & \ddots & \vdots   
\end{pmatrix} & \begin{pmatrix}
\vdots & \vdots & \ddots & \vdots   
\end{pmatrix} & \ddots & \begin{pmatrix}
\vdots & \vdots & \ddots & \vdots   
\end{pmatrix}  \\ \begin{pmatrix}
0 & 0 & \dots & 0 \\ 0 & 0 & \dots & 0 \\ \vdots & \vdots & \ddots & \vdots  \\ 0 & 0 & \dots & 0 \\ 
\end{pmatrix} & \begin{pmatrix}
0 & 0 & \dots & 0 \\ 0 & 0 & \dots & 0 \\ \vdots & \vdots & \ddots & \vdots  \\ 0 & 0 & \dots & 0 \\ 
\end{pmatrix} & \begin{pmatrix}
  \dots  \\   \dots   \\ \vdots    \\   \dots   \\ 
\end{pmatrix} & \begin{pmatrix}
b^{1\hspace{0.05 cm}q'-q_{1}}_{11} & b^{1\hspace{0.05 cm}q'-q_{1}}_{12} & \dots & b^{1\hspace{0.05 cm}q'-q_{1}}_{1q_{1}} \\ b^{1\hspace{0.05 cm}q'-q_{1}}_{11} & b^{1\hspace{0.05 cm}q'-q_{1}}_{12} & \dots & b^{1\hspace{0.05 cm}q'-q_{1}}_{1q_{1}} \\ \vdots & \vdots & \ddots & \vdots  \\ b^{1\hspace{0.05 cm}q'-q_{1}}_{11} & b^{1\hspace{0.05 cm}q'-q_{1}}_{12} & \dots & b^{1\hspace{0.05 cm}q'-q_{1}}_{1q_{1}}   
\end{pmatrix}   
\end{pmatrix},\\& \mathcal{D}_{2}=  \left(\overline{\tilde{b}^{1k}_{1i}}\cdot  \tilde{b}^{1k'}_{1i'}     \right) _{ 1 \leq k,k'\leq  q'-q_{1}\atop{1 \leq i,i'\leq   q_{1}}}.\end{split}
 \end{equation*}

Finally, we focus on the matrix
 \begin{equation*}\begin{pmatrix} R_{1}
& \mbox{O}_{N\times N} & \dots & \mbox{O}_{N\times N} \\ \mbox{O}_{N\times N}  & R_{2} & \dots &  \mbox{O}_{N\times N} \\  \vdots & \vdots & \ddots & \vdots \\ \mbox{O}_{N\times N} & \mbox{O}_{N\times N} & \dots &  R_{q_{1}}
\end{pmatrix}.
\end{equation*} 
    
In particular, we consider the sets of  rows  \begin{equation}\begin{split}& \mbox{$\left\{\alpha_{1}\left(i\right),\alpha_{2}\left(i\right),\dots,\alpha_{N}\left(i\right)\right\}_{i=1,\dots,q_{1}}$ in $\mathbb{C}^{ q_{1} N}$ such that:}  \\&\quad\quad\quad\hspace{0.2 cm} \left(\left< \alpha_{u}\left(i\right), \alpha_{l}\left(j\right)\right>\right)_{1 \leq u,l \leq  N}^{1 \leq i,j \leq q_{1}}=\left(\delta_{u}^{l}\cdot\delta_{i}^{j}\right)_{1 \leq u,l \leq  N}^{1 \leq i,j \leq q_{1}} .\end{split} \label{for}\end{equation}  

The procedure of Baouendi-Huang\cite{BH} applied for (\ref{for})  provides the   sets of orthonormal rows   
$$\mbox{$\left\{\alpha_{1}\left(i\right), \dots,\alpha_{N}\left(i\right),\alpha_{N+1}^{\star}\left(i\right),\dots,\alpha_{2 N}^{\star}\left(i\right)\right\}_{i=1,\dots,q_{1}}$   using a set of matrices $\left\{\tilde{\mathcal{A}}_{i}\right\}_{i=1,\dots,q_{1}}$.}$$  

Defining $\tau^{\star} \left(\tilde{Z}\right)=\tilde{Z}^{\star}$, we use  the composition  $F_{1}:=\tau^{\star}  \circ F_{1}$,
given   the matrix $\tilde{Z}^{\star}$ of   rows   defined by
$$\mbox{$\left(\tilde{Z}^{\star}_{i}\right)_{1 \leq i \leq q_{1}}=\left({\tilde{Z}}_{i}\tilde{A}_{i}^{-1}\right)_{1 \leq i \leq q_{1}}$ and  $\left( \tilde{Z}^{\star}_{q+j}\right)_{1 \leq j \leq q'-q_{1}}=\left({\tilde{Z}}_{q+j}\right)_{1 \leq j \leq q'-q_{1}}$.}$$ 

   Therefore $q=q_{1}$, because   $\left(F,G\right)$ is a local Embedding.
\end{proof}

\section{Analogues of the   normalizations $(2.5)$ from Huang\cite{huang1}}

In order to understand   the structure of  the   Embedding (\ref{v11}), we define the product of matrices
\begin{equation} A\otimes W=\left(\displaystyle\sum_{k,l=1}^{q}a_{kl}^{ij}w_{kl}\right)_{1\leq i\leq q'\atop 1\leq j\leq N'},\quad\mbox{given the matrix}\hspace{0.1 cm}A=\left(a^{\alpha}_{\beta}\right)^{1\leq \alpha \leq q'N'}_{1\leq \beta \leq  qq'}\in\mathcal{M}_{qq'\times qN'}\left(\mathbb{C}\right).\label{edef000}
\end{equation}   

\subsection{Parameters} In particular, we write 
 \begin{equation} \begin{pmatrix}
 F_{1}\left(Z,W\right)\\ F_{2}\left(Z,W\right)
 \end{pmatrix} =\begin{pmatrix} Z \\ 0  \end{pmatrix}  + A\otimes W+ 
 \mbox{O}_{\mbox{wt}}\left(3\right). \label{90A}
\end{equation}

Moreover, the   Embedding (\ref{v11}) depends on  the following parameters
 \begin{equation}  r^{ij}_{ab}= \frac{1}{2}\left(\frac{\partial^{2}g_{ji}}{\partial w_{aa} \partial w_{bb}}(0) + \overline{\frac{\partial^{2}g_{ij}}{\partial w_{aa} \partial w_{bb}}(0)} \right) ,\quad\mbox{for all $a,b=1,\dots, q$ and $i,j=1,\dots, q'$.}\label{09}
 \end{equation}

The  matrix $R=\left(r^{ij}_{ab}\right)_{1 \leq  a,b \leq q}^{1 \leq  i,j \leq q'}$  induces by (\ref{GG2}),(\ref{GG222}) and (\ref{edef}) the matrix  $R\left(W\right)$ defined by 

   \begin{equation*}   \begin{pmatrix} \begin{pmatrix}r^{11}_{11}\cdot w_{11}        & r^{11}_{12}\cdot w_{12}        & \dots & r^{11}_{1q}\cdot w_{1q}   \\   r^{12}_{11}\cdot w_{11}        & r^{12}_{12}\cdot w_{12}        & \dots & r^{12}_{1q}\cdot w_{1q}   \\   \vdots & \vdots &   \ddots & \vdots \\ r^{1q}_{11}\cdot w_{11}             & r^{1q}_{12}\cdot w_{12}             &\dots & r^{1q}_{1q}\cdot w_{1q}  
    \end{pmatrix}  &  \begin{pmatrix}r^{11}_{21}\cdot w_{21}        & r^{11}_{21}\cdot w_{22}        & \dots & r^{11}_{2q}\cdot w_{2q}   \\   r^{12}_{21}\cdot w_{21}        & r^{12}_{21}\cdot w_{22}        & \dots & r^{12}_{2q}\cdot w_{2q}   \\ \vdots &  \vdots &  \ddots & \vdots \\ r^{1q}_{21}\cdot w_{21}             & r^{1q}_{22}\cdot w_{22}             &\dots & r^{1q}_{2q}\cdot w_{2q}  
    \end{pmatrix}  &  \begin{pmatrix}  \dots     \\  \dots     \\  \vdots   \\  \dots    
    \end{pmatrix}  &   \begin{pmatrix} r^{11}_{q1}\cdot w_{q1}       &  r^{11}_{q2}\cdot w_{q2}       &  \dots & r^{11}_{qq}\cdot w_{qq}       \\  r^{12}_{q1}\cdot w_{q1}       &  r^{12}_{q2}\cdot w_{q2}       &  \dots & r^{12}_{qq}\cdot w_{qq}       \\\vdots       & \vdots       &\ddots& \vdots \\ r^{1q}_{q1}\cdot w_{q1}          &  r^{1q}_{q2}\cdot w_{q2}          & \dots & r^{1q}_{qq}\cdot w_{qq}  
    \end{pmatrix}      \\   \begin{pmatrix}r^{21}_{11}\cdot w_{11}        & r^{21}_{12}\cdot w_{12}        & \dots & r^{21}_{1q}\cdot w_{1q}   \\  r^{22}_{11}\cdot w_{11}        & r^{22}_{12}\cdot w_{12}        & \dots & r^{22}_{1q}\cdot w_{1q}   \\   \vdots &  \vdots &  \ddots & \vdots \\ r^{2q}_{11}\cdot w_{11}             & r^{2q}_{12}\cdot w_{12}             &\dots & r^{2q}_{1q}\cdot w_{1q}  
    \end{pmatrix}  &  \begin{pmatrix}r^{21}_{21}\cdot w_{21}        & r^{21}_{22}\cdot w_{22}        & \dots & r^{21}_{2q}\cdot w_{2q}   \\  r^{22}_{21}\cdot w_{21}        & r^{22}_{22}\cdot w_{22}        & \dots & r^{22}_{2q}\cdot w_{2q}   \\   \vdots &  \vdots & \ddots & \vdots \\ r^{2q}_{21}\cdot w_{21}             & r^{2q}_{22}\cdot w_{22}             &\dots & r^{2q}_{2q}\cdot w_{2q}  
    \end{pmatrix}  &  \begin{pmatrix}  \dots     \\  \dots     \\   \vdots   \\  \dots    
    \end{pmatrix}  &   \begin{pmatrix} r^{21}_{q1}\cdot w_{q1}       & r^{21}_{q2}\cdot w_{q2}       &   \dots & r^{21}_{qq}\cdot w_{qq}       \\ r^{22}_{q1}\cdot w_{q1}       & r^{22}_{q2}\cdot w_{q2}       &   \dots & r^{22}_{qq}\cdot w_{qq}       \\\vdots       & \vdots       & \ddots& \vdots \\ r^{2q}_{q1}\cdot w_{q1} &  r^{2q}_{q2}\cdot w_{q2}          &   \dots &  r^{2q}_{qq}\cdot w_{qq}  
    \end{pmatrix}      \\ {\begin{pmatrix}   \vdots & \vdots &  \ddots & \vdots
\end{pmatrix}} &  {\begin{pmatrix}   \vdots & \vdots &  \ddots & \vdots
\end{pmatrix}} & \begin{pmatrix}\vdots\end{pmatrix} & {\begin{pmatrix}   \vdots & \vdots &  \ddots & \vdots
\end{pmatrix}} \\ \begin{pmatrix}r^{q'1}_{11}\cdot w_{11} &  r^{q'1}_{12}\cdot w_{12} & \dots & r^{q'1}_{1q}\cdot w_{1q}      \\   r^{q'2}_{11}\cdot w_{11} &  r^{q'2}_{12}\cdot w_{12} & \dots & r^{q'2}_{1q}\cdot w_{1q}      \\  \vdots & \vdots & \ddots& \vdots \\ r^{q'q'}_{11}\cdot w_{11}   &  r^{q'q'}_{12}\cdot w_{12}   &  \dots & r^{q'q'}_{1q} \cdot w_{1q}
    \end{pmatrix} &  \begin{pmatrix}r^{q'1}_{21}\cdot w_{21} & r^{q'1}_{22}\cdot w_{22} &  \dots & r^{q'1}_{2q}\cdot w_{2q}      \\  r^{q'2}_{21}\cdot w_{21} & r^{q'2}_{22}\cdot w_{22} &  \dots & r^{q'2}_{2q}\cdot w_{2q}      \\  \vdots & \vdots &  \ddots& \vdots \\ r^{q'q'}_{21}\cdot w_{21}   &  r^{q'q'}_{22}\cdot w_{22}   &  \dots & r^{q'q'}_{2q} \cdot w_{2q}
    \end{pmatrix} & \begin{pmatrix}
\dots   \\ \dots   \\\vdots \\ \dots
\end{pmatrix} & \begin{pmatrix} r^{q'1}_{q'1}\cdot w_{q1}  & r^{q'1}_{q'2}\cdot w_{q2}  & \dots & r^{q'1}_{q'q}\cdot w_{qq}  \\  r^{q'2}_{q'1}\cdot w_{q1}  & r^{q'2}_{q'2}\cdot w_{q2}  & \dots & r^{q'2}_{q'q}\cdot w_{qq}  \\   \vdots &  \vdots &\ddots& \vdots \\ r^{q'q'}_{q'1}\cdot w_{q1}     & r^{q'q'}_{q'2}\cdot w_{q2}     & \dots & r^{q'q'}_{q'q} \cdot w_{qq}
    \end{pmatrix}  
    \end{pmatrix}.  
\end{equation*}

In order to eliminate the parameters $A$ and $R$,   we follow   Baouendi-Huang\cite{BH} and Chern-Moser\cite{CM}:
\subsection{Elimination of the  matrix R}  We search   for a    transformation   as
  \begin{equation}T_{1}\left(W',Z'\right)=\left(\frac{1}{I_{q'^{2}}+R\left(W\right)} \otimes W', V\left(W',Z'\right)\right),\quad\mbox{for a fit matrix $V\left(W',Z'\right)$.}\label{be45}
    \end{equation}
    
The transformation (\ref{be45})  preserves  the   $\mathcal{BSD}$-Model $\mbox{Im}W'=Z'\overline{Z'}^{t}$:
\begin{equation}\frac{ \frac{1}{I_{q'^{2}}+R\left(W\right)} \otimes W'-\left(\overline{\frac{1}{ I_{q'^{2}}+R\left(W\right)} \otimes W'}\right)^{t}}{2\sqrt{-1}} = \left(V\left(W',Z'\right) \right)\cdot\overline{\left(V\left(W',Z'\right) \right)}^{t}. \label{LU}
 \end{equation}  

We make fit computations in (\ref{LU}). Given the model $\mbox{Im}W'=Z'\overline{Z'}^{t}$, we obtain
\begin{equation*} \tilde{A} \otimes \frac{  \left(\overline{I_{q'^{2}}+R\left(W\right)  }\right)^{t} \otimes W'-\left(\overline{  \left(\overline{I_{q'^{2}}+R\left(W\right)  }\right)^{t}\otimes W'}\right)^{t} }{2\sqrt{-1}}  =\left(V\left(W',Z'\right) \right)\cdot\overline{\left(V\left(W',Z'\right) \right)}^{t}, 
  \end{equation*}
  because the matrix $\tilde{A}\left(W',\overline{W}'\right)$ is hermitian:
  \begin{equation}
\tilde{A}\left(W',\overline{W}'\right)=\frac{1}{\left(I_{q'^{2}}+R\left(W\right)  \right)\cdot \left(\overline{I_{q'^{2}}+R\left(W\right) }\right)^{t}}.\label{bv}
\end{equation}

In the light of (\ref{edef000}) and  (\ref{09}), we have  
\begin{equation*} 
 \left(\overline{I_{q'^{2}}+R\left(W\right)  }\right)^{t} \otimes W'-\left(\overline{  \left(\overline{I_{q'^{2}}+R \left(W\right)}\right)^{t}\otimes W'}\right)^{t}  =W'-\overline{W'}^{t}+\overline{R\left(W\right)}^{t}\otimes W'-\overline{\left(\overline{R \left(W\right)  }^{t}\otimes W'\right)}^{t}  =W'-\overline{W'}^{t}, 
\end{equation*}
according to  the following computations
\begin{equation*}\begin{split} \overline{\left(\overline{R \left(W\right)  }^{t}\otimes W'\right)}^{t}&=\left(\overline{\displaystyle\sum_{k,l=1}^{q}\overline{r^{ji}_{lk} w_{lk}}w_{kl}}\right)_{1 \leq i,j\leq {q'}^{2}}^{t}= \left(\displaystyle\sum_{k,l=1}^{q} r^{ji}_{lk} \overline{w_{kl}}w_{lk}\right)^{t}_{1 \leq i,j\leq {q'}^{2}}\\&       =\left(\displaystyle\sum_{k,l=1}^{q} \overline{r^{ij}_{kl} w_{kl}}w_{lk}\right)^{t}_{1 \leq i,j\leq {q'}^{2}}=\left(\displaystyle\sum_{k,l=1}^{q} \overline{r^{ij}_{lk} w_{lk}}w_{kl}\right)_{1 \leq i,j\leq {q'}^{2}}=\overline{R\left(W\right)}^{t}\otimes W'.\end{split}
\end{equation*}

Given the model $\mbox{Im}W'=Z'\overline{Z'}^{t}$, the  equation (\ref{LU}) is equivalent to  
\begin{equation} \tilde{A}\left(W',\overline{W}'\right) \otimes \left< Z', Z'\right>=\left(V\left(W' ,Z'\right)\right)\cdot\overline{\left(V\left(W',Z'\right)\right)}^{t}. \label{LU1}
 \end{equation}  
 
In particular, we consider the  formal expansion
\begin{equation*}\frac{1}{ I_{q'^{2}}+R\left(W\right)}= \left( r_{\alpha'}^{\beta'}\left(W'\right)\right)_{1\leq \alpha',\beta' \leq {q'}^{2}}. 
\end{equation*}  
     
 In consequence, it suffices to  consider
 
   \begin{equation}V\left(W',Z'\right)=  \begin{pmatrix}\displaystyle\sum_{k=1}^{q'}r_{k1}^{11}\left(W'\right)\cdot  Z'_{k}    & \displaystyle\sum_{k=1}^{q'}r_{k2}^{11}\left(W'\right)\cdot  Z'_{k}   & \dots & \displaystyle\sum_{k=1}^{q'}r_{kq'}^{11}\left(W'\right)\cdot  Z'_{k}      \\   \vdots & \vdots & \ddots & \vdots \\ \displaystyle\sum_{k=1}^{q'}r_{k1}^{1q'}\left(W'\right)\cdot  Z'_{k}    & \displaystyle\sum_{k=1}^{q'}r_{k2}^{1q'}\left(W'\right)\cdot  Z'_{k}   & \dots & \displaystyle\sum_{k=1}^{q'}r_{kq'}^{1q'}\left(W'\right)\cdot  Z'_{k}\\
  \vdots & \vdots & \ddots & \vdots \\ \displaystyle\sum_{k=1}^{q'}r_{k1}^{q'1}\left(W'\right)\cdot  Z'_{k}    & \displaystyle\sum_{k=1}^{q'}r_{k2}^{q'1}\left(W'\right)\cdot  Z'_{k}   & \dots & \displaystyle\sum_{k=1}^{q'}r_{kq'}^{q'1}\left(W'\right)\cdot  Z'_{k}      \\   
  \vdots & \vdots & \ddots & \vdots \\   \displaystyle\sum_{k=1}^{q'}r_{k1}^{q'q'}\left(W'\right)\cdot  Z'_{k}    & \displaystyle\sum_{k=1}^{q'}r_{k2}^{q'q'}\left(W'\right)\cdot  Z'_{k}   & \dots & \displaystyle\sum_{k=1}^{q'}r_{k'q'}^{q'q'}\left(W'\right)\cdot  Z'_{k}   
    \end{pmatrix}.\label{yp}
   \end{equation}

We compute
  \begin{equation}\begin{split}& \partial_{w_{aa}}\left(\frac{1}{I_{q'^{2}}+R\left(G\left(W,Z\right)\right)}\cdot G\left(W,Z\right)\right)=\frac{1}{I_{q'^{2}}+R\left(G\left(W,Z\right)\right) }\cdot \partial_{w_{aa}}\left( G\left(W,Z\right)\right) +\\& \quad\quad\quad\quad\quad \quad\quad\quad\quad\quad\quad\quad\quad\quad\quad\quad\quad\quad\quad \partial_{w_{aa}}\left(\frac{1}{I_{q'^{2}}+R\left(G\left(W,Z\right)\right) }\right)\cdot G\left(W,Z\right),\quad\mbox{for all $a=1,\dots,q$.}\end{split}\label{GU1}
  \end{equation}

According to the differentiation rule of a product of matrices applied in
 the left-hand side in (\ref{GU1}), we compute
 \begin{equation}\begin{split}&\hspace{0.26 cm}
\partial_{w_{bb}}\left(\frac{1}{I_{q'^{2}}+R\left(G\left(W,Z\right)\right)}\cdot\partial_{w_{aa}}\left( G\left(W,Z\right)\right)\right)=\partial_{w_{bb}}\left(\frac{1}{I_{q'^{2}}+R\left(G\left(W,Z\right)\right)}\right)\cdot\partial_{w_{aa}}\left( G\left(W,Z\right)\right)+  \\&\quad\quad\quad\quad\quad\quad\quad\quad \quad\quad\quad\quad\quad\quad \quad\quad\quad\quad \quad\quad\quad\quad\quad\quad \quad\quad\quad\quad\quad\quad \hspace{0.1 cm}\quad\quad \quad\quad \frac{\partial^{2}_{w_{bb} w_{aa}}\left( G\left(W,Z\right) \right)}{I_{q'^{2}}+R\left(G\left(W,Z\right)\right)},\quad \mbox{for all $a,b=1,\dots,q$,} \\& \partial_{w_{bb}}\left(\partial_{w_{aa}}\left(\frac{1}{I_{q'^{2}}+R\left(G\left(W,Z\right)\right)}\right)\cdot G\left(W,Z\right)\right)=   \partial^{2}_{w_{bb} w_{aa}}\left(\frac{1}{I_{q'^{2}}+R\left(G\left(W,Z\right)\right)}\right) \cdot G\left(W,Z\right) +  \\&\quad\quad\quad\quad\quad \quad\quad \quad\quad \quad\quad \quad\quad \quad\quad \quad\quad  \quad\quad \quad\quad   \partial_{w_{aa}}\left(\frac{1}{I_{q'^{2}}+R\left(G\left(W,Z\right)\right) }\right)\cdot \partial_{w_{bb}}\left( G\left(W,Z\right)\right), \quad\quad \hspace{0.2 cm} \mbox{for all $a,b=1,\dots,q$.} \end{split}\label{GU2}
  \end{equation} 

Finally, we consider $\left(G^{\star},F^{\star}\right)=T
_{1}\circ \left(G,F\right) $ using  (\ref{be45}).   We obtain 
\begin{equation}\left(  \left.\left(\overline{\frac{\partial^{2} \left(g^{\star}_{ij}\left(Z,W\right) \right)}{\partial w_{aa}\partial w_{bb}   }} +   \frac{\partial^{2} \left(g^{\star}_{ji}\left(Z,W\right) \right)}{\partial w_{aa}\partial w_{bb}   }\right)\right\vert_{Z= O_{q\times N} \atop{W= O_{q\times q}}}\right)_{1 \leq i,j \leq q'} =\mbox{O}_{q' \times q'}, \quad\mbox{ for all $a,b= 1,\dots,q$,} \label{99se}
\end{equation}
according to (\ref{GU2})   and to the following evaluations
  \begin{equation*}\begin{split}&  \hspace{0.2 cm}\quad\quad\quad \left. \frac{\partial^{2}_{w_{bb} w_{aa}}\left( G\left(W,Z\right) \right)}{I_{q'^{2}}+R\left(G\left(W,Z\right)\right)} \right\vert_{Z= O_{q\times N} \atop{W= O_{q\times q}}}=\left. \partial^{2}_{w_{bb} w_{aa}}\left( G\left(W,Z\right) \right) \right\vert_{Z= O_{q\times N} \atop{W= O_{q\times q}}},\quad\quad\quad \hspace{0.05 cm}  \mbox{for all $a,b=1,\dots,q$,}\\& \left. \partial_{w_{bb}}\left(\frac{1}{I_{q'^{2}}+R\left(G\left(W,Z\right)\right)}\right)\right\vert_{ W= O_{q\times q}}= -\left.\partial_{w_{bb}}\left(R\left(G\left(W,Z\right)\right)\right)\right\vert_{ W= O_{q\times q}},\quad\quad \hspace{0.13 cm} \mbox{for all $b=1,\dots,q$,} \end{split} \end{equation*}
  \begin{equation*}\left.\partial^{2}_{w_{bb} w_{aa}}\left(\frac{1}{I_{q'^{2}}+R\left(G\left(W,Z\right)\right)}\right) \cdot G\left(W,Z\right)\right\vert_{Z= O_{q\times N} \atop{W= O_{q\times q}}}=  O_{q\times q} ,\quad \mbox{for all $a,b=1,\dots,q$.}
  \end{equation*}
\begin{equation*}\left.\partial_{w_{aa}}\left( G\left(W,Z\right)\right)\right\vert_{Z= O_{q\times N} \atop{W= O_{q\times q}}}=\begin{pmatrix}
 0& \dots & 0& \dots & 0 \\ \vdots & \ddots & \vdots& \ddots & \vdots \\ 0& \dots & 1& \dots & 0 \\ \vdots & \ddots & \vdots & \ddots & \vdots \\ 0& \dots & 0 & \dots & 0
 \end{pmatrix}  \hspace{0.1 cm}\mbox{and}\hspace{0.1 cm}\left.\partial_{w_{bb}}\left( G\left(W,Z\right)\right)\right\vert_{Z= O_{q\times N} \atop{W= O_{q\times q}}}=\begin{pmatrix}
 0& \dots & 0& \dots & 0 \\ \vdots & \ddots & \vdots& \ddots & \vdots \\ 0& \dots & 1& \dots & 0 \\ \vdots & \ddots & \vdots & \ddots & \vdots \\ 0& \dots & 0 & \dots & 0
 \end{pmatrix},  
  \end{equation*}
where  the $\left(a,a\right)$-entry and the $\left(b,b\right)$-entry are $1$, otherwise the entries are vanishing,   for all $a,b=1,\dots,q$.
 
 \subsection{Elimination of  the Matrix $A$:} We search  for a   transformation as 
 
   \begin{equation} T_{2}\left(W',Z'\right)=\left(\frac{1}{I_{{q'}^{2}} +\mathcal{L}\left(W',Z'\right)}\otimes W',V\left(W',Z', Z'-A\otimes W' \right)\right),\quad\mbox{for   fit matrices      $ V\left(W',Z',Z'\right)$ and $\mathcal{L}\left(W',Z'\right)$.} \label{be4}    \end{equation}
 
 The matrix  $\mathcal{L}\left(W',Z'\right)$ is chosen such that
\begin{equation}  \mathcal{L}\left(W',Z'\right)=I_{q'}+2\sqrt{-1} C\left(Z'\right)-2\sqrt{-1}B\left(W'\right),\quad\mbox{for fit matrices  $C\left(Z'\right)$ and $B\left(W'\right)$.}\label{fut2}
\end{equation}
 
The transformation (\ref{be4})   preserves  the   $\mathcal{BSD}$-Model $\mbox{Im}W'=Z'\overline{Z'}^{t}$:
\begin{equation}\frac{  \frac{1}{I_{{q'}^{2}} +\mathcal{L}\left(W',Z'\right)}\otimes W'-\overline{\frac{1}{I_{{q'}^{2}} +\mathcal{L}\left(W',Z'\right)}\otimes W'}^{t}}{2\sqrt{-1}}   = \left(V\left(W',Z', Z'-A\otimes W \right)\right)\cdot  \overline{\left(V\left(W',Z' , Z'-A\otimes W \right)\right)}^{t}.  \label{lol1}
\end{equation}

We make fit computations in (\ref{lol1}). Given the model $\mbox{Im}W'=Z'\overline{Z'}^{t}$, we obtain
\begin{equation}  \tilde{A}\left(W',Z',\overline{W'},\overline{Z'}\right) \otimes\frac{  W'-\overline{W'}^{t}+\overline{\mathcal{L}\left(W',Z'\right)}^{t}\otimes W'- \mathcal{L}\left(W',Z'\right)\otimes \overline{W'}^{t} }{2\sqrt{-1}}=  \left(V\left(W',Z', Z'-A\otimes W \right)\right)\cdot  \hspace{0.1 cm} \overline{\left(V\left(W',Z', Z'-A\otimes W \right)\right)}^{t},    \label{lol4}
\end{equation}   because the matrix $ \tilde{A}\left(W',Z',\overline{W'},\overline{Z'}\right)$ is hermitian:
\begin{equation} \tilde{A}\left(W',Z',\overline{W'},\overline{Z'}\right)
=\frac{1}{\left(I_{{q'}^{2}} +2\sqrt{-1} C\left(Z'\right)-2\sqrt{-1}B\left(W'\right)\right)\cdot\overline{\left(I_{{q'}^{2}} +2\sqrt{-1} C\left(Z'\right)-2\sqrt{-1}B\left(W'\right)\right)}^{t}}.\label{bv1} 
\end{equation}

 Given the model $\mbox{Im}W'=Z'\overline{Z'}^{t}$, we search for a matrix $V\left(W',Z'\right)$ such that
  \begin{equation}  \tilde{A}\left(W',Z',\overline{W'},\overline{Z'}  \right)\otimes\frac{ W'- \overline{ W'} ^{t}}{2\sqrt{-1}}= \left(V\left(W',Z',Z'\right)\right)\cdot\left(\overline{V\left(W',Z', Z'\right)}\right)^{t}.\label{new1}
  \end{equation}

In particular, it suffices to  consider
 \begin{equation}V\left(W',Z',Z'\right)=  \begin{pmatrix}\displaystyle\sum_{k=1}^{q'}l_{k1}^{11}\left(W',Z'\right)\cdot  Z'_{k}    & \displaystyle\sum_{k=1}^{q'}l_{k2}^{11}\left(W',Z'\right)\cdot  Z'_{k}   & \dots & \displaystyle\sum_{k=1}^{q'}l_{kq'}^{11}\left(W',Z'\right)\cdot  Z'_{k}       \\     \vdots & \vdots & \ddots & \vdots \\ \displaystyle\sum_{k=1}^{q'}l_{k1}^{1q'}\left(W',Z'\right)\cdot  Z'_{k}    & \displaystyle\sum_{k=1}^{q'}l_{k2}^{1q'}\left(W',Z'\right)\cdot  Z'_{k}   & \dots & \displaystyle\sum_{k=1}^{q'}l_{kq'}^{1q'}\left(W',Z'\right)\cdot  Z'_{k}\\
  \vdots & \vdots & \ddots & \vdots \\ \displaystyle\sum_{k=1}^{q'}l_{k1}^{q'1}\left(W',Z'\right)\cdot  Z'_{k}    & \displaystyle\sum_{k=1}^{q'}l_{k2}^{q'1}\left(W',Z'\right)\cdot  Z'_{k}   & \dots & \displaystyle\sum_{k=1}^{q'}l_{kq'}^{q'1}\left(W',Z'\right)\cdot  Z'_{k}                   \\ 
  \vdots & \vdots & \ddots & \vdots \\   \displaystyle\sum_{k=1}^{q'}l_{k1}^{N'q'}\left(W',Z'\right)\cdot  Z'_{k}    & \displaystyle\sum_{k=1}^{q'}l_{k2}^{N'q'}\left(W',Z'\right)\cdot  Z'_{k}   & \dots & \displaystyle\sum_{k=1}^{q'}l_{kq'}^{N'q'}\left(W',Z'\right)\cdot  Z'_{k}   
    \end{pmatrix},\label{yp1}
   \end{equation}
where we have used the formal expansion
\begin{equation*}\frac{1}{ I_{q'^{2}}+\mathcal{L}\left(W',Z'\right)}= \left( l_{\alpha'}^{\beta'}\left(W',Z'\right)\right)_{1\leq \alpha' \leq {q'}^{2}\atop{1\leq \beta' \leq {q'}{N'}}}. 
\end{equation*}

Given the model $\mbox{Im}W'=Z'\overline{Z'}^{t}$,  we search for   the matrix     
$$ C\left(Z'\right)=\left(c_{kl}^{ij}\left(Z'\right)\right)^{1\leq i,j\leq q'^{2}}_{1\leq k,l\leq q'^{2}}\hspace{0.1 cm}\mbox{such that} \hspace{0.1 cm}C\left(Z'\right)\otimes \overline{ W'}^{t}  =   Z' \cdot\left(\overline{  A\otimes W'}\right)^{t}.
$$

In particular, we obtain
\begin{equation*} \begin{split}&\displaystyle\sum_{k,l=1}^{q'^{2}}c_{kl}^{ij}\left(Z'\right)\overline{w'_{lk}} =\displaystyle\sum_{u=1}^{N'}z'_{iu} \cdot \left(\overline{\displaystyle\sum_{k,l=1}^{q}a_{kl}^{uj}w'_{kl}}\right), \hspace{0.1 cm}\mbox{for all $i,j=1,\dots,q'$,}\quad\mbox{or equivalently:}\\&\quad\quad\quad\quad\left( c_{kl}^{ij}\left(Z'\right)\right)^{1\leq i,j\leq q'^{2}}_{1\leq k,l\leq q'^{2}}=\left(\displaystyle\sum_{u=1}^{N'}\overline{a_{lk}^{uj}}\cdot z'_{iu}\right)^{1\leq i,j\leq q'^{2}}_{1\leq k,l\leq q'^{2}}.  \end{split}
   \end{equation*}
  
In consequence,  we    have   
  $$ \tilde{A}\left(W',Z',\overline{W'},\overline{Z'}\right)\otimes\left(  C\left(Z'\right)\otimes \overline{ W'} ^{t}\right) = \left(V\left(W',Z', Z'\right)\right)\cdot\left(\overline{ V\left(W',Z',A\otimes W'\right)}\right)^{t} .$$

Given the model $\mbox{Im}W'=Z'\overline{Z'}^{t}$,  we search for  the matrix  
as  
$$ B\left(W'\right)=\left( \displaystyle\sum_{k',l'=1}^{q'^{2}}\left(b_{kl}^{ij}\right)_{k'l'}w'_{k'l'}  \right)^{1\leq i,j\leq q'^{2}}_{1\leq k,l\leq q'^{2}} \hspace{0.1 cm}\mbox{such that} \hspace{0.1 cm}  B\left(W'\right)  \otimes \overline{ W'}^{t} + \overline{\left(B\left(W'\right)\right)^{t} }  \otimes   W' =     \left(A\otimes W'\right)\cdot\left(\overline{  A\otimes W'}\right)^{t}.
 $$
 
In particular, we obtain
$$ \left(\displaystyle\sum_{k,l=1}^{q'^{2}}\displaystyle\sum_{k',l'=1}^{q'^{2}}\left(\left(b_{kl}^{ij}\right)_{k'l'}w'_{k'l'}\right)\overline{w'_{lk}}+\displaystyle\sum_{k,l=1}^{q'^{2}}\overline{\left(\displaystyle\sum_{k',l'=1}^{q'^{2}}\left(b_{lk}^{ji}\right)_{k'l'}w'_{k'l'}\right)} w'_{kl}\right)_{1\leq i,j\leq q'^{2}} =\displaystyle\sum_{r=1}^{N'}\left(\displaystyle\sum_{k,l=1}^{q^{2}}a_{kl}^{ir}w'_{kl}\right) \overline{\left(\displaystyle\sum_{k,l=1}^{q^{2}}a_{kl}^{rj}w'_{kl}\right)},$$
or equivalently, we obtain
   \begin{equation*} \left(b_{kl}^{ij}\right)_{k'l'}+\overline{\left(b_{l'k'}^{ji}\right)_{lk}} =\left\{\begin{split}&  \displaystyle\sum_{r=1}^{N'} \left( a_{k'l'}^{ir}   \cdot\overline{a_{kl}^{rj}}+ a_{kl}^{ir}   \cdot\overline{a_{k'l'}^{rj}}\right) ,\quad\hspace{0.2 cm}\mbox{for all $i,j,k',l'=1,\dots,q^{2}$ and $k',l'=1,\dots,q'^{2}$, }\\& \quad\quad \quad \quad\quad \quad \quad 0,\quad\quad\quad\quad \quad \quad \quad\mbox{for all $i,j,k,l=1,\dots,q'^{2}$ and  $k',l'=q^{2}+1,\dots,q'^{2}$.} \end{split}\right.
      \end{equation*}

In consequence,  we    have   
 $$   \tilde{A}\left(W',Z',\overline{W'},\overline{Z'}\right)\otimes\left(  B\left(W'\right)  \otimes \overline{ W'}^{t} + \overline{\left(B\left(W'\right)\right)^{t} }  \otimes   W' \right) =  \left(V\left(W',Z',  A\otimes W'\right)\right) \cdot\left(\overline{ V\left(W',Z',A\otimes W'\right)}\right)^{t}.$$      

Finally, we consider     $\left(G^{\star\star},F^{\star\star}\right)=T_{2}\circ \left(G ,F \right)$ using  (\ref{be4}). We obtain 
\begin{equation} \left(\left.\frac{\partial f^{\star\star}_{il}\left(Z,W\right)}{\partial w_{ab}}\right\vert_{Z= O_{q\times N} \atop{W= O_{q\times q}}}\right)_{1 \leq i \leq q'\atop{1 \leq l \leq N'}}=\mbox{O}_{q' \times N'} ,\quad\mbox{for all  $a,b= 1,\dots,q$.}\label{99te}\end{equation}

The  normalisations (\ref{99se}) and  (\ref{99te}) are essential in order to study the following:

\subsection{Rank} Given the $\mathcal{BSD}$-Models $\mathcal{M}$ and  $\mathcal{M}'$ defined by (\ref{models}), we  examine using (\ref{v}) the   following equations:
\begin{equation}
\begin{split}& \frac{G^{\star\star}_{11}\left(Z,W\right)-\overline{G^{\star\star}_{11}\left(Z,W\right)}}{2\sqrt{-1}}=F_{1}^{\star\star}\left(Z,W\right)\overline{ F_{1}^{\star\star}\left(Z,W\right) } ,\quad\quad \frac{G^{\star\star}_{12}\left(Z,W\right)-\overline{G^{\star\star}_{21}\left(Z,W\right)}}{2\sqrt{-1}}=F_{1}^{\star\star}\left(Z,W\right) \overline{ F_{2}^{\star\star}\left(Z,W\right)},\\& \frac{G^{\star\star}_{21}\left(Z,W\right)-\overline{G^{\star}_{12}\left(Z,W\right)}}{2\sqrt{-1}}=F_{2}^{\star\star}\left(Z,W\right) \overline{ F_{1}^{\star\star}\left(Z,W\right) },\quad\quad \frac{G^{\star\star}_{22}\left(Z,W\right)- \overline{G^{\star\star}_{22}\left(Z,W\right)}}{2\sqrt{-1}} =F_{2}^{\star\star}\left(Z,W\right) \overline{ F_{2}^{\star\star}\left(Z,W\right) } . \end{split}\label{qqq1}
\end{equation}

The equations (\ref{qqq1}) are    studied in order to compute the formal mapping (\ref{v11}) from 
 their   diagonal entries.   We extract  terms of total degree   $4$ from (\ref{qqq1}) in order to implement linear changes of coordinates using (\ref{edef}).  Following Huang-Ji\cite{HJ}, we use    the   notation
\begin{equation}F_{1}^{\star\star}\left(Z,W\right)=\left(f^{\star\star},\varphi^{\star\star}\right)\left(Z,W\right).
\label{extra}
\end{equation} 
 
We work  with formal power series in $\left(Z,W\right)$ defining: the weight of any entry of $Z$ is $1$, and the weight to any entry of $W$ is $2$.  Its   weighted degree, denoted by $n$,   is the minimum of the weighted degrees of   its  homogeneous terms from its formal expansion. 
We write    $H\left(Z,W\right)=\mbox{O}(n)$, if $H\left(Z,W\right)$ is  a formal power series of weighted degree $n$.  More precisely, it follows that
  
 \bp\label{propo11}There  exist  $k_{0}\in\{0,1,\dots,q\}$, a set of polynomials of degree $2$, denoted by $\left\{\mathcal{L}_{il}^{k}\left(Z_{i},Z_{k}\right)\right\}_{1\leq i \leq k_{0}\atop{1 \leq k  \leq q \atop{1 \leq l \leq N}}}$ and  a set    of polynomials of degree $1$, denoted by  $\left\{a_{ku}^{il}\left(Z_{i}\right)\right\}_{1\leq i \leq k_{0}\atop{1 \leq k  \leq q \atop{1 \leq l \leq N}}}$,     such that 
\begin{equation*}\begin{split}&\quad\left( f_{il}^{\star\star }\left(Z,W\right)\right)_{ 1\leq i\leq q\atop{1\leq l\leq N}} =Z  + \begin{pmatrix}
\displaystyle\sum_{k,u=1\atop{(k,u)\in\mathcal{S}}}^{q}\left(a_{ku}^{il}\left(Z_{i}\right)w_{ku}\right)_{ 1\leq i\leq k_{0}\atop{1\leq l\leq N}} \\ \mbox{O}_{\left(q-k_{0}\right)\times N}\end{pmatrix}+\mbox{O}_{\mbox{wt}}(4),  \\& \quad\left( \varphi^{\star\star }_{il}\left(Z,W\right)\right)_{1\leq i\leq q\atop{1\leq l\leq N}}=\begin{pmatrix}
\left(  \displaystyle\sum_{k=1}^{q}\mathcal{L}_{il}^{k}\left(Z_{i},Z_{k}\right) \right)_{ 1\leq i\leq k_{0}\atop{1\leq l\leq N}}   
\\ \mbox{O}_{\left(q-k_{0}\right)\times N}
\end{pmatrix}    +\mbox{O}_{\mbox{wt}}(3),\end{split}
   \end{equation*}      
where we have used the set    
     \begin{equation*}\mathcal{S}:=\left\{(k,i)\in \left\{1,\dots,k_{0}\right\}\times \left\{1,\dots,q\right\} \left\vert \left(a_{kk}^{il}\left(Z_{i}\right)\right)_{1 \leq l  \leq N }  \not\equiv (0,0,\dots,0)   \right.\right\}.
\end{equation*}
 \ep 
\begin{proof} We   write (\ref{extra})  by (\ref{nota}),(\ref{10v})  and   (\ref{lili1})        as

\begin{equation}\left. \begin{split}&\left( g^{\star \star}_{ij}\left(Z,W\right)\right)_{1\leq i,j  \leq q}= \left(w_{ij}+A_{ij}\left(Z\right)+\displaystyle\sum_{k,u=1}^{q}b^{ij}_{ku}\left(Z\right)w_{ku}+  \displaystyle\sum_{k,u=1\atop{k',u'=1}}^{q}D^{ij}_{kuk'u'}w_{ku}w_{k'u'}+\mbox{O}_{\mbox{wt}}(5)\right)_{1\leq i,j  \leq q},    \\&      \left( f^{\star\star}_{il}\left(Z,W\right)\right)_{1\leq i  \leq q\atop{1\leq l  \leq N}} =\left(z_{il}+\displaystyle\sum_{k,u=1}^{q}a_{ku}^{il}\left(Z\right)w_{ku}+ b_{il}\left(Z\right)+\mbox{O}_{\mbox{wt}}(4)\right)_{1\leq i  \leq q\atop{1\leq l  \leq N}},\\&  \left(\varphi_{il}^{\star\star }\left(Z,W\right)\right)_{1\leq i  \leq q\atop{1\leq l  \leq N}} = \left(\left(\varphi_{il}^{\star\star }\left(Z\right)\right)^{(2)} +\mbox{O}_{\mbox{wt}}(3)\right)_{1\leq i  \leq q\atop{1\leq l  \leq N}}  , \end{split}\right.\label{982}
\end{equation}
where we have considered  using (\ref{nota})  homogeneous polynomials  of    degree $1$ and  of degree $2$ in $Z$,  denoted by 
$$\mbox{ $\left( b_{il}\left(Z\right)\right)_{1\leq i  \leq q\atop{1\leq l  \leq N}}, \left(A_{ij}\left(Z\right)\right)_{1\leq i,j  \leq q}, \left(a_{ku}^{il}\left(Z\right)\right)^{1 \leq  k,u \leq q}_{1\leq i  \leq q\atop{1\leq l  \leq N}} , \left( b^{ij}_{ku}\left(Z\right)\right)^{1 \leq  k,u \leq q}_{1\leq i,j  \leq q}$ and $ \left( \left(\varphi_{il}^{\star\star }\left(Z\right)\right)^{(2)}\right)_{1\leq i  \leq q\atop{1\leq l  \leq N}}$. }$$

 We study the sum of   terms of degree at most $4$ in $\left(Z,\overline{Z}\right)$ from the  entries of (\ref{qqq1}) using  (\ref{982}). We obtain
\begin{equation}\begin{split}& \mbox{Im} \left\{A_{ij}\left(Z\right)+\displaystyle\sum_{k,u=1}^{q}b_{ku}^{ij}\left(Z\right)w_{ku}+  \displaystyle\sum_{k,u=1\atop{k',u'=1}}^{q}D^{ij}_{kuk'u'}w_{ku}w_{k'u'}  \right\}= 2\mbox{Re} \left\{ \displaystyle\sum_{l=1}^{N}\overline{z}_{il}\cdot   \left( b_{jl}\left(Z\right)+\right.\right.\\&\quad\quad\quad\quad\quad\quad\quad\quad\quad\quad\quad\quad\quad \quad\quad\hspace{0.2 cm} \left.\left.\displaystyle\sum_{k,u=1}^{q}a_{ku}^{jl}\left(Z\right)w_{ku}\right)\right\}+  \displaystyle\sum_{l=1}^{N} \left(\varphi_{il}^{\star\star }\left(Z\right)\right)^{(2)}\overline{\left(\varphi_{jl}^{\star\star }\left(Z\right)\right)^{(2)}},\quad\quad\mbox{for all $i,j=1,\dots,q$.} \end{split} \label{val1999se} \end{equation}

It suffices to identify linear terms in (\ref{val1999se}), in order to obtain     $\left( A_{ij}\left(Z\right)\right)_{1 \leq i,j \leq q}=\mbox{O}_{q\times q}$.   The equation (\ref{val1999se}) is equivalent to
\begin{equation} S_{1}-S_{2}+S_{3}+S_{4}+S_{5}-S_{6}-S_{7}      =  \displaystyle\sum_{l=1}^{N}\left(\varphi_{il}^{\star\star }\left(Z\right)\right)^{(2)}\overline{\left(\varphi_{jl}^{\star\star }\left(Z\right)\right)^{(2)}}+2\mbox{Re} \left\{ \displaystyle\sum_{l=1}^{N}\overline{z}_{il}b_{jl}\left(Z\right) \right\},\quad\mbox{for all $i,j=1,\dots,q$,}\label{val1999} 
\end{equation} 
where we have used the following sums
\begin{equation*}S_{1}=\frac{1}{2\sqrt{-1}} \left( \displaystyle\sum_{k,u=1}^{q}\left(b_{kk}^{ij}\left(Z\right)\left(\Re w_{kk}+\sqrt{-1}\left<Z_{k},Z_{k}\right>\right)+  D^{ij}_{kkuu}\left(\Re w_{kk}+\sqrt{-1}\left<Z_{k},Z_{k}\right>\right)\left(\Re w_{uu}+\sqrt{-1}\left<Z_{u},Z_{u}\right>\right)\right)\right) ,
\end{equation*}
\begin{equation*}S_{2}=\frac{1}{2\sqrt{-1}} \left( \overline{ \displaystyle\sum_{k,u=1}^{q}\left(b_{kk}^{ij}\left(Z\right)\left(\Re w_{kk}+\sqrt{-1}\left<Z_{k},Z_{k}\right>\right)+  D^{ij}_{kkuu}\left(\Re w_{kk}+\sqrt{-1}\left<Z_{k},Z_{k}\right>\right)\left(\Re w_{uu}+\sqrt{-1}\left<Z_{u},Z_{u}\right>\right)\right)}\right) ,
\end{equation*}
\begin{equation*}S_{3}=\frac{1}{2\sqrt{-1}}\left( \displaystyle\sum_{k,u=1\atop k\neq u}^{q}\left(b_{ku}^{ij}\left(Z\right)w_{ku}- \overline{b_{ku}^{ij}}\left(Z\right)\left(w_{uk}-2\sqrt{-1}\left<Z_{u},Z_{k} \right>    \right)\right)\right),
\end{equation*}
\begin{equation*}S_{4}=\frac{1}{2\sqrt{-1}}   \left(\displaystyle\sum_{k,u,k',u'=1\atop k\neq u,\hspace{0.1 cm} k'\neq u'}^{q}\left(D^{ij}_{kuk'u'}w_{ku}w_{k'u'}- \overline{D^{ij}_{kuk'u'}}\left( w_{uk}-2\sqrt{-1}\left<Z_{u},Z_{k} \right>  \right)\left( w_{u'k'}-2\sqrt{-1}\left<Z_{u'},Z_{k'} \right>  \right)\right)\right),
\end{equation*}
\begin{equation*}S_{5}= \frac{1}{2\sqrt{-1}}\left(\displaystyle\sum_{k,k',u'=1\atop k'\neq u'}^{q}\left(D^{ij}_{kkk'u'}\left(\Re w_{kk}+\sqrt{-1}\left<Z_{k},Z_{k}\right>\right)w_{k'u'}-  \overline{D^{ij}_{kkk'u'}}\left(\Re w_{kk}+\sqrt{-1}\left<Z_{k},Z_{k}\right>\right)\left(w_{u'k'}-2\sqrt{-1}\left<Z_{u'},Z_{k'} \right>\right)  \right)\right),
\end{equation*}\begin{equation*}S_{6}=\displaystyle\sum_{l=1}^{N}\overline{z}_{il}\left( \displaystyle\sum_{k,u=1\atop k\neq u}^{q}\left(a_{ku}^{jl}\left(Z\right)w_{ku}+a_{kk}^{jl}\left(Z\right)\left(\Re w_{kk}+\sqrt{-1}\left<Z_{k},Z_{k}\right>\right)\right) \right) ,
\end{equation*}\begin{equation*}S_{7}= \displaystyle\sum_{l=1}^{N}z_{il}\left( \displaystyle\sum_{k,u=1\atop k\neq u}^{q}\left(\overline{a_{ku}^{jl}\left(Z\right)}\left(w_{uk} -2\sqrt{-1}\left<Z_{u},Z_{k} \right>  \right)  + \overline{a_{kk}^{jl}\left(Z\right)\left(\Re w_{kk}+\sqrt{-1}\left<Z_{k},Z_{k}\right>\right)}\right) \right).
\end{equation*}

Since (\ref{99se}) holds and   the right-hand side  from (\ref{val1999})  does not depend on     products as  $\left(\Re w_{uu}\right) \cdot\left(\Re w_{u'u'}\right)$, for all $u,u'=1,\dots,q$, it suffices to focus on  the  difference $S_{1}-S_{2}$  in the left-hand side from  (\ref{val1999}). We obtain 
\begin{equation} \left(  D^{ij}_{kkuu   }\right)_{1\leq i,j  \leq q}=O_{q\times q},\quad    \mbox{for all $k,u=1,\dots,q$.} \label{kkqC} \end{equation}

Since the right-hand side  in (\ref{val1999})  does not depend on   terms involving $Z$ multiplied by  $  w_{ku}$ and $\left(\Re w_{uu}\right)$, for all $k,u=1,\dots,q$ such that $k \neq u $, it suffices to focus on  the   sums $S_{1}$, $S_{2}$ and $S_{3}$  in the left-hand side from (\ref{val1999}). We obtain
\begin{equation} \left(    b_{ku}^{ij}\left(Z\right)\right)_{1\leq i,j  \leq q}=O_{q\times q},\quad\mbox{for all $k,u=1,\dots,q$.} \label{tel1} \end{equation}
  
 Since the left-hand side  from (\ref{val1999})  does not depend on   terms involving terms of bi-degree $(1,2)$ in $\left(Z,\overline{Z}\right)$,  it suffices to focus on  the second sum  from the right-hand side
 in (\ref{val1999}). We obtain
\begin{equation}   \left(  b_{il}\left(Z\right)\right)_{ 1\leq i  \leq q\atop{1\leq l  \leq N}}=O_{q\times N} . \label{tel2} \end{equation}
 
Since the right-hand side from (\ref{val1999})  does not depend on   products as
 $w_{u'k'}\cdot \left<Z_{u},Z_{k}\right>$,  for all $k,u,k',u'=1,\dots,q$ such that $k \neq u $, $k'\neq u'$, we focus on the fourth sum from left-hand side from (\ref{val1999}). We obtain
\begin{equation}  \left(   D^{ij}_{ku k'u'  }\right)_{ 1\leq i,j \leq q}=O_{q\times q} , \quad  \mbox{for all $ k,u,k',u'=1,\dots,q$ such that $k \neq u $, $k'\neq u'$.} \label{kkqA} \end{equation}
 
Since the right-hand side from (\ref{val1999})  does not depend on   products as
$\left(\Re w_{kk}\right)\cdot w_{k'u'}$, for all $k,k',u'=1,\dots,q$ such that  $k'\neq u'$,  
we obtain
\begin{equation} \left(    D^{ij}_{kk k'u'  }\right)_{ 1\leq i,j \leq q}=\left(  \overline{D^{ij}_{kk u'k'  }}\right)_{ 1\leq i,j \leq q} ,\quad     \mbox{for all $k,k',u'=1,\dots,q$ such that $k'\neq u'$.} \label{kkqBB} \end{equation}

Next, we study the coefficient of the  product of terms as
 $\left<Z_{k},Z_{k}\right>\cdot \left<Z_{u'},Z_{k'} \right>$, for all $k,k',u'=1,\dots,q$ such that $k'\neq u'$, in the both sides from (\ref{val1999}) in the light of (\ref{99se}),(\ref{kkqC}),(\ref{kkqA})  and (\ref{kkqBB}). We obtain
\begin{equation} \left(D^{ij}_{kk k'u'  }\right)_{1\leq i,j  \leq q}=O_{q\times q} ,\quad     \mbox{for all $k,k',u'=1,\dots,q$ such that $k'\neq u'$.} \label{kkqB} \end{equation}

We identify    the coefficients of 
$w_{ku}$ and $\Re w_{kk}$ in $S_{6}$ and $S_{7}$ from (\ref{val1999}), for all $k,u=1,\dots,q$ such that $k\neq q$.     We obtain
\begin{equation}
\displaystyle\sum_{l=1}^{N}z_{il}\overline{a_{ku}^{il}\left(Z\right)}+\displaystyle\sum_{l=1}^{N}\overline{z}_{jl}a_{uk}^{jl}\left(Z\right)=0, \quad    \mbox{for all $i,j,k,u=1,\dots,q$.} \label{kk} 
\end{equation}

In particular, we write
\begin{equation}
\begin{split}&\left(a_{ku}^{il}\left(Z\right)\right)_{1\leq i  \leq q\atop{1\leq l  \leq N}} =\displaystyle\sum_{i'=1}^{q}\left(a_{ku}^{il}\left(Z_{i'}\right)\right)_{1\leq i  \leq q\atop{1\leq l  \leq N}},\quad\mbox{then the  equation (\ref{kk})  is equivalent to} \\&   \quad\quad\quad\quad\quad\quad\quad\quad\quad\quad\quad\quad\quad\quad
\quad\quad\quad\quad\quad\hspace{0.2 cm}
\displaystyle\sum_{l=1}^{N}\displaystyle\sum_{i'=1}^{q}z_{il}\overline{a_{ku}^{jl}\left(Z_{i'}\right)}+\sum_{i'=1}^{q}\displaystyle\sum_{l=1}^{N}\overline{z}_{jl}a_{uk}^{il}\left(Z_{i'}\right)=0, \quad    \mbox{for all $i,j,k,u=1,\dots,q$.} \label{kkkk} \end{split}
\end{equation}

It consequence, we obtain
\begin{equation}
\begin{split}&\left(a_{ku}^{il}\left(Z\right)\right)_{1\leq i \leq q\atop{1\leq l  \leq N}}^{1\leq k,u  \leq q}=\left(a_{ku}^{il}\left(Z_{i}\right)\right)_{1\leq i \leq q\atop{1\leq l  \leq N}}^{1\leq k,u  \leq q},\quad\mbox{then the  equation (\ref{val1999}) is equivalent to
} \\&   \quad\quad\hspace{0.1 cm}  2\sqrt{-1}\left( \displaystyle\sum_{k,u=1 }^{q} \left(\displaystyle\sum_{l=1}^{N}z_{il}\overline{a_{ku}^{jl}\left(Z_{j}\right)}  \right)\left<Z_{u},Z_{k}\right>\right)  =\left(   \displaystyle\sum_{l=1}^{N}\left(\varphi_{il}^{\star\star }\left(Z\right)\right)^{(2)} \overline{\left(\varphi_{jl}^{\star\star }\left(Z\right)\right)^{(2)}}\right),\quad\mbox{for all $i,j=1,\dots,q$, }\end{split}   \label{nuuit}   
\end{equation}
 where we  have used   the polynomial expansion
 \begin{equation}\left(\left(\varphi_{il}^{\star \star}\left(Z\right)\right)^{(2)}\right)_{1 \leq i \leq q \atop{1 \leq l \leq N}}=\left(\displaystyle\sum_{u,k=1}^{q} \displaystyle\sum_{\alpha,\beta=1\atop{\alpha \leq\beta}}^{N}b_{\alpha\beta\atop{uk}}^{il}z_{u\alpha}z_{k\beta}\right)_{1 \leq i \leq q \atop{1 \leq l \leq N}}
=\left( \displaystyle\sum_{u,k=1\atop{u \leq k}}^{q}\left(\left(\varphi_{il}^{uk}\right)^{\star \star}\left(Z_{u},Z_{k}\right)\right)\right)_{1 \leq i \leq q \atop{1 \leq l \leq N}}.\label{722s}
\end{equation}  

It suffices to study terms of similar bi-degree in the both sides of (\ref{nuuit}) in order to obtain that
\begin{equation*}\begin{split}\mbox{O}_{q\times N} &= \left(\left(\left(\varphi_{il}^{22}\right)^{\star \star}\left(Z_{2},Z_{2}\right)\right)\right)_{1 \leq i \leq q \atop{1 \leq l \leq N}}= \left(\left(\left(\varphi_{il}^{23}\right)^{\star \star}\left(Z_{2},Z_{3}\right)\right) \right)_{1 \leq i \leq q \atop{1 \leq l \leq N}}=\dots=\left(\left(\left(\varphi_{il}^{2q}\right)^{\star \star}\left(Z_{2},Z_{q}\right)\right) \right)_{1 \leq i \leq q \atop{1 \leq l \leq N}} \\&  \hspace{0.2 cm}\quad\quad\quad\quad\quad\quad\quad\quad \quad\quad\quad\quad\quad\quad\quad\quad\quad\quad\vdots \hspace{0.2 cm}\quad\quad\quad\quad\quad\quad\quad \quad\quad\quad\quad\quad\quad\quad\quad\quad\hspace{0.2 cm} \vdots\\& \hspace{0.2 cm}\quad\quad\quad\quad\quad\quad\quad\quad\quad\quad\quad\quad\quad=\left(\left(\left(\varphi_{il}^{q2}\right)^{\star \star}\left(Z_{q},Z_{2}\right)\right) \right)_{1 \leq i \leq q \atop{1 \leq l \leq N}}=\dots=\left(\left(\left(\varphi_{il}^{qq}\right)^{\star \star}\left(Z_{q},Z_{q}\right)\right) \right)_{1 \leq i \leq q \atop{1 \leq l \leq N}} .\end{split}
\end{equation*}

Finally, we take $k_{0}\in \left\{0,1,\dots,q\right\}$ such that 
 \begin{equation*}\begin{split}& \left(\varphi_{i1}^{\star \star}\left(Z\right),\dots, \varphi_{iN}^{\star \star}\left(Z\right)\right)\not\equiv (0,0,\dots,0),\quad\mbox{for all $i=1,\dots, k_{0}$,} \\& \left(\varphi_{i1}^{\star \star}\left(Z\right),\dots, \varphi_{iN}^{\star \star}\left(Z\right)\right)\equiv (0,0,\dots,0),\quad\mbox{for all $i=k_{0}+1,\dots, q$.}\end{split}
  \end{equation*}
   \end{proof}

\section{Application of  the Moving Point Trick from  Baouendi-Ebenfelt-Huang\cite{BEH}}   

In order to further normalize   the     Embedding (\ref{v11}),  we write the formal expansions  
\begin{equation}   G\left(Z,W\right) =\left(\displaystyle\sum_{ J\in\mathbb{N}^{q^{2}}\atop{ I\in\mathbb{N}^{q N } } }g_{ij}^{ IJ }\left(Z\right)W^{J}\right)_{1\leq i,j\leq q'}   \hspace{0.1 cm}\mbox{and}\hspace{0.1 cm} F\left(Z,W\right) =     \left(\displaystyle\sum_{ J\in\mathbb{N}^{q^{2}}\atop{  I\in\mathbb{N}^{qN}}  }f_{kl}^{ IJ }\left(Z\right)W^{J}\right)_{  1\leq k \leq q'\atop 1\leq  l\leq N'}.\label{v11V1}  \end{equation} 

 Their entries are formal power series. In particular, the coefficients of $W$ are defined by  
\begin{equation} \left(g_{ij}^{ IJ }\left(Z\right)\right)_{1\leq i,j\leq q'}  =\displaystyle\sum_{ I\in\mathbb{N}^{ qN}  }\left( c^{ IJ }_{ij}Z^{I}\right)_{1\leq i,j\leq q'}   \hspace{0.1 cm}\mbox{and}\hspace{0.1 cm} 
\left( f_{kl}^{ IJ  }\left(Z\right)\right)_{1\leq k\leq q'\atop{1\leq l\leq N'}}  =\displaystyle\sum_{  I\in\mathbb{N}^{ qN}}\left( d^{ IJ }_{kl}Z^{I}\right)_{1\leq k\leq q'\atop{1\leq l\leq N'}},\quad\hspace{0.12 cm}\mbox{for all  $J\in\mathbb{N}^{q^{2}}$.} \label{7771} 
\end{equation}  
  
We    extract from (\ref{v})  the terms of total degree $d$ in $\left(Z,\overline{Z}\right)$.    Given the model $\Im W=Z\overline{Z}^{t}$, we obtain
\begin{equation}  \displaystyle\sum_{ J\in\mathbb{N}^{q^{2}}, \hspace{0.05 cm}  I\in\mathbb{N}^{ qN}\atop \left|I\right|+2\left|J\right|=d} \left(\frac{g_{ij}^{ IJ }\left(Z\right)W^{J}-  \overline{g_{ji}^{ IJ }\left(Z\right)W^{J} } }{2\sqrt{-1}}\right)_{1\leq i,j\leq q'} = \displaystyle\sum_{l=1}^{N'}\left(\displaystyle\sum_{\substack{ J_{1},J_{2}\in\mathbb{N}^{q^{2}}, \hspace{0.05 cm} I_{1}, I_{2}\in\mathbb{N}^{qN} \\ \left|I_{1}\right|+2\left|J_{1}\right|+ \left|I_{2}\right|+2\left|J_{2}\right|=d   }}\left( f_{il}^{ I_{1}J_{1} }\left(Z\right)W^{J_{1}}\right)\left(\overline{f_{jl}^{ I_{2}J_{2} }\left(Z\right)W^{J_{2}}}\right)^{t}\right)_{1\leq i,j\leq q'}  . \label{yu123} 
\end{equation}
   
  Since the right-hand side from (\ref{yu123}) is hermitian, we obtain
\begin{equation} \left( g_{ij}^{ 0J }\left(Z\right)\right)_{1\leq i,j\leq q'}  =\left(\overline{g_{ji}^{ 0J }\left(Z\right) }\right)_{1\leq i,j\leq q'} ,\quad\mbox{for all   $J\in\mathbb{N}^{q^{2}}$.}\label{corr} 
\end{equation}
 
In particular, (\ref{corr}) describes  real-valued constants. Following Baouendi-Huang\cite{BH} and  Huang\cite{huang1}, it follows that 
\bl Up to compositions with holomorphic  automorphisms of $\mathcal{M}'$, we have \begin{equation}G\left(Z,W\right)=\begin{pmatrix}
 W & \mbox{O}_{q \times \left(q'-q\right) } \\ \mbox{O}_{\left(q'-q\right)\times q} & \mbox{O}_{\left(q'-q\right)\times \left(q'-q\right)}
\end{pmatrix}\hspace{0.1 cm}\mbox{and}\hspace{0.1 cm}F\left(\mbox{O}_{q\times N},W\right)=\mbox{O}_{q\times N}.\label{vvvq}\end{equation} \el 
\begin{proof} We define     the following transformations
\begin{equation*} \begin{split}&\quad
\sigma^{0}_{\left(Z_{0},W_{0}\right)}\left(Z,W\right)=\left(Z+Z_{0},W+W_{0}+2\sqrt{-1} \left< Z,Z_{0}\right>\right),\\& \tau_{\left(Z_{0},W_{0}\right)}^{\left(F,G\right)}\left(Z^{\star},W^{\star}\right)=\left(Z^{\star} -F\left(Z_{0},W_{0}\right), W^{\star}-\overline{G\left(Z_{0},W_{0}\right)}^{t}-2\sqrt{-1} \left<Z^{\star}, F\left(Z_{0},W_{0}\right) \right> \right).\end{split} \end{equation*}

Given  the point $P=\left(Z_{0},W_{0}\right)\in \mathcal{M}$ close to origin, we define
 $ \left(F,G\right)_{P}=\tau _{P}^{\left(F,G\right)}\circ \left(F,G\right)\circ \sigma_{P}^{0}=\left(F_{P}, G_{P}\right)$. In consequence, we obtain
$$\mbox{$\sigma^{0}_{P}(0)=P$, $\tau_{\left(F,G\right)\left(P\right)}^{\left(F,G\right)} \left(\left(F,G\right)\left(P\right)\right) =0$  and  $\det \left(\frac{\partial G_{11}\left(W\right)}{\partial W }\right)(0)\neq 0$.}$$

We consider using  (\ref{be4}) the  transformation
$\left(\tilde{G},\tilde{F}\right)=T_{2}\circ \left(G,F\right)$, 
for $T_{2}=T_{2}\left(P\right)$. It is  defined by  substractions of  terms in $W$   from the $F$-component of   the  Embedding  (\ref{v11}).  Varying the point $P\in \mathcal{M}$,   we obtain
 \begin{equation}G\left(\mbox{O}_{q\times N},W\right)=\begin{pmatrix}
 W & \mbox{O}_{q \times \left(q'-q\right) } \\ \mbox{O}_{\left(q'-q\right)\times q} & \mbox{O}_{\left(q'-q\right)\times \left(q'-q\right)}
\end{pmatrix}\hspace{0.1 cm}\mbox{and}\hspace{0.1 cm}F\left(\mbox{O}_{q\times N},W\right)=\mbox{O}_{q\times N}.\label{abc}\end{equation}
 
 The decisive argument is motivated by Hamada\cite{H}:   
 we study      the coefficients  of    the  terms  of the following type
\begin{equation*}\begin{split}Z^{I}&\left(\mbox{Re} w_{11}\right)^{j_{11}} w_{12}^{j_{12}}\dots 
 w_{1q}^{j_{1q}}\cdot\\&\overline{w}_{12}^{j_{21}}\left( \Re w_{22}\right)^{j_{22}}\dots  w_{2q}^{j_{2q}}\cdot\\&\quad \vdots\quad\quad \vdots\quad    \quad\hspace{0.1 cm}\quad\ddots\quad  \vdots \\& \overline{ w}_{1q}^{j_{q1}} \overline{  w}_{2q}^{j_{q2}}\dots \left(  \Re w_{qq}\right)^{j_{qq}},\quad\mbox{for all $I\in\mathbb{N}^{qN}$ and $J\in\mathbb{N}^{q^{2}}$.}\end{split}\end{equation*}

We identify    coefficients of   homogeneous terms in the left-hand side from (\ref{yu123}) using (\ref{abc}). Therefore  (\ref{vvvq}) holds.  
\end{proof}

\section{Application of the Method of Hamada\cite{H}}

In order to achieve the normal form (\ref{clase}), we study homogeneous terms in order to identify the final change of coordinates:
\subsection{Formal Expansions} 
Before beginning, we  write the following weighted formal expansions
\begin{equation}\begin{split}&\quad\left( f_{il}^{\star\star }\left(Z,W\right)\right)_{ 1\leq i\leq q\atop{1\leq l\leq N}} =Z  + \begin{pmatrix}
\displaystyle\sum_{k,u=1\atop{(k,u)\in\mathcal{S}}}^{q}\left(a_{ku}^{il}\left(Z_{i}\right)w_{ku}\right)_{ 1\leq i\leq k_{0}\atop{1\leq l\leq N}} \\ \mbox{O}_{\left(q-k_{0}\right)\times N}\end{pmatrix}+\displaystyle\sum_{\alpha \geq 4}\left( f_{kl}^{\left(\alpha\right)}\left(Z,W\right)\right)_{1\leq k  \leq q\atop{1\leq l  \leq N}} ,  \\& \quad\left( \varphi^{\star\star }_{il}\left(Z,W\right)\right)_{1\leq i\leq q\atop{1\leq l\leq N}}=\begin{pmatrix}
\left(  \displaystyle\sum_{k=1}^{q}\mathcal{L}_{il}^{k}\left(Z_{i},Z_{k}\right) \right)_{ 1\leq i\leq k_{0}\atop{1\leq l\leq N}}   
\\ \mbox{O}_{\left(q-k_{0}\right)\times N}
\end{pmatrix}   +\displaystyle\sum_{\alpha \geq 3}\left(\varphi_{kl}^{\left(\alpha\right)}\left(Z,W\right)\right)_{1\leq k  \leq q\atop{1\leq l  \leq N}},\end{split}\label{dub1}
   \end{equation} 
according to the assumptions related to (\ref{bub1}),  where we have used
 according to (\ref{nota}) formal expansions of the following type
 \begin{equation*}\begin{split}&  \left( f_{kl}^{\left(\alpha\right)}\left(Z,W\right)\right)_{1\leq k  \leq q\atop{1\leq l  \leq N}}=\displaystyle\sum_{\substack{ J  \in\mathbb{N}^{q^{2}}, \hspace{0.05 cm} I  \in\mathbb{N}^{qN} \\ \left|I \right|+2\left|J \right| =\alpha  }}\left( d_{kl}^{ I J  }\right)_{1\leq k  \leq q\atop{1\leq l  \leq N}} Z^{I}W^{J },\\&\left(\varphi_{kl}^{\left(\alpha\right)}\left(Z,W\right)\right)_{1\leq k  \leq q\atop{1\leq l  \leq N}}=\displaystyle\sum_{\substack{ J  \in\mathbb{N}^{q^{2}}, \hspace{0.05 cm} I  \in\mathbb{N}^{qN} \\ \left|I \right|+2\left|J \right| =\alpha  }}\left(\tilde{d}_{kl}^{ I J  }\right)_{1\leq k  \leq q\atop{1\leq l  \leq N}} Z^{I}W^{J }.\end{split}  
\end{equation*}

Given the model  (\ref{models}), from the equation $F_{2}\left(Z,W\right)\overline{F_{2}\left(Z,W\right)^{t}}=\mbox{O}_{\left( q'-q\right)\times \left( q'-q\right)}$ we obtain
\begin{equation}F_{2}\left(Z,W\right)=0.
\label{tre}
\end{equation}

 Next, we consider  $\alpha\in\{4m,4m+1,4m+2,4m+3 \hspace{0.1 cm}\vert \hspace{0.1 cm} m\in\mathbb{N}^{\star} \}$ in order to linearise the diagonal entries  in   (\ref{v})  according to     similar  computations   from  Hamada\cite{H}, where $\alpha\in\mathbb{N}^{\star}$.   In particular, we consider the following notations
 
\begin{equation*}\begin{split}&\hspace{0.25 cm}
\mbox{$J'=\left(j'_{ij}\right)_{1\leq i,j \leq q}\in\mathbb{N}^{q^{2}}$ any multi-index of degree $1$ derived from the multi-index $J=\left(j_{ij}\right)_{1\leq i,j \leq q}\in\mathbb{N}^{q^{2}}$,}\\&\hspace{0.17 cm}\mbox{$J''=\left(j''_{ij}\right)_{1\leq i,j \leq q}\in\mathbb{N}^{q^{2}}$   any multi-index of degree $2$ derived from the multi-index $J=\left(j_{ij}\right)_{1\leq i,j \leq q}\in\mathbb{N}^{q^{2}}$,} \\& \mbox{  $J'''=\left(j'''_{ij}\right)_{1\leq i,j \leq q}\in\mathbb{N}^{q^{2}}$ any multi-index of degree $1$ derived from the multi-index $J=\left(j_{ij}\right)_{1\leq i,j \leq q}\in\mathbb{N}^{q^{2}}$.}\end{split}
\end{equation*}

In consequence, we consider the Multinomial Coefficients denoted by
\begin{equation*}\begin{split}&
\hspace{0.1 cm}n\left(J^{'}_{0}\right) ,\quad\hspace{0.2 cm}\mbox{where}\hspace{0.1 cm}J^{'}_{0}\in\mathbb{N}^{q^{2}}\hspace{0.1 cm}\mbox{is a complement multi-index of degree $1$,} \\&
\hspace{0.1 cm} n\left(J^{''}_{0}\right),\quad\hspace{0.1 cm}\mbox{where}\hspace{0.1 cm}J^{''}_{0}\in\mathbb{N}^{q^{2}}\hspace{0.1 cm}\mbox{is a complement multi-index of degree $2$,}\\&
\hspace{0.1 cm} n\left(J^{'''}_{0}\right),\quad\hspace{0.1 cm}\mbox{where}\hspace{0.1 cm}J^{'''}_{0}\in\mathbb{N}^{q^{2}}\hspace{0.1 cm}\mbox{is a complement multi-index of degree $3$,}\end{split}
\end{equation*}
 associated to the   multi-index   $J=\left(j_{ij}\right)_{1\leq i,j \leq q}\in\mathbb{N}^{q^{2}}$. Given  $l_{1},\dots,l_{k}\in\mathbb{N}$   such that
$l_{1}+\cdots+l_{k}=n$, we define
 $$C^{n}_{l_{1},\dots,l_{k}}
=\frac{n!}{l_{1}!\cdots l_{k}!}.
$$
 
In order to consider formal powers series using the language of matrices, we use the   matrices
 \begin{equation}
W_{0}=\begin{pmatrix}\left\langle Z_{1}, Z_{1}  \right\rangle & \left\langle Z_{1}, Z_{2}  \right\rangle & \dots &
\left\langle Z_{1}, Z_{q}  \right\rangle \\  \left\langle Z_{2}, Z_{1}  \right\rangle & \left\langle Z_{2}, Z_{2}  \right\rangle & \dots &
\left\langle Z_{2}, Z_{q}  \right\rangle \\ \vdots & \vdots & \ddots & \vdots  \\ \left\langle Z_{q}, Z_{1}  \right\rangle & \left\langle Z_{q}, Z_{2}  \right\rangle & \dots &
\left\langle Z_{q}, Z_{q}  \right\rangle \end{pmatrix}\hspace{0.1 cm}\mbox{and}\hspace{0.1 cm}    \mathcal{W}=\begin{pmatrix}\Re w_{11} &   w_{12} & \dots &
  w_{1q} \\   \overline{w}_{12} & \Re w_{22} & \dots &
\overline{w}_{2q} \\ \vdots & \vdots & \ddots & \vdots  \\  \overline{w}_{1q} &   \overline{w}_{2q} & \dots &
\Re w_{qq}\end{pmatrix}.\label{dubluwe}
\end{equation}

\subsection{Cases $\alpha=4m+1 $ and $\alpha=4m+3 $} According to appropriate extractions of terms from (\ref{yu123}), we obtain
\begin{equation}\Re \left\{\left.\left(\displaystyle\sum_{l=1}^{N}\overline{z}_{kl} f^{(4m)}_{kl}\left(Z,W\right)\right\vert_{\left(Z,W\right)\in\mathcal{M}} \right)_{k_{0}+1\leq k \leq q} \right\}=\mbox{O}_{\left(q-k_{0}\right)\times 1},   \label{eKu1}\end{equation}
\begin{equation}\Re \left\{\left(\left.\displaystyle\sum_{l=1}^{N}\overline{z}_{kl} f^{(4m+2)}_{kl}\left(Z,W\right) \right\vert_{\left(Z,W\right)\in\mathcal{M}}\right)_{k_{0}+1\leq k \leq k_{0}} \right\}  =\mbox{O}_{\left(q-k_{0}\right)\times 1},   \label{eKu1r}\end{equation}
 where  we have used the following sums
\begin{equation*} \begin{split}&  \left( f^{(4m )}_{kl} \left(Z,W\right)\right)_{  1\leq k \leq q\atop{k_{0}+1\leq l  \leq N}}=\displaystyle\sum_{p=1}^{2m}\displaystyle\sum_{\substack{ J  \in\mathbb{N}^{q^{2}}    \\   \left|J \right| =2m-p  }}\left(\left(P_{ J}^{(2p)}\left(Z\right)\right)_{kl}\right)_{  1\leq k \leq q\atop{1\leq l  \leq N}}
W^{J},\\&     \left( f^{(4m+2 )}_{kl} \left(Z,W\right)\right)_{  1\leq k \leq q\atop{ 1\leq l  \leq N}}=\displaystyle\sum_{p=1}^{2m+1}\displaystyle\sum_{\substack{ J  \in\mathbb{N}^{q^{2}}    \\   \left|J \right| =2m+1-p  }}\left(\left(P_{ J}^{(2p)}\left(Z\right)\right)_{kl}\right)_{  1\leq k \leq q\atop{1\leq l  \leq N}}
W^{J}.   \end{split} 
\end{equation*} 

 Comparing similar terms in   (\ref{eKu1}) and (\ref{eKu1r}) for $W=\tilde{W}$, it follows that
 \begin{equation}   \left( f^{(4m )}_{kl}\left(Z,W\right) \right)_{k_{0}+1\leq k  \leq q\atop{1\leq l  \leq N}}=\mbox{O}_{\left(q-k_{0}\right)\times N} \hspace{0.1 cm}\mbox{and}\hspace{0.1 cm}  \left( f^{(4m+2 )}_{kl}\left(Z,W\right) \right)_{ k_{0}+1\leq k  \leq q\atop{1\leq l  \leq N}}=\mbox{O}_{\left(q-k_{0}\right)\times N}. \label{rezu111r}
\end{equation}

\subsection{Case $\alpha=4m+2$}According to appropriate extractions of terms from (\ref{yu123}), we obtain 
 \begin{equation}   2\Re \left\{\left(\left. \displaystyle\sum_{l=1}^{N}\overline{z}_{kl} f^{(4m+1)}_{kl}\left(Z,W\right)\right\vert_{\left(Z,W\right)\in\mathcal{M}}\right)_{ k_{0}+1\leq k \leq q} \right\} +  \displaystyle\sum_{l=1}^{N}\left(\left.\left| \varphi^{(2m+1)}_{kl}\left(Z,W\right)\right|^{2}\right\vert_{\left(Z,W\right)\in\mathcal{M}}   \right)_{ k_{0}+1\leq k \leq q}=\mbox{O}_{ \left(q-k_{0}\right) \times 1},  \label{ecu3} \end{equation}
where  we have used the following sums
 \begin{equation*}    \begin{split}& \left( f^{(4m+1)}_{kl} \left(Z,W\right)\right)_{   1\leq k \leq q\atop{1\leq l  \leq N}}=\displaystyle\sum_{p=0}^{2m}\displaystyle\sum_{\substack{ J  \in\mathbb{N}^{q^{2}}    \\   \left|J \right| =2m-p  }}\left(\left(P_{ J}^{(2p+1)}\left(Z\right)\right)_{kl}\right)_{  1\leq k \leq q\atop{1\leq l  \leq N}}
W^{J},\\&  \left( \varphi ^{(2m+1)}_{kl}\left(Z,W\right)\right)_{   1\leq k \leq q\atop{1\leq l  \leq N}} =\displaystyle\sum_{p=0}^{m}\displaystyle\sum_{\substack{ J  \in\mathbb{N}^{q^{2}}    \\   \left|J \right| =m-p  }}\left(\left(Q_{ J}^{(2p+1)}\left(Z\right)\right)_{kl}\right)_{   1\leq k \leq q\atop{1\leq l  \leq N}} 
 W^{J}.\end{split}  
\end{equation*}

Extracting multiples of terms of bi-degree $(1,1)$ in $\left(Z,\overline{Z}\right)$ from (\ref{ecu3}), we obtain

\begin{equation} \begin{split}&   2\Re\left\{ \left(\left.\displaystyle\sum_{l=1}^{N} \displaystyle\sum_{\substack{ J  \in\mathbb{N}^{q^{2}}    \\   \left|J \right| =2m  }} \left(P_{ J}^{(1)}\left(Z\right)\right)_{kl}\overline{z}_{kl} 
 \mathcal{W}^{J}\right\vert_{\left(Z,W\right)\in\mathcal{M}} \right)_{k_{0}+1\leq k \leq q} \right\} +\\&  \quad\quad\left( \displaystyle\sum_{l=1}^{N}\left. \left(    \displaystyle\sum_{\substack{ J  \in\mathbb{N}^{q^{2}}    \\   \left|J \right| = m  }}\left(Q_{ J}^{(1)}\left(Z\right)\right)_{kl} 
 \mathcal{W}^{J}  \right)\cdot  \overline{ \left(      \displaystyle\sum_{\substack{ J  \in\mathbb{N}^{q^{2}}    \\   \left|J \right| = m  }} \left(Q_{ J}^{(1)}\left(Z\right)\right)_{kl}  
\mathcal{W}^{J}  \right)}^{t}\right\vert_{\left(Z,W\right)\in\mathcal{M}}  \right)_{ k_{0}+1\leq k \leq q}
=\mbox{O}_{\left(q-k_{0}\right)\times 1}. \end{split}  \label{suplim1}
 \end{equation}

Extracting multiples of terms of bi-degree $(3,3)$ in $\left(Z,\overline{Z}\right)$ from (\ref{ecu3}),  we obtain 

\begin{equation}\begin{split} &   2\Re \left\{\left( \left.  \displaystyle\sum_{l=1}^{N}    \displaystyle\sum_{\substack{ J  \in\mathbb{N}^{q^{2}}    \\   \left|J \right| =2m   }}\left(-n\left(J^{''}_{0}\right)\right)\left(P_{ J}^{(1)}\left(Z\right)\right)_{kl}\overline{z}_{kl} W_{0}^{J^{''}_{0}}  \mathcal{W}^{J ''} \right\vert_{\left(Z,W\right)\in\mathcal{M}}\right)_{ k_{0}+1\leq k \leq q}
  \right\}  +\\&      \quad\quad\quad\quad\hspace{0.15 cm}\quad  \left(\left.\displaystyle\sum_{l=1}^{N} \left(    \displaystyle\sum_{\substack{ J  \in\mathbb{N}^{q^{2}}    \\   \left|J \right| = m  }}\left(n\left(J^{'}_{0}\right) Q_{ J}^{(1)}\left(Z\right)\right)_{kl} W_{0}^{J^{'}_{0}}   \mathcal{W}^{J'}    \right)\cdot \right.   \left. \overline{\left(       \displaystyle\sum_{\substack{ J \in\mathbb{N}^{q^{2}}    \\   \left|J \right| = m  }}  \left(Q_{ J}^{(1)}\left(Z\right)\right)_{kl}   
 \mathcal{W}^{J}  \right)}^{t} \right\vert_{\left(Z,W\right)\in\mathcal{M}}\right)_{ k_{0}+1\leq k \leq q}+\\&     \quad\quad\quad\quad\quad\quad\quad\quad\quad\quad \displaystyle\sum_{l=1}^{N}\left(\left. \left(    \displaystyle\sum_{\substack{ J  \in\mathbb{N}^{q^{2}}    \\   \left|J \right| = m-1  }} \left(Q_{ J}^{(3)}\left(Z\right)\right)_{kl}   
\mathcal{W}^{J}  \right)  \cdot       \overline{ \left(    \displaystyle\sum_{\substack{ J  \in\mathbb{N}^{q^{2}}    \\   \left|J \right| = m-1  }} \left(Q_{ J}^{(3)}\left(Z\right)\right)_{kl}  
\mathcal{W}^{J}  \right)}^{t} \right\vert_{\left(Z,W\right)\in\mathcal{M}}\right)_{ k_{0}+1\leq k \leq q} =\mbox{O}_{\left(q-k_{0}\right)\times 1}.\end{split} \label{suplim2}
 \end{equation} 

Subtracting suitable terms between (\ref{suplim1}) and (\ref{suplim2}), we obtain 
 \begin{equation*} 2\Re \left\{  \left(  \left.   \displaystyle\sum_{l=1}^{N}  \displaystyle\sum_{\substack{ J  \in\mathbb{N}^{q^{2}}    \\   \left|J \right| =2m  }} \left(P_{ J}^{(1)}\left(Z\right)\right)_{kl}\overline{z}_{kl}    \mathcal{W}^{J}\right\vert_{\left(Z,W\right)\in\mathcal{M}}
   \right)_{ k_{0}+1\leq k \leq q} \right\} =\mbox{O}_{\left(q-k_{0}\right) \times 1}.  
 \end{equation*}

Extracting multiples of terms of bi-degree $(2,2)$ in $\left(Z,\overline{Z}\right)$ from (\ref{ecu3}), we obtain
\begin{equation*}  2\Re \left\{\sqrt{-1} \left(  \left.   \displaystyle\sum_{l=1}^{N}  \displaystyle\sum_{\substack{ J   \in\mathbb{N}^{q^{2}}    \\   \left|J \right| =2m   }} \left(P_{ J}^{(1)}\left(Z\right)\right)_{kl}\overline{z}_{kl}W_{0}^{J^{'}_{0}}   \mathcal{W}^{J'}\right\vert_{\left(Z,W\right)\in\mathcal{M}}
   \right)_{ k_{0}+1\leq k \leq q} \right\} =\mbox{O}_{\left(q-k_{0}\right) \times 1}.  
 \end{equation*}
 
   In particular, it follows that
\begin{equation}    \left(\left( P_{ J}^{( 1)}\left(Z\right) \right)_{kl} \right)_{ k_{0}+1\leq k  \leq q\atop{1\leq l  \leq N}}=\mbox{O}_{\left(q-k_{0}\right)\times N},\quad\mbox{for all    $J  \in\mathbb{N}^{q^{2}}$  with $\left|J \right| =2m$.}
\label{ecuA5se}
\end{equation} 

Returning to (\ref{suplim1}), it follows that
 \begin{equation}    \left(\left(   Q_{ J}^{( 1)}\left(Z\right) \right)_{kl}  \right)_{ k_{0}+1\leq k  \leq q\atop{1\leq l  \leq N}} =\mbox{O}_{\left(q-k_{0}\right)\times N},\quad\mbox{for all   $J  \in\mathbb{N}^{q^{2}}$  with $\left|J \right| =m$.} 
\label{ecuA5}
\end{equation}

Extracting multiples of terms of bi-degree $(3,1)$ in $\left(Z,\overline{Z}\right)$ from (\ref{ecu3}), we obtain
\begin{equation*} 2\Re \left\{\left(\left.\displaystyle\sum_{l=1}^{N} \displaystyle\sum_{\substack{ J  \in\mathbb{N}^{q^{2}}    \\   \left|J \right| =2m-1  }} \left(P_{ J}^{(3)}\left(Z\right)\right)_{kl}\overline{z}_{kl}   \mathcal{W}^{J}\right\vert_{\left(Z,W\right)\in\mathcal{M}} \right)_{k_{0}+1\leq k \leq q} \right\} =\mbox{O}_{\left(q-k_{0}\right)\times 1},    \end{equation*}

Moving forward, it follows that

\begin{equation*}\begin{split}&\left.2\Re \left\{\left(\left.\displaystyle\sum_{l=1}^{N} \displaystyle\sum_{\substack{ J  \in\mathbb{N}^{q^{2}}    \\   \left|J \right| =2m-2  }} \left(P_{ J}^{(5)}\left(Z\right)\right)_{kl}\overline{z}_{kl}   \mathcal{W}^{J}\right\vert_{\left(Z,W\right)\in\mathcal{M}}\right)_{k_{0}+1\leq k \leq q}\right\}+\displaystyle\sum_{l=1}^{N}\left(\displaystyle\sum_{\substack{ J  \in\mathbb{N}^{q^{2}}    \\   \left|J \right| = m -1 }} \left(Q_{ J}^{(3)}\left(Z\right)\right)_{kl}     \mathcal{W}^{J}\right)\cdot \right.  \\& \quad\quad\quad\quad\quad\quad\quad\quad\quad\quad\quad \hspace{0.26 cm}\quad\quad\quad\quad\quad\quad\quad\quad\quad\quad\quad\quad\quad\quad  \left.\left. \overline{ \left(\displaystyle\sum_{\substack{ J  \in\mathbb{N}^{q^{2}}    \\   \left|J \right| = m -1 }} \left(Q_{ J}^{(3)}\left(Z\right)\right)_{kl}     \mathcal{W}^{J}\right)} \right\vert_{\left(Z,W\right)\in\mathcal{M}} \right)_{k_{0}+1\leq k \leq q} =\mbox{O}_{\left(q-k_{0}\right)\times 1}, \end{split}   \end{equation*}
 $$\vdots\quad\quad\quad\quad\quad\quad\quad\quad\quad\quad\quad\quad    \vdots\quad\quad\quad\quad \quad\quad\quad\quad \quad\quad\quad\quad \vdots\quad\quad\quad\quad \quad\quad\quad\quad \quad\quad\quad\quad\vdots\quad\quad\quad\quad\quad \quad\quad\quad \quad\quad\quad\quad\vdots $$

In consequence, it follows that
\begin{equation}\begin{split}& \left(  \left( P_{ J}^{(2p+1)}\left(Z\right)   \right)_{kl}\right)_{ k_{0}+1\leq k  \leq q\atop{1\leq l  \leq N}} =\mbox{O}_{\left(q-k_{0}\right)\times N},\quad\mbox{for all $p=1,\dots, 2m$ and  $J  \in\mathbb{N}^{q^{2}}$  with $\left|J \right| =2m-p$,}\\& \left(  \left( Q_{ J}^{(2p+1)}\left(Z\right)\right)_{kl}\right)_{k_{0}+1 \leq k  \leq q\atop{1\leq l  \leq N}} =\mbox{O}_{\left(q-k_{0}\right)\times N} ,\quad\mbox{for all $p=1,\dots, m$ and  $J  \in\mathbb{N}^{q^{2}}$  with $\left|J \right| =m-p$.}  \end{split}
\label{ecuA3sec}
\end{equation}  

\subsection{Case $\alpha=4m+4$} According to appropriate extractions of terms from (\ref{yu123}), we obtain 
\begin{equation}  2\Re \left\{ \left(\left.\displaystyle\sum_{l=1}^{N}\overline{z}_{kl} f^{(4m+3)}_{kl}\left(Z,W\right)\right\vert_{\left(Z,W\right)\in\mathcal{M}}\right) _{ k_{0}+1\leq k \leq q}\right\}+ \left( \left. \displaystyle\sum_{l=1}^{N}\left| \varphi^{(2m+2)}_{kl}\left(Z,W\right)\right|^{2}\right\vert_{\left(Z,W\right)\in\mathcal{M}}\right)_{ k_{0}+1\leq k \leq q}=\mbox{O}_{\left(q-k_{0}\right)\times 1},  \label{ecu56} \end{equation}
  where  we have used the following sums
 \begin{equation*}    \begin{split}& \left( f^{(4m+3)}_{kl} \left(Z,W\right)\right)_{   1\leq k \leq q\atop{1\leq l  \leq N}}
=\displaystyle\sum_{p=0}^{2m+1}\displaystyle\sum_{\substack{ J  \in\mathbb{N}^{q^{2}}    \\   \left|J \right| =2m+1-p  }}\left(\left(P_{ J}^{(2p+1)}\left(Z\right)\right)_{kl}\right)_{   1\leq k \leq q\atop{1\leq l  \leq N}}
 W^{J} ,  \\&  \left( \varphi^{(2m+2)}_{kl}\left(Z,W\right)\right)_{  1\leq k \leq q\atop{1\leq l  \leq N}}
 =\displaystyle\sum_{p=1}^{m+1}\displaystyle\sum_{\substack{ J  \in\mathbb{N}^{q^{2}}    \\   \left|J \right| =m+1-p  }}\left(\left(Q_{ J}^{(2p)}\left(Z\right)\right)_{kl}\right)_{  1\leq k \leq q\atop{1\leq l  \leq N}}
 W^{J} .\end{split}   
\end{equation*}
 
Extracting multiples of terms of bi-degree $(1,1)$ in $\left(Z,\overline{Z}\right)$ from (\ref{ecu56}),  it follows that \begin{equation} 2\Re \left\{    \displaystyle\sum_{l=1}^{N}\left( \left.\displaystyle\sum_{\substack{ J  \in\mathbb{N}^{q^{2}}    \\   \left|J \right| =2m+1  }}\left( P_{ J}^{(1)}\left(Z\right)\right)_{kl}\overline{z}_{kl} 
 \mathcal{W}^{J}  \right\vert_{\left(Z,W\right)\in\mathcal{M}}  \right)_{ k_{0}+1\leq k \leq q} \right\} =\mbox{O}_{\left(q-k_{0}\right)\times 1}. \label{wr}  
 \end{equation}
 
Extracting multiples of terms of bi-degree $(4,4)$ in $\left(Z,\overline{Z}\right)$ from (\ref{ecu56}),  we obtain
 \begin{equation}\begin{split}&   \left(2\Re \left\{ \left. \displaystyle\sum_{l=1}^{N}  \displaystyle\sum_{\substack{ J   \in\mathbb{N}^{q^{2}}    \\   \left|J   \right| =2m+1  }}\left(-n\left(J_{0}^{'}\right)\right)\sqrt{-1}\left(P_{ J }^{(1)}\left(Z\right)\right)_{kl}\overline{z}_{kl} W_{0}^{J^{'}_{0}}  \mathcal{W}^{J^{'} } \right\vert_{\left(Z,W\right)\in\mathcal{M}}  
  \right\}\right)_{ k_{0}+1\leq k \leq q}+\\&   \quad\quad\quad\quad\quad\quad\quad \left( \left. \displaystyle\sum_{l=1}^{N}  \left(       \displaystyle\sum_{\substack{ J  \in\mathbb{N}^{q^{2}}    \\   \left|J \right| = m   }} n\left(J^{'}_{0}\right)\left(Q_{ J}^{(2)}\left(Z\right)\right)_{kl} W_{0}^{J_{0}^{'}}      \mathcal{W}^{J'}  \right)    \cdot      
                       \overline{   \left(    \displaystyle\sum_{\substack{ J  \in\mathbb{N}^{q^{2}}    \\   \left|J \right| = m -1 }}  \left(Q_{ J}^{(2)}\left(Z\right)\right)_{kl}    
  \mathcal{W}^{J}  \right) }^{t}\right\vert_{\left(Z,W\right)\in\mathcal{M}}  \right)_{ k_{0}+1\leq k \leq q}  +\\&        \quad\quad\quad\quad\quad\quad\quad\quad\quad\quad\quad\quad   \left(\left.\displaystyle\sum_{l=1}^{N}\left(    \displaystyle\sum_{\substack{ J  \in\mathbb{N}^{q^{2}}    \\   \left|J \right| = m -1 }} \left(Q_{ J}^{(4)}\left(Z\right)\right)_{kl}  
     \mathcal{W}^{J}  \right) \cdot         
    \overline{ \left(      \displaystyle\sum_{\substack{ J  \in\mathbb{N}^{q^{2}}    \\   \left|J \right| = m -1 }} \left(Q_{ J}^{(4)}\left(Z\right)\right)_{kl}   
  \mathcal{W}^{J}   \right)}^{t}\right\vert_{\left(Z,W\right)\in\mathcal{M}}  \right)_{ k_{0}+1\leq k \leq q}    =\mbox{O}_{\left(q-k_{0}\right)\times 1}. \end{split}\label{sup}
 \end{equation}

Extracting multiples of terms of bi-degree $(2,2)$ in $\left(Z,\overline{Z}\right)$ from (\ref{ecu56}),  we obtain

\begin{equation}\begin{split}&    2\Re \left\{\displaystyle\sum_{l=1}^{N}   \left.\left( \displaystyle\sum_{\substack{ J    \in\mathbb{N}^{q^{2}}    \\   \left|J   \right| =2m+1  }} n\left(J^{'}_{0}\right)\sqrt{-1}\left(P_{ J}^{(1)}\left(Z\right)\right)_{kl}\overline{z}_{kl} W_{0}^{J_{0}} 
 \mathcal{W}^{J '}\right\vert_{\left(Z,W\right)\in\mathcal{M}} \right)_{ k_{0}+1\leq k \leq q}\right\}+ \\& \quad\hspace{0.2 cm}    \left( \displaystyle\sum_{l=1}^{N}  \left(     \displaystyle\sum_{\substack{ J  \in\mathbb{N}^{q^{2}}    \\   \left|J \right| = m   }} \left(Q_{ J}^{(2)}\left(Z\right)\right)_{kl}   
  \mathcal{W}^{J}   \right)\cdot       \left.\overline{ \left(     \displaystyle\sum_{\substack{ J  \in\mathbb{N}^{q^{2}}    \\   \left|J \right| = m   }}  \left(Q_{ J}^{(2)}\left(Z\right)\right)_{kl} \right) 
 \mathcal{W}^{J}  }^{t}\right\vert_{\left(Z,W\right)\in\mathcal{M}}  \right)_{ k_{0}+1\leq k \leq q}
    =\mbox{O}_{\left(q-k_{0}\right)\times 1}. \end{split}  \label{sup1}
 \end{equation}

Subtracting suitable terms between (\ref{sup1}) and (\ref{sup}), it follows that
 \begin{equation} 2\Re \left\{ \left.   \displaystyle\sum_{l=1}^{N}\left( \displaystyle\sum_{\substack{ J  \in\mathbb{N}^{q^{2}}    \\   \left|J \right| =2m+1  }}\left( P_{ J}^{(1)}\left(Z\right)\right)_{kl}\overline{z}_{kl} 
 \mathcal{W}^{J}\right\vert_{\left(Z,W\right)\in\mathcal{M}} \right)_{ k_{0}+1\leq k \leq q} \right\} =\mbox{O}_{\left(q-k_{0}\right)\times 1}. \label{wr1}  
 \end{equation}
 
Since (\ref{wr}) and (\ref{wr1}) hold, it follows that
 \begin{equation}   \left(\left( P_{ J}^{( 1)}\left(Z\right) \right)_{kl} \right)_{ k_{0}+1\leq k  \leq q\atop{1\leq l  \leq N}}=\mbox{O}_{ \left(q-k_{0}\right)\times N},\quad\mbox{for all    $J  \in\mathbb{N}^{q^{2}}$  with $\left|J \right| =2m+1$.} 
\label{ecuA8}
\end{equation}

Returning to (\ref{sup1}), it follows that
 \begin{equation}   \left(\left(   Q_{ J}^{( 2)}\left(Z\right)   \right)_{kl}\right)_{ k_{0}+1\leq k  \leq q\atop{1\leq l  \leq N}} =\mbox{O}_{ \left(q-k_{0}\right)\times N},\quad\mbox{for all   $J  \in\mathbb{N}^{q^{2}}$  with $\left|J \right| =m+1$.} 
\label{ecuA8se}    
\end{equation} 

Extracting multiples of terms of bi-degree $(4,1)$ in $\left(Z,\overline{Z}\right)$ from (\ref{ecu56}), we obtain
\begin{equation*}\begin{split}&\left.2\Re \left\{\left(\left.\displaystyle\sum_{l=1}^{N} \displaystyle\sum_{\substack{ J  \in\mathbb{N}^{q^{2}}    \\   \left|J \right| =2m  }} \left(P_{ J}^{(3)}\left(Z\right)\right)_{kl}\overline{z}_{kl}    \mathcal{W}^{J}\right\vert_{\left(Z,W\right)\in\mathcal{M}}\right)_{k_{0}+1\leq k \leq q}\right\}\right.+\\&\left.\left(\displaystyle\sum_{l=1}^{N}\left(\displaystyle\sum_{\substack{ J  \in\mathbb{N}^{q^{2}}    \\   \left|J \right| = m  }} \left(Q_{ J}^{(2)}\left(Z\right)\right)_{kl}     \mathcal{W}^{J}\right)\cdot    \overline{ \left(\displaystyle\sum_{\substack{ J  \in\mathbb{N}^{q^{2}}    \\   \left|J \right| = m -1 }} \left(Q_{ J}^{(2)}\left(Z\right)\right)_{kl}     \mathcal{W}^{J}\right)} \right\vert_{\left(Z,W\right)\in\mathcal{M}}\right)_{k_{0}+1\leq k \leq q}  =\mbox{O}_{\left(q-k_{0}\right)\times 1}. \end{split}   \end{equation*}

Moving forward, it follows that
\begin{equation*} 2\Re \left\{\left(\left.\displaystyle\sum_{l=1}^{N} \displaystyle\sum_{\substack{ J  \in\mathbb{N}^{q^{2}}    \\   \left|J \right| =2m -1 }} \left(P_{ J}^{(5)}\left(Z\right)\right)_{kl}\overline{z}_{kl}    \mathcal{W}^{J}\right\vert_{\left(Z,W\right)\in\mathcal{M}}  \right)_{k_{0}+1\leq k \leq q} \right\} =\mbox{O}_{\left(q-k_{0}\right)\times 1},    \end{equation*}
 $$    \vdots\quad\quad\quad\quad \quad\quad\quad\quad \quad\quad\quad\quad \vdots\quad\quad\quad\quad \quad\quad\quad\quad \quad\quad\quad\quad\vdots  $$
 
In consequence,   we obtain
\begin{equation}\begin{split}&\left(  \left(  P_{ J}^{(2p+1)}\left(Z\right)   \right)_{kl}\right)_{ k_{0}+1\leq k \leq q\atop{1\leq l  \leq N}} =\mbox{O}_{ \left(q-k_{0}\right)\times N},\quad\quad\mbox{for all $p=0,\dots, 2m+1$ and  $J  \in\mathbb{N}^{q^{2}}$  with $\left|J \right| =2m+1-p$,}\\& \left(\left( Q_{ J}^{(2p )}\left(Z\right) \right)_{kl} \right)_{ k_{0}+1\leq k  \leq q\atop{1\leq l  \leq N}}=\mbox{O}_{ \left(q-k_{0}\right)\times N} , \quad\quad\quad \mbox{for all $p=1,\dots, m$ and  $J  \in\mathbb{N}^{q^{2}}$  with $\left|J \right| =m+1-p$.}\end{split}
\label{ecuA6se}
\end{equation}

\section{Changes of Coordinates from Huang-Ji\cite{HJ}}   In order to consider diagonalisations of matrices, we  consider the matrices 
 \begin{equation}\begin{split}&\left(\mathcal{B}_{ku}^{i}\right)_{1 \leq k,u  \leq q}^{1 \leq  i \leq q}=-2\sqrt{-1}\left(\left(a_{ku}^{i1},a_{ku}^{i2},\dots,a_{ku}^{iN}\right)\right)_{1 \leq k,u  \leq q}^{1 \leq  i \leq q},\quad\mbox{then the  equation (\ref{kk})  is equivalent to}\\&\quad\quad\quad\quad\quad\quad\quad\quad\quad\quad\quad\quad\quad\quad\quad\quad\quad\quad\quad\quad\quad\quad\quad\quad\quad\quad\quad\quad\quad\quad
 \left<Z_{i},\mathcal{B}_{ku}^{j}\left(Z_{j}\right)\right>=\left<\mathcal{B}_{uk}^{i}\left(Z_{i}\right),Z_{j}\right>,\quad \mbox{for all
 $ k,u,i,j=1,\dots,q$. }\end{split}  \label{5050seee}
\end{equation} 
   
It suffices to assume   $i= k=1$.   In particular, we consider  $U$   a unitary matrix defining the diagonalisation of the hermitian matrix $\mathcal{B}^{1}_{11}$ in order to consider new coordinates preserving the model $\mathcal{M}$. We define

  \begin{equation}\left(W,Z\right):= \begin{pmatrix} W & \begin{pmatrix}
\begin{pmatrix} z_{11},z_{12},\dots, z_{1N} \end{pmatrix}\cdot U^{-1} \\ \begin{pmatrix} z_{21},z_{22},\dots, z_{2N} \end{pmatrix}\cdot U^{-1} \\ \begin{matrix}
 \vdots & \ddots & \vdots
\end{matrix} \\  \begin{pmatrix}  z_{q1},z_{q2},\dots, z_{qN}  \end{pmatrix}\cdot U^{-1}
 \end{pmatrix}\end{pmatrix}.\label{Q11}
\end{equation}

It suffices to apply the procedure of Huang-Ji\cite{HJ}    in the light of (\ref{nuuit}) and (\ref{722s}).
 It exists the   matrix $ \tilde{B}$ such that
  \begin{equation*}\tilde{B} \overline{\tilde{B} }^{t}= \alpha \cdot I_{N} \hspace{0.1 cm}\mbox{and}\hspace{0.1 cm}\begin{pmatrix}\varphi_{i1}^{\star \star}\left(Z_{1}\right)\\ \varphi_{i2}^{\star \star}\left(Z_{1}\right)\\ \vdots \\ \varphi_{iN}^{\star \star}\left(Z_{1}\right)\end{pmatrix}=z_{11}\tilde{B}  Z_{1}^{t}.
  \end{equation*}
 
We move forward following Huang-Ji\cite{HJ} in order to implement the   coordinates
\begin{equation}   \left(W',Z'\right):=\left(W',\begin{pmatrix}\begin{pmatrix}
 z'_{11},z'_{12},\dots, z'_{1N}
 \end{pmatrix}\\ \begin{pmatrix}
 z'_{21},z'_{22},\dots, z'_{2N}
 \end{pmatrix}  \\   \vdots \\ \begin{pmatrix}
 z'_{q'1},z'_{q'2},\dots, z'_{q'N}
 \end{pmatrix}
\end{pmatrix},\begin{pmatrix}\begin{pmatrix}
 z'_{1\hspace{0.05 cm}  q+1},z'_{1\hspace{0.05 cm}  q+2},\dots, z'_{1\hspace{0.05 cm} N}
 \end{pmatrix}\cdot \frac{\tilde{B} }{\sqrt{\alpha } }\\  \begin{pmatrix}
 z'_{2\hspace{0.05 cm}  q+1},z'_{2\hspace{0.05 cm}  q+1},\dots, z'_{2\hspace{0.05 cm}  N}
 \end{pmatrix} \\\vdots   \\ \begin{pmatrix}
 z'_{q'\hspace{0.05 cm}  q+1},z'_{q'\hspace{0.05 cm}  q+2},\dots, z'_{q'\hspace{0.05 cm}  N}
 \end{pmatrix}
\end{pmatrix} \right). \label{Q1}
\end{equation}

Finally, we consider the   change of coordinates 
\begin{equation}
\left(W',Z'\right)=\left(\begin{pmatrix}  \alpha  \cdot w'_{11}   &  \sqrt{\alpha}  \cdot   w'_{12}&   \dots &   \sqrt{\alpha}  \cdot  w'_{1\hspace{0.05 cm} q'}\\         \sqrt{\alpha}  \cdot w'_{2\hspace{0.05 cm}1}  &    w'_{2\hspace{0.05 cm}1}&    \dots &  w'_{2\hspace{0.05 cm} q'} \\ \vdots    & \vdots   &\ddots &\vdots  \\     \sqrt{\alpha}  \cdot  w'_{q'1}&w'_{q'\hspace{0.05 cm} 2} & \dots & w'_{q'\hspace{0.05 cm} q'}\end{pmatrix},\begin{pmatrix}Z'_{1}\cdot \sqrt{\alpha } \\   Z'_{2} \\ \vdots \\ Z'_{q'}\end{pmatrix}\right),\label{Q2}
\end{equation}
where have used the identification

\begin{equation}
\left(\left(w_{ij}\right)_{1\leq i,j \leq q} ,\left(z_{ij}\right)_{1\leq i  \leq q\atop{1\leq  j \leq N}} \right):=\left(\begin{pmatrix}
\frac{w_{11}}{ \alpha }& \frac{w_{12}}{\sqrt{\alpha } } & \dots & \frac{w_{1q}}{\sqrt{\alpha } } \\  \frac{w_{21}}{\sqrt{\alpha } }&  w_{22}  & \dots &  w_{2q}  \\ \vdots & \vdots & \ddots & \vdots \\ \frac{w_{q1}}{\sqrt{\alpha } }&  w_{q2}  & \dots &  w_{qq} 
\end{pmatrix},\begin{pmatrix}\frac{Z_{1}}{\sqrt{\alpha }}\\  Z_{2} \\\vdots \\  Z _{q}  
\end{pmatrix}\right).\label{Q3}
\end{equation}

  It suffices to assume  $(k,i)\in \mathcal{S}$ only for $k=1,\dots,q_{0}$, where $1\leq q_{0}\leq q$, in order to consider the sets of coordinates $\left\{\left(V_{ki},U_{ki}\right)\right\}_{(k,i)\in\mathcal{S}}\subset\mathcal{M}_{ q^{2}\times  q^{2} }\left(\mathbb{C}\right)\times \mathcal{M}_{Nq\times Nq }\left(\mathbb{C}\right)$ and $\left\{\left(\tilde{V}_{ki},\tilde{U}_{ki}\right)\right\}_{(k,i)\in\mathcal{S}}\subset\mathcal{M}_{ {q}^{2}\times  {q}^{2} }\left(\mathbb{C}\right)\times \mathcal{M}_{Nq\times Nq}\left(\mathbb{C}\right)$
according to (\ref{Q11}),(\ref{Q1}),(\ref{Q2}) and (\ref{Q3}),   such that
\begin{equation} \left( \tilde{U}_{k_{0}i_{0}}\otimes  \left( g^{\star \star}_{ij}\left(V_{k_{0}i_{0}}\otimes Z,U_{k_{0}i_{0}}\otimes W\right)\right)_{1\leq i,j  \leq q}\right)_{1\leq i,j  \leq q}= \tilde{U}_{k_{0}i_{0}}\otimes \left( U_{k_{0}i_{0}}\otimes W\right) + \mbox{O}_{\mbox{wt}}(5),\label{seco1}\end{equation}
\begin{equation}\begin{split} &\left(\tilde{V}_{k_{0}i_{0}}\otimes \left(  \left( F^{\star\star} \left(V_{k_{0}i_{0}}\otimes Z,U_{k_{0}i_{0}}\otimes W\right)\right)_{1\leq i\leq q\atop{1\leq l\leq N}}   \right)\right)_{i=i_{0}\atop{ 1\leq l  \leq N }} =\left(\left(V_{k_{0}i_{0}}\otimes Z\right)_{11} ,\left(V_{k_{0}i_{0}}\otimes Z\right)_{12},\dots, \left(V_{k_{0}i_{0}}\otimes Z\right)_{1N}\right)\\& \quad\quad\quad\quad\quad+ \left( \frac{\sqrt{-1}}{2}\left(V_{k_{0}i_{0}}\otimes Z\right)\left( U_{k_{0}i_{0}}\otimes W\right)_{11},0,\dots, 0\right)    +\left( \displaystyle\sum_{k,u=1\atop{k\neq k_{0}}}^{q}\left(\tilde{a}_{ku}^{i_{0}l}\left(V_{k_{0}i_{0}}\otimes Z\right)\left( U_{k_{0}i_{0}}\otimes W\right)_{ku}\right)_{}\right)_{  1\leq l  \leq N  }+  \mbox{O}_{\mbox{wt}}(4),   \end{split}\label{seco11}
\end{equation}
where the linear form $  \tilde{a}_{ku}^{i_{0}l}\left(V_{k_{0}i_{0}}\otimes Z\right)$ is    derived from  the linear form $  a_{ku}^{i_{0}l}\left(V_{k_{0}i_{0}}\otimes Z\right)$ according to (\ref{Q1}),(\ref{Q2}) and (\ref{Q3}), for all $l=1,\dots,N$ and $k,u=1,\dots,q$ with $k\neq k_{0}$, and  
\begin{equation}\begin{split}& \left(\tilde{V}_{k_{0}i_{0}}\otimes \left(  \left( F^{\star\star} \left(V_{k_{0}i_{0}}\otimes Z,U_{k_{0}i_{0}}\otimes W\right)\right)   \right)_{1\leq i\leq q\atop{1\leq l\leq N}}\right)_{i=i_{0}\atop{ N+1\leq l  \leq 2N}} =\left(V_{k_{0}i_{0}}\otimes Z\right)_{11}\cdot\left(\left(V_{k_{0}i_{0}}\otimes Z\right)_{11},\left(V_{k_{0}i_{0}}\otimes Z\right)_{12},\dots, \left(V_{k_{0}i_{0}}\otimes Z\right)_{1N}\right)\\& \quad\quad\quad\quad\quad\quad \quad\quad\quad\quad\quad\quad\quad\quad\quad   +  \displaystyle\sum_{l=1}^{N}\displaystyle\sum_{k=1\atop{k\neq k_{0}}}^{q} \left(\tilde{\varphi}_{1l}^{1k}\right)\left(\left(V_{k_{0}i_{0}}\otimes Z\right)_{1},\left(V_{k_{0}i_{0}}\otimes Z\right)_{k}\right) \overline{\left( \left(\tilde{\varphi}_{1l}^{1k}\right)\left(\left(V_{k_{0}i_{0}}\otimes Z\right)_{1},\left(V_{k_{0}i_{0}}\otimes Z\right)_{k}\right)\right)}  +  \mbox{O}_{\mbox{wt}}(5),   \end{split}\label{seco111}
\end{equation}
where the polynomial $ \left(\tilde{\varphi}_{i_{0}l}^{1k}\right)\left(\left(V_{k_{0}i_{0}}\otimes Z\right)_{1},\left(V_{k_{0}i_{0}}\otimes Z\right)_{k}\right)   $  is    derived from  the polynomial $ \left(\varphi_{1l}^{i_{0}k}\right)^{\star \star}\left(\left(V_{k_{0}i_{0}}\otimes Z\right)_{1},\left(V_{k_{0}i_{0}}\otimes Z\right)_{k}\right)   $ according to (\ref{Q1}),(\ref{Q2}) and (\ref{Q3}), for all $l=1,\dots,N$ and $k =1,\dots,q$ with $k\neq k_{0}$, for all $i_{0}=1,\dots,q_{0}$.

They are crucial in order to write  the following expressions
 \begin{equation}\left(\left(\varphi_{i1}^{\star \star}\left(Z\right)\right)^{(2)},\left(\varphi_{i2}^{\star \star}\left(Z\right)\right)^{(2)},\dots, \left(\varphi_{iN}^{\star \star}\left(Z\right)\right)^{(2)}\right) = \displaystyle\sum_{k=1}^{q_{0}}A_{ik}\left(Z_{i}\right)\left(B^{(i)}_{k1}\left(Z_{k}\right),B^{(i)}_{k2}\left(Z_{k}\right),\dots,B^{(i)}_{kN}\left(Z_{k}\right)\right), \label{bub1} 
\end{equation}
where $\left\{A_{ik}\left(Z_{i}\right)\right\}_{ 1\leq k,i\leq q_{0} }$ are linear forms in $Z_{i}$ and $\left\{B^{(i)}_{kl}\left(Z_{k}\right)\right\}_{1\leq k,i\leq q_{0}\atop{1\leq l \leq N}}$ are linear forms in $Z_{k}$, for all $k,i=1,\dots,q_{0}$.

We reiterate the previous computations.  Certain  conditions of compatibility are satisfied according to identifications of similar terms from substitutions of different coordinates in   (\ref{nuuit}), because
 $(k,i)\in \mathcal{S}$ if and only if $(k,j) \in \mathcal{S}$, for all $i,j=1,\dots, q$ with $i\neq j$. In particular: 
 
 The matrix $W_{0}$ from (\ref{dubluwe}) reduces in the coordinates $\left(V_{11}\otimes Z,U_{11}\otimes W\right)$ to the matrix 
 \begin{equation*}
\left(U_{11}\otimes W\right)_{0}=\begin{pmatrix}\left\langle \left(V_{11}\otimes Z\right)_{1},  \left(V_{11}\otimes Z\right)_{1}  \right\rangle & \left\langle \left(V_{11}\otimes Z\right)_{1}, \left(V_{11}\otimes Z\right)_{2}  \right\rangle & \dots &
\left\langle \left(V_{11}\otimes Z\right)_{1}, \left(V_{11}\otimes Z\right)_{q}  \right\rangle \\  \left\langle \left(V_{11}\otimes Z\right)_{2}, \left(V_{11}\otimes Z\right)_{1}  \right\rangle & \left\langle \left(V_{11}\otimes Z\right)_{2}, \left(V_{11}\otimes Z\right)_{2}  \right\rangle & \dots &
\left\langle \left(V_{11}\otimes Z\right)_{2}, \left(V_{11}\otimes Z\right)_{q}  \right\rangle \\ \vdots & \vdots & \ddots & \vdots  \\ \left\langle \left(V_{11}\otimes Z\right)_{q}, \left(V_{11}\otimes Z\right)_{1}  \right\rangle & \left\langle \left(V_{11}\otimes Z\right)_{q}, \left(V_{11}\otimes Z\right)_{2}  \right\rangle & \dots &
\left\langle \left(V_{11}\otimes Z\right)_{q}, \left(V_{11}\otimes Z\right)_{q}  \right\rangle \end{pmatrix}.
\end{equation*}

The matrix $\mathcal{W}$ from (\ref{dubluwe}) reduces in the coordinates $\left(V_{11}\otimes Z,U_{11}\otimes W\right)$  to the matrix 
 \begin{equation*}
    \hat{\mathcal{W}}:=\begin{pmatrix}\Re \left(U_{11}\otimes W\right)_{11} &   \left(U_{11}\otimes W\right)_{12} & \dots &
  \left(U_{11}\otimes W\right)_{1q} \\   \overline{\left(U_{11}\otimes W\right)}_{12} & \Re \left(U_{11}\otimes W\right)_{22} & \dots &
\overline{\left(U_{11}\otimes W\right)}_{2q} \\ \vdots & \vdots & \ddots & \vdots  \\  \overline{\left(U_{11}\otimes W\right)}_{1q} &   \overline{\left(U_{11}\otimes W\right)}_{2q} & \dots &
\Re \left(U_{11}\otimes W\right)_{qq}\end{pmatrix}.
\end{equation*}

For $k_{0}=i_{0}=1$, we obtain the following cases:
 \subsection{Case $\alpha=4m+1 $ } We write the the formal  expansion 
\begin{equation*}   \left( \left(\tilde{V}_{11}\otimes \left(f_{kl}\right)_{1\leq k\leq q\atop{1\leq l\leq 2N}}\right)^{(4m )}_{kl} \left(V_{11}\otimes Z,U_{11}\otimes W\right)\right)_{  1\leq k \leq q\atop{k_{0}+1\leq l  \leq N}}=\displaystyle\sum_{p=1}^{2m}\displaystyle\sum_{\substack{ J  \in\mathbb{N}^{q^{2}}    \\   \left|J \right| =2m-p  }}\left(\left(\tilde{P}_{ J}^{(2p)}\left(V_{11}\otimes Z\right)\right)_{kl}\right)_{  1\leq k \leq q\atop{1\leq l  \leq N}}
\left(U_{11}\otimes W\right)^{J}.  
\end{equation*}
  
 The analogue of the equation  (\ref{eKu1})  is the equation
\begin{equation} \begin{split}& 2\Re \left\{\left. \displaystyle\sum_{l=2}^{N}\overline{\left(V_{11}\otimes Z\right)}_{1l} \left(\tilde{V}_{11}\otimes \left(f_{kl}\right)_{1\leq k\leq q\atop{1\leq l\leq 2N}}\right)^{(4m)}_{1l}\left(V_{11}\otimes Z,U_{11}\otimes W
\right) \right\vert_{\left(Z,W\right)\in\mathcal{M}\atop{\left(V_{11}\otimes Z\right)_{11}=0}}  \right\}+\\&  2\Re\left\{ \left.\displaystyle\sum_{l=1}^{N}\left(\overline{\displaystyle\sum_{k,u=1\atop{(k,u)\in\mathcal{S}\atop{k,u\neq 1}}}^{q}a_{ku}^{1l} \left( \left(V_{11}\otimes Z\right)_{1l}
\right)_{i} \left(U_{11}\otimes W
\right)_{ku}} \right) \cdot\right.\right.\\& \left.\left.\quad\quad\quad\quad \quad\quad \quad\quad  \left(\tilde{V}_{11}\otimes \left(f_{kl}\right)_{1\leq k\leq q\atop{1\leq l\leq 2N}}\right)^{(4m-2)}_{1l}\left(V_{11}\otimes Z,U_{11}\otimes W
\right)  \right\vert_{\left(Z,W\right)\in\mathcal{M}\atop{\left(V_{11}\otimes Z\right)_{11}=0}} \right\}=0. \end{split}  \label{eKu1rRR1} \end{equation}

 \subsection{Case $\alpha=4m+3 $}We write the the formal expansion
 \begin{equation*}    \left( \left(\tilde{V}_{11}\otimes \left(f_{kl}\right)_{1\leq k\leq q\atop{1\leq l\leq 2N}}\right)^{(4m+2 )}_{kl} \left(V_{11}\otimes Z,U_{11}\otimes W\right)\right)_{  1\leq k \leq q\atop{ 1\leq l  \leq N}}=\displaystyle\sum_{p=1}^{2m+1}\displaystyle\sum_{\substack{ J  \in\mathbb{N}^{q^{2}}    \\   \left|J \right| =2m+1-p  }}\left(\left(\tilde{P}_{ J}^{(2p)}\left(V_{11}\otimes Z\right)\right)_{kl}\right)_{  1\leq k \leq q\atop{1\leq l  \leq N}}
\left(U_{11}\otimes W\right)^{J}.   
\end{equation*}

 The analogue of the equation (\ref{eKu1r}) is the equation
\begin{equation}\begin{split}& 2\Re \left\{\left. \displaystyle\sum_{l=2}^{N}\overline{\left(V_{11}\otimes Z\right)}_{1l} \left(\tilde{V}_{11}\otimes \left(f_{kl}\right)_{1\leq k\leq q\atop{1\leq l\leq 2N}}\right)^{(4m+2)}_{1l}\left(V_{11}\otimes Z,U_{11}\otimes W
\right) \right\vert_{\left(Z,W\right)\in\mathcal{M}\atop{\left(V_{11}\otimes Z\right)_{11}=0}}  \right\}+\\& 2\Re\left\{  \left. \displaystyle\sum_{l=1}^{N}\left(\overline{\displaystyle\sum_{k,u=1\atop{(k,u)\in\mathcal{S}\atop{k,u\neq 1}}}^{q}a_{ku}^{1l}\left(  \left(V_{11}\otimes Z\right)_{1l}
\right)_{i} \left(U_{11}\otimes W
\right)_{ku}} \right)\cdot\right.\right.\\& \left.\left.\quad\quad\quad\quad \quad\quad \quad\quad \quad\quad  \left(\tilde{V}_{11}\otimes \left(f_{kl}\right)_{1\leq k\leq q\atop{1\leq l\leq 2N}}\right)^{(4m)}_{1l}\left(V_{11}\otimes Z,U_{11}\otimes W
\right)\right\vert_{\left(Z,W\right)\in\mathcal{M}\atop{\left(V_{11}\otimes Z\right)_{11}=0}}  \right\}=0.\end{split}  \label{eKu1rRR2}\end{equation}

Comparing similar terms in   (\ref{eKu1rRR1}) and (\ref{eKu1rRR2}) for $W=\hat{\tilde{W}}$, it follows that
 \begin{equation} \begin{split}&  \left( \left(\left(\tilde{V}_{11}\otimes \left(f_{kl}\right)_{1\leq k\leq q\atop{1\leq l\leq 2N}}\right)^{(4m )}_{1l}\left(V_{11}\otimes Z,U_{11}\otimes W
\right)\right\vert_{\left(V_{11}\otimes Z\right)_{11}=0}  \right)_{ 2\leq l  \leq N } \quad\mbox{depending on}\\& \quad\quad\quad\quad\quad\quad\quad\quad\quad\quad\quad\quad\quad   \left( \left(\left(\tilde{V}_{11}\otimes \left(f_{kl}\right)_{1\leq k\leq q\atop{1\leq l\leq 2N}}\right)^{(4m-2 )}_{1l}\left(V_{11}\otimes Z,U_{11}\otimes W
\right)\right\vert_{\left(V_{11}\otimes Z\right)_{11}=0}  \right)_{ 2\leq l  \leq N }, \end{split} \label{lugu1}
\end{equation} 

 \begin{equation} \begin{split}&   \left( \left(\left(\tilde{V}_{11}\otimes \left(f_{kl}\right)_{1\leq k\leq k_{0}\atop{1\leq l\leq 2N}}\right)^{(4m+2 )}_{1l}\left(V_{11}\otimes Z,U_{11}\otimes W
\right)\right\vert_{\left(V_{11}\otimes Z\right)_{11}=0} \right)_{  2\leq l  \leq N }\quad\mbox{depending on}\\& \quad\quad\quad\quad\quad\quad\quad\quad\quad\quad\quad\quad\quad   \left( \left(\left(\tilde{V}_{11}\otimes \left(f_{kl}\right)_{1\leq k\leq q\atop{1\leq l\leq 2N}}\right)^{(4m-2 )}_{1l}\left(V_{11}\otimes Z,U_{11}\otimes W
\right)\right\vert_{\left(V_{11}\otimes Z\right)_{11}=0}  \right)_{ 2\leq l  \leq N }.\end{split} \label{lugu2}
\end{equation}

\subsection{Case $\alpha=4m+2$} We write the the formal expansions
 \begin{equation*}    \begin{split}& \left(  \left(\tilde{V}_{11}\otimes \left(f_{kl}\right)_{1\leq k\leq q\atop{1\leq l\leq 2N}}\right)^{(4m+1)}_{kl} \left(V_{11}\otimes Z,U_{11}\otimes W\right)\right)_{   1\leq k \leq q\atop{1\leq l  \leq 2N}}=\displaystyle\sum_{p=0}^{2m}\displaystyle\sum_{\substack{ J  \in\mathbb{N}^{q^{2}}    \\   \left|J \right| =2m-p  }}\left(\left(\tilde{P}_{ J}^{(2p+1)}\left(V_{11}\otimes Z\right)\right)_{kl}\right)_{  1\leq k \leq q\atop{1\leq l  \leq N}}
\left(U_{11}\otimes W\right)^{J},\\&  \left(  \left(\tilde{V}_{11}\otimes \left(f_{kl}\right)_{1\leq k\leq q\atop{1\leq l\leq 2N}}\right) ^{(2m+1)}_{kl}\left(V_{11}\otimes Z,U_{11}\otimes W\right)\right)_{   1\leq k \leq q\atop{N+1\leq l  \leq N}} =\displaystyle\sum_{p=0}^{m}\displaystyle\sum_{\substack{ J  \in\mathbb{N}^{q^{2}}    \\   \left|J \right| =m-p  }}\left(\left(\tilde{Q}_{ J}^{(2p+1)}\left(V_{11}\otimes Z\right)\right)_{kl}\right)_{   1\leq k \leq q\atop{1\leq l  \leq N}} 
 \left(U_{11}\otimes W\right)^{J}.\end{split}   
\end{equation*}

The analogue of the equation (\ref{ecu3}) is the equation
 \begin{equation}\begin{split}& 2\Re \left\{ \left.\displaystyle\sum_{l=1}^{N}\overline{\left(V_{11}\otimes Z\right)}_{1l} \left(\tilde{V}_{11}\otimes \left(f_{kl}\right)_{1\leq k\leq k_{0}\atop{1\leq l\leq 2N}}\right)^{(4m+1)}_{1l}\left(V_{11}\otimes Z,U_{11}\otimes W
\right)\right\vert_{\left(Z,W\right)\in\mathcal{M}\atop{{\left(V_{11}\otimes Z\right)_{11}=0}}} \right\}+\\&     \displaystyle\sum_{\alpha+\beta=4m+2\atop{l=1,\dots,N} }   \left( \left(\tilde{V}_{11}\otimes \left(f_{kl}\right)_{1\leq k\leq k_{0}\atop{1\leq l\leq 2N}}\right)^{(\alpha)}_{1\hspace{0.05 cm}1+l}\left(V_{11}\otimes Z,U_{11}\otimes W
\right)\right)\cdot   \\& \quad\quad\quad\quad\quad\quad\quad\quad\quad\quad        \left.\overline{\left( \left(\tilde{V}_{11}\otimes \left(f_{kl}\right)_{1\leq k\leq q\atop{1\leq l\leq 2N}}\right)^{(\beta)}_{1\hspace{0.05 cm}1+l}\left(V_{11}\otimes Z,U_{11}\otimes W
\right)\right)}   \right\vert_{\left(Z,W\right)\in\mathcal{M}\atop{{\left(V_{11}\otimes Z\right)_{11}=0}}}  + \\&     2\Re\left\{  \displaystyle\sum_{l=1}^{N}\left(\overline{\displaystyle\sum_{k,u=1\atop{(k,u)\in\mathcal{S}\atop{k,u\neq 1}}}^{q}a_{ku}^{1l}\left(  \left(V_{11}\otimes Z\right)_{1l}
\right)_{i} \left(U_{11}\otimes W
\right)_{ku}} \right)\right. \cdot\\& \quad\quad\quad\quad\quad\quad\quad\quad\quad \left.  \left. \left(\tilde{V}_{11}\otimes \left(f_{kl}\right)_{1\leq k\leq q\atop{1\leq l\leq 2N}}\right)^{(4m-3)}_{1l}\left(V_{11}\otimes Z,U_{11}\otimes W
\right)\right\vert_{\left(Z,W\right)\in\mathcal{M}\atop{{\left(V_{11}\otimes Z\right)_{11}=0}}} \right\} =0, \end{split}\label{ecu3ror} \end{equation}

Extracting multiples of terms of bi-degree $(1,1)$ in $\left(\left(V_{11}\otimes Z\right),\overline{\left(V_{11}\otimes Z\right)}\right)$ from (\ref{ecu3ror}), we obtain
\begin{equation} \begin{split}&        2\Re\left\{\left. \displaystyle\sum_{l=1}^{N} \displaystyle\sum_{\substack{ J  \in\mathbb{N}^{q^{2}}    \\   \left|J \right| =2m  }} \left(\tilde{P}_{ J}^{(1)}\left(V_{11}\otimes Z\right)\right)_{1l}\overline{\left(V_{11}\otimes Z\right)}_{1l} 
\hat{\tilde{W}}^{J}\right\vert_{\left(Z,W\right)\in\mathcal{M}\atop{\left(V_{11}\otimes Z\right)_{11}=0}} \right\}  +  \\&   \displaystyle\sum_{\alpha+\beta=2m\atop{l=1,\dots,N} }    \left(    \displaystyle\sum_{\substack{ J  \in\mathbb{N}^{q^{2}}    \\   \left|J \right| = m  }}\left(\tilde{Q}_{ J}^{(1)}\left(V_{11}\otimes Z\right)\right)_{1l} 
 \hat{\tilde{W}}^{J}  \right)\cdot      \left. \overline{ \left(      \displaystyle\sum_{\substack{ J  \in\mathbb{N}^{q^{2}}    \\   \left|J \right| = m  }} \left(\tilde{Q}_{ J}^{(1)}\left(V_{11}\otimes Z\right)\right)_{1l}  
\hat{\tilde{W}}^{J}  \right)}^{t}\right\vert_{\left(Z,W\right)\in\mathcal{M}\atop{\left(V_{11}\otimes Z\right)_{11}=0}} +  \\&  
      2\Re\left\{ \displaystyle\sum_{l=1}^{N} \displaystyle\sum_{\substack{ J  \in\mathbb{N}^{q^{2}}    \\   \left|J \right| =2m -3 }} \left(\tilde{P}_{ J}^{(1)}\left(V_{11}\otimes Z\right)\right)_{1l}\cdot \left. \left(\overline{\displaystyle\sum_{k,u=1\atop{(k,u)\in\mathcal{S}\atop{k,u\neq 1}}}^{q}a_{ku}^{1l} \left( \left(V_{11}\otimes Z\right)_{1l}
\right)_{i} \left(U_{11}\otimes W
\right)_{ku}} \right) 
\hat{\tilde{W}}^{J}\right\vert_{\left(Z,W\right)\in\mathcal{M}\atop{\left(V_{11}\otimes Z\right)_{11}=0}} \right\}   +\\& \quad\quad\quad\quad \left. \displaystyle\sum_{\alpha+\beta=2m-6\atop{l=1,\dots,N} }    \left(    \displaystyle\sum_{\substack{ J  \in\mathbb{N}^{q^{2}}    \\   \left|J \right| = \alpha }}\left(\tilde{Q}_{ J}^{(1)}\left(V_{11}\otimes Z\right)\right)_{1l}  
 \hat{\tilde{W}}^{J}  \right)\cdot        \overline{ \left(      \displaystyle\sum_{\substack{ J  \in\mathbb{N}^{q^{2}}    \\   \left|J \right| = \beta }} \left(\tilde{Q}_{ J}^{(1)}\left(V_{11}\otimes Z\right)\right)_{1l}   
\hat{\tilde{W}}^{J}  \right)}^{t}\right\vert_{\left(Z,W\right)\in\mathcal{M}\atop{\left(V_{11}\otimes Z\right)_{11}=0}}   
  =0. \end{split}  \label{suplim1ror}
 \end{equation}

Extracting multiples of terms of bi-degree $(3,3)$ in $\left(V_{11}\otimes Z,\overline{V_{11}\otimes Z}\right)$ from (\ref{ecu3ror}),  we obtain 
 
 \begin{equation}\begin{split} &        2\Re \left\{\left. \left( \displaystyle\sum_{l=1}^{N}     \displaystyle\sum_{\substack{ J  \in\mathbb{N}^{q^{2}}    \\   \left|J \right| =2m  }}\left(-n\left(J^{''}_{0}\right)\left(\tilde{P}_{ J}^{(1)}\left(V_{11}\otimes Z\right)\right)_{1l}\overline{\left(V_{11}\otimes Z\right)}_{1l} \left(U_{11}\otimes W\right)_{0}^{J^{''}_{0}}\right)  \hat{\tilde{W}}^{J ''} \right)
 \right\vert_{\left(Z,W\right)\in\mathcal{M}\atop{\left(V_{11}\otimes Z\right)_{11}=0}} \right\}  +\\&       \displaystyle\sum_{l=1}^{N}  \left(     \displaystyle\sum_{\substack{ J  \in\mathbb{N}^{q^{2}}    \\   \left|J \right| = m  }} n\left(J^{'}_{0}\right)\left(\tilde{Q}_{ J}^{(1)}\left(V_{11}\otimes Z\right)\right)_{1l} \left(U_{11}\otimes W\right)_{0}^{J^{'}_{0}}\hat{\tilde{W}}^{J'} \right)    \cdot    
 \left.\overline{\left(       \displaystyle\sum_{\substack{ J \in\mathbb{N}^{q^{2}}    \\   \left|J \right| = m  }}  \left(\tilde{Q}_{ J}^{(1)}\left(V_{11}\otimes Z\right)\right)_{1l}  
 \hat{\tilde{W}}^{J}  \right)}^{t}\right\vert_{\left(Z,W\right)\in\mathcal{M}\atop{\left(V_{11}\otimes Z\right)_{11}=0}}+  \\&    \left.  \displaystyle\sum_{l=1}^{N} \left(      \displaystyle\sum_{\substack{ J  \in\mathbb{N}^{q^{2}}    \\   \left|J \right| = m-1  }} \left(\tilde{Q}_{ J}^{(3)}\left(V_{11}\otimes Z\right)\right)_{1l}   
\hat{\tilde{W}}^{J}  \right)  \cdot    \overline{ \left(    \displaystyle\sum_{\substack{ J  \in\mathbb{N}^{q^{2}}    \\   \left|J \right| = m-1  }} \left(\tilde{Q}_{ J}^{(3)}\left(V_{11}\otimes Z\right)\right)_{1l}  
\hat{\tilde{W}}^{J}  \right)}^{t} \right\vert_{\left(Z,W\right)\in\mathcal{M}\atop{\left(V_{11}\otimes Z\right)_{11}=0}} +2\Re \left\{   \displaystyle\sum_{l=1}^{N}     \displaystyle\sum_{\substack{ J  \in\mathbb{N}^{q^{2}}    \\   \left|J \right| = 2m-3  }} \right.      \\&       \left(-n\left(J^{''}_{0}\right)\right)\left(\tilde{P}_{ J}^{(1)}\left(V_{11}\otimes Z\right)\right)_{1l}   \cdot      \left. \left.       \left.\left(\overline{\displaystyle\sum_{k,u=1\atop{(k,u)\in\mathcal{S}\atop{k,u\neq 1}}}^{q}a_{ku}^{1l}\left(  \left(V_{11}\otimes Z\right)_{1l}
\right)_{i} \left(U_{11}\otimes W
\right)_{ku}} \right)     \left(U_{11}\otimes W\right)_{0}^{J_{0}}  \hat{\tilde{W}}^{J ''} \right)\right\vert_{\left(Z,W\right)\in\mathcal{M}\atop{\left(V_{11}\otimes Z\right)_{11}=0}}
  \right\}  +\\&    \left.   \displaystyle\sum_{\alpha+\beta=2m-6\atop{l=1,\dots,N} }
      \left(    \displaystyle\sum_{\substack{ J  \in\mathbb{N}^{q^{2}}    \\   \left|J \right| = \alpha  }} n\left(J^{'}_{0}\right)\left(\tilde{Q}_{ J}^{(1)}\left(V_{11}\otimes Z\right)\right)_{1l} \left(U_{11}\otimes W\right)_{0}^{J_{0}}   \hat{\tilde{W}}^{J'}    \right)\cdot
       \overline{\left(       \displaystyle\sum_{\substack{ J \in\mathbb{N}^{q^{2}}    \\   \left|J \right| = \beta  }}\left( \left(\tilde{Q}_{ J}^{(1)}\left(V_{11}\otimes Z\right)\right)_{1l}\right) 
 \hat{\tilde{W}}^{J}  \right)}^{t}\right\vert_{\left(Z,W\right)\in\mathcal{M}\atop{\left(V_{11}\otimes Z\right)_{11}=0}} +  \\&     \quad\quad\quad\quad\quad\quad\quad\quad\quad\quad\quad \left. \displaystyle\sum_{\alpha+\beta=2m-8\atop{l=1,\dots,N} }
   \left(      \displaystyle\sum_{\substack{ J  \in\mathbb{N}^{q^{2}}    \\   \left|J \right| = \alpha  }} \left(\tilde{Q}_{ J}^{(3)}\left(V_{11}\otimes Z\right)\right)_{1l}    
\hat{\tilde{W}}^{J}  \right)  \cdot    \overline{ \left(    \displaystyle\sum_{\substack{ J  \in\mathbb{N}^{q^{2}}    \\   \left|J \right| = \beta  }} \left(\tilde{Q}_{ J}^{(3)}\left(V_{11}\otimes Z\right)\right)_{1l}   
 \hat{\tilde{W}}^{J}  \right)}^{t}\right\vert_{\left(Z,W\right)\in\mathcal{M}\atop{\left(V_{11}\otimes Z\right)_{11}=0}}=0  . \end{split} \label{suplim2ror}
 \end{equation}
 
Extracting multiples of terms of bi-degree $(2,2)$ in $\left(V_{11}\otimes Z,\overline{V_{11}\otimes Z}\right)$ from (\ref{ecu3ror}), we obtain
  \begin{equation}\begin{split}& 2\Re \left\{\left.\sqrt{-1} \left(     \displaystyle\sum_{l=1}^{N}  n\left(J^{'}_{0}\right)\displaystyle\sum_{\substack{ J  \in\mathbb{N}^{q^{2}}    \\   \left|J \right| =2m  }} \left(\tilde{P}_{ J}^{(1)}\left(V_{11}\otimes Z\right)\right)_{1l}\overline{\left(V_{11}\otimes Z\right)}_{1l}\left(U_{11}\otimes W\right)_{0}^{J^{'}_{0}}  \left(U_{11}\otimes \mathcal{W}\right)^{J^{'}}
   \right)\right\vert_{\left(Z,W\right)\in\mathcal{M}\atop{\left(V_{11}\otimes Z\right)_{11}=0}} \right\}+ \\& 2\Re \left\{\sqrt{-1} \left(     \displaystyle\sum_{l=1}^{N}  n\left(J^{'}_{0}\right)\displaystyle\sum_{\substack{ J  \in\mathbb{N}^{q^{2}}    \\   \left|J \right| =2m-3  }}     \left(\tilde{P}_{ J}^{(1)}\left(V_{11}\otimes Z\right)\right)_{1l}\cdot  \right.\right.\\& \left.\left.\quad\quad\quad\quad \quad\quad \hspace{0.1 cm}   \left.  \left(\overline{\displaystyle\sum_{k,u=1\atop{(k,u)\in\mathcal{S}\atop{k,u\neq 1}}}^{q}a_{ku}^{1l}\left(\left( \left(V_{11}\otimes Z\right)_{1l}
\right)_{i}\right) \cdot  \left(U_{11}\otimes W
\right)_{ku}} \right) \left(U_{11}\otimes W\right)_{0}^{J^{'}_{0}} \left(U_{11}\otimes \mathcal{W}\right)^{J'}
   \right)\right\vert_{\left(Z,W\right)\in\mathcal{M}\atop{\left(V_{11}\otimes Z\right)_{11}=0}} \right\}=0.\end{split}\label{pi}  \end{equation}

 Subtracting suitable terms between (\ref{suplim1ror}) and (\ref{suplim2ror}), we obtain 
 
  \begin{equation*}\begin{split}&      2\Re \left\{ \left.      \displaystyle\sum_{l=1}^{N}  \displaystyle\sum_{\substack{ J  \in\mathbb{N}^{q^{2}}    \\   \left|J \right| =2m  }} \left(\tilde{P}_{ J}^{(1)}\left(V_{11}\otimes Z\right)\right)_{kl}\left(\overline{V_{11}\otimes Z}\right)_{kl}    \hat{\tilde{W}}^{J}
  \right\vert_{\left(Z,W\right)\in\mathcal{M}\atop{\left(V_{11}\otimes Z\right)_{11}=0}} \right\} +  2\Re \left\{     \displaystyle\sum_{l=1}^{N}  \displaystyle\sum_{\substack{ J  \in\mathbb{N}^{q^{2}}    \\   \left|J \right| =2m-3  }}    \left(\tilde{P}_{ J}^{(1)}\left(V_{11}\otimes Z\right)\right)_{kl}\cdot\right. \\&     \quad\quad\quad\quad\quad\quad\quad\quad\quad\quad\quad\quad\quad\quad\quad\quad\quad\quad\left.   \left.  \left(\overline{\displaystyle\sum_{k,u=1\atop{(k,u)\in\mathcal{S}\atop{k,u\neq 1}}}^{q}a_{ku}^{1l}\left(\left( \left(V_{11}\otimes Z\right)_{1l}
\right)_{i}\right) \cdot  \left(U_{11}\otimes W
\right)_{ku}} \right)_{kl} \cdot  \hat{\tilde{W}}^{J}\right\vert_{\left(Z,W\right)\in\mathcal{M}\atop{\left(V_{11}\otimes Z\right)_{11}=0}}
     \right\}=0 .\end{split}   
 \end{equation*} 
 
In particular, we compute
\begin{equation}\begin{split}&   \left(  \left(\left.\left( \tilde{P}_{ J}^{( 1)}\left(V_{11}\otimes Z\right) \right)_{1l}\right\vert_{\left(V_{11}\otimes Z\right)_{11}=0} \right)_{  1\leq l  \leq N }\right)_{J  \in\mathbb{N}^{q^{2}}\atop{\left|J \right| =2m }}  \mbox{depending on} \\& \quad\quad\quad\quad\quad\quad\quad\quad\quad\quad\quad\quad\quad\quad
\quad\quad\quad\quad\quad\hspace{0.1 cm}  \left(  \left(\left.\left( \tilde{P}_{ J}^{( 1)}\left(V_{11}\otimes Z\right) \right)_{1l}\right\vert_{\left(V_{11}\otimes Z\right)_{11}=0} \right)_{  1\leq l  \leq N }\right)_{J  \in\mathbb{N}^{q^{2}}\atop{\left|J \right| <2m-3}}.\end{split} 
 \label{lugu3} 
\end{equation}

 Extracting multiples of terms of bi-degree $(3,1)$ in $\left(V_{11}\otimes Z,\overline{V_{11}\otimes Z}\right)$ from (\ref{ecu3ror}), we obtain
 \begin{equation*}\begin{split} & 2\Re \left\{\left. \displaystyle\sum_{l=1}^{N} \displaystyle\sum_{\substack{ J  \in\mathbb{N}^{q^{2}}    \\   \left|J \right| =2m-1  }} \left(\tilde{P}_{ J}^{(3)}\left(V_{11}\otimes Z\right)\right)_{kl}\overline{\left(V_{11}\otimes Z\right)}_{kl}   \hat{\tilde{W}}^{J}  \right\vert_{\left(Z,W\right)\in\mathcal{M}\atop{\left(V_{11}\otimes Z\right)_{11}=0}} \right\}+  2\Re \left\{ \displaystyle\sum_{l=1}^{N} \displaystyle\sum_{\substack{ J  \in\mathbb{N}^{q^{2}}    \\   \left|J \right| =2m-4  }}   \left(\tilde{P}_{ J}^{(3)}\left(V_{11}\otimes Z\right)\right)_{kl}\cdot \right.  \\& \quad\quad\quad\quad\hspace{0.1 cm}\quad\quad\quad\quad\quad\quad\quad\quad\quad\quad\quad\quad\quad\quad\quad\quad  \left. \left.\left(\overline{\displaystyle\sum_{k,u=1\atop{(k,u)\in\mathcal{S}\atop{k,u\neq 1}}}^{q}a_{ku}^{1l} \left( \left(V_{11}\otimes Z\right)_{1l}
\right)_{i}  \cdot  \left(U_{11}\otimes W
\right)_{ku}} \right)_{kl}  \cdot \hat{\tilde{W}}^{J}\right\vert_{\left(Z,W\right)\in\mathcal{M}\atop{\left(V_{11}\otimes Z\right)_{11}=0}}   \right\}=0,  \end{split}  \end{equation*}

Moving forward, it follows that
\begin{equation*}\begin{split}&  2\Re \left\{\left. \displaystyle\sum_{l=1}^{N} \displaystyle\sum_{\substack{ J  \in\mathbb{N}^{q^{2}}    \\   \left|J \right| =2m-2  }} \left(\tilde{P}_{ J}^{(5)}\left(V_{11}\otimes Z\right)\right)_{kl}\overline{\left(V_{11}\otimes Z\right)}_{kl}    \hat{\tilde{W}}^{J} \right\vert_{\left(Z,W\right)\in\mathcal{M}\atop{\left(V_{11}\otimes Z\right)_{11}=0}} \right\} +2\Re \left\{\displaystyle\sum_{l=1}^{N} \displaystyle\sum_{\substack{ J  \in\mathbb{N}^{q^{2}}    \\   \left|J \right| =2m-5  }}       \left(\tilde{P}_{ J}^{(5)}\left(V_{11}\otimes Z\right)\right)_{kl}\cdot  \right.\\&   \quad\quad\quad\quad\quad\quad\quad\quad\quad\quad\quad\quad\quad\quad\quad\quad\quad\quad   \left.  \left. \left(\overline{\displaystyle\sum_{k,u=1\atop{(k,u)\in\mathcal{S}\atop{k,u\neq 1}}}^{q}a_{ku}^{1l} \left( \left(V_{11}\otimes Z\right)_{1l}
\right)_{i}  \cdot  \left(U_{11}\otimes W
\right)_{ku}} \right)_{kl}   \cdot  \hat{\tilde{W}}^{J} \right\vert_{\left(Z,W\right)\in\mathcal{M}\atop{\left(V_{11}\otimes Z\right)_{11}=0}} \right\}+\\& \quad \hspace{0.03 cm}  \left.     \displaystyle\sum_{l=1}^{N}\left(\displaystyle\sum_{\substack{ J  \in\mathbb{N}^{q^{2}}    \\   \left|J \right| = m-2 }} \left(\tilde{Q}_{ J}^{(3)}\left(V_{11}\otimes Z\right)\right)_{kl}    \hat{\tilde{W}}\right)\cdot        \overline{ \left(\displaystyle\sum_{\substack{ J  \in\mathbb{N}^{q^{2}}    \\   \left|J \right| = m-2 }} \left(\tilde{Q}_{ J}^{(3)}\left(V_{11}\otimes Z\right)\right)_{kl}    \hat{\tilde{W}}^{J}\right)} \right\vert_{\left(Z,W\right)\in\mathcal{M}\atop{\left(V_{11}\otimes Z\right)_{11}=0}}+\\&   \left. \displaystyle\sum_{\alpha+\beta=2m-2\atop{l=1,\dots,N} } \left(\displaystyle\sum_{\substack{ J  \in\mathbb{N}^{q^{2}}    \\   \left|J \right| = \alpha }} \left(\tilde{Q}_{ J}^{(3)}\left(V_{11}\otimes Z\right)\right)_{kl}   \hat{\tilde{W}}\right)\cdot       \overline{ \left(\displaystyle\sum_{\substack{ J  \in\mathbb{N}^{q^{2}}    \\   \left|J \right| = \beta }} \left(\tilde{Q}_{ J}^{(3)}\left(V_{11}\otimes Z\right)\right)_{kl}    \hat{\tilde{W}}^{J}\right)}  \right\vert_{\left(Z,W\right)\in\mathcal{M}\atop{\left(V_{11}\otimes Z\right)_{11}=0}} =0, \end{split}   \end{equation*}
 $$\vdots\quad\quad\quad\quad\quad\quad\quad\quad\quad\quad\quad\quad    \vdots\quad\quad\quad\quad \quad\quad\quad\quad   \quad\quad\quad\quad \quad\quad\quad\quad\vdots\quad\quad\quad\quad\quad \quad\quad\quad \quad\quad\quad\quad\vdots $$

 In consequence, we compute

\begin{equation}\begin{split}& \left(\left(\left.  \left( \tilde{P}_{ J}^{(2p+1)}\left(V_{11}\otimes Z\right)   \right)_{1l}\right\vert_{\left(V_{11}\otimes Z\right)_{11}=0}\right)_{  1\leq l  \leq N }\right)_{1 \leq p\leq 2m\atop{J  \in\mathbb{N}^{q^{2}}\atop{\left|J \right| =2m-p}}}   \mbox{depending on} \\& \quad\quad\quad\quad\quad\quad\quad\quad\quad\quad\quad\quad\quad \quad\quad\quad \quad\quad\quad \quad \hspace{0.2 cm} \left(\left(\left.  \left( \tilde{P}_{ J}^{(2p+1)}\left(V_{11}\otimes Z\right)   \right)_{1l}\right\vert_{\left(V_{11}\otimes Z\right)_{11}=0}\right)_{  1\leq l  \leq N }\right)_{1 \leq p< 2m\atop{J  \in\mathbb{N}^{q^{2}}\atop{\left|J \right| <2m-p}}},    \\&  \left(  \left(\left.  \left( \tilde{Q}_{ J}^{(2p+1)}\left(V_{11}\otimes Z\right)   \right)_{1l}\right\vert_{\left(V_{11}\otimes Z\right)_{11}=0}\right)_{  1\leq l  \leq N }\right)_{1 \leq p\leq m\atop{J  \in\mathbb{N}^{q^{2}}\atop{\left|J \right| =m-p}}}  \mbox{depending on} \\& \quad\quad\quad\quad\quad\quad\quad\quad\quad\quad\quad\quad\quad \\&\quad\quad\quad  \quad\quad\quad  \quad\quad\quad  \quad\quad\quad  \quad\quad\quad  \quad\quad\quad \quad\quad     \left(  \left(\left.  \left( \tilde{Q}_{ J}^{(2p+1)}\left(V_{11}\otimes Z\right)   \right)_{1l}\right\vert_{\left(V_{11}\otimes Z\right)_{11}=0}\right)_{  1\leq l  \leq N }\right)_{1 \leq p< m\atop{J  \in\mathbb{N}^{q^{2}}\atop{\left|J \right| <m-p}}} .    \end{split}
 \label{lugu4}
\end{equation}   
  
 \subsection{Case $\alpha=4m+4$} We write the following formal expansions
  \begin{equation*}
 \begin{split}& \left(  \left(\tilde{V}_{11}\otimes \left(f_{kl}\right)_{1\leq k\leq q\atop{1\leq l\leq 2N}}\right)^{(4m+3)}_{kl} \left(V_{11}\otimes Z,U_{11}\otimes W\right)\right)_{   1\leq k \leq q\atop{1\leq l  \leq 2N}}
=\displaystyle\sum_{p=0}^{2m+1}\displaystyle\sum_{\substack{ J  \in\mathbb{N}^{q^{2}}    \\   \left|J \right| =2m+1-p  }}\left(\left(\tilde{P}_{ J}^{(2p+1)}\left(Z\right)\right)_{kl}\right)_{   1\leq k \leq q\atop{1\leq l  \leq N}}
 \left(U_{11}\otimes W\right)^{J}  ,\\&\left(  \left(\tilde{V}_{11}\otimes \left(f_{kl}\right)_{1\leq k\leq q\atop{1\leq l\leq 2N}}\right)^{(2m+2)}_{kl}\left(V_{11}\otimes Z,U_{11}\otimes W\right)\right)_{  1\leq k \leq q\atop{1\leq l  \leq N}}
 =\displaystyle\sum_{1=0}^{m+1}\displaystyle\sum_{\substack{ J  \in\mathbb{N}^{q^{2}}    \\   \left|J \right| =m+1-p  }}\left(\left(\tilde{Q}_{ J}^{(2p)}\left(V_{11}\otimes Z\right)\right)_{kl}\right)_{  1\leq k \leq q\atop{1\leq l  \leq N}}
 \left(U_{11}\otimes W\right)^{J} .\end{split} 
 \end{equation*}

 The analogue of the equation (\ref{ecu56}) is the equation
\begin{equation}\begin{split}& \left.2 \Re \left\{ \displaystyle\sum_{l=1}^{N}\overline{\left(V_{11}\otimes Z\right)}_{1l} \left(\tilde{V}_{11}\otimes \left(f_{kl}\right)_{1\leq k\leq q\atop{1\leq l\leq wN}}\right)^{(4m+3)}_{1l}\left(V_{11}\otimes Z,U_{11}\otimes W
\right)\right\vert_{\left(Z,W\right)\in\mathcal{M}\atop{\left(V_{11}\otimes Z\right)_{11}=0}} \right\} +\\&     \displaystyle\sum_{\alpha+\beta=4m+4 \atop{l=1,\dots,N}}   \left( \left(\tilde{V}_{11}\otimes \left(f_{kl}\right)_{1\leq k\leq q\atop{1\leq l\leq 2N}}\right)^{(\alpha)}_{1\hspace{0.05 cm}1+l}\left(V_{11}\otimes Z,U_{11}\otimes W
\right)\right) \cdot\\& \quad\quad\quad\quad\quad\quad\quad\quad\quad\quad \left.\overline{\left( \left(\tilde{V}_{11}\otimes \left(f_{kl}\right)_{1\leq k\leq q\atop{1\leq l\leq 2N}}\right)^{(\beta)}_{1\hspace{0.05 cm}1+l}\left(V_{11}\otimes Z,U_{11}\otimes W
\right)\right)}   \right\vert_{\left(Z,W\right)\in\mathcal{M}\atop{{\left(V_{11}\otimes Z\right)_{11}=0}}}    + \\&         \Re\left\{ \displaystyle\sum_{l=1}^{N}\left(\overline{\displaystyle\sum_{k,u=1\atop{(k,u)\in\mathcal{S}\atop{k,u\neq 1}}}^{q}a_{ku}^{1l}\left(\left( \left(V_{11}\otimes Z\right)_{1l}
\right)_{i}\right)\left(U_{11}\otimes W
\right)_{ku}} \right) \cdot \right.\\& \quad\quad\quad\quad\quad\quad\quad\quad\quad\quad \left.\left. \left(\tilde{V}_{11}\otimes \left(f_{kl}\right)_{1\leq k\leq q\atop{1\leq l\leq 2N}}\right)^{(4m-3)}_{1l}\left(V_{11}\otimes Z,U_{11}\otimes W
\right)   \right\vert_{\left(Z,W\right)\in\mathcal{M}\atop{{\left(V_{11}\otimes Z\right)_{11}=0}}}\right\}  =0. \end{split}\label{ecu56ror} \end{equation}

 Extracting multiples of terms of bi-degree $(1,1)$ in $\left(\left(V_{11}\otimes Z\right),\overline{\left(V_{11}\otimes Z\right)}\right)$ from (\ref{ecu56ror}),  it follows that 
  
  \begin{equation*}\begin{split}&   2\Re \left\{  \left.  \displaystyle\sum_{l=1}^{N}\left( \displaystyle\sum_{\substack{ J  \in\mathbb{N}^{q^{2}}    \\   \left|J \right| =2m+1  }}\left( \tilde{P}_{ J}^{(1)}\left(V_{11}\otimes Z\right)\right)_{1l}\overline{\left(V_{11}\otimes Z\right)}_{1l} 
 \hat{\tilde{W}}^{J}\right) \right\vert_{\left(Z,W\right)\in\mathcal{M}\atop{\left(V_{11}\otimes Z\right)_{11}=0}} \right\}   +     
    2\Re \left\{    \displaystyle\sum_{l=1}^{N}\left( \displaystyle\sum_{\substack{ J  \in\mathbb{N}^{q^{2}}    \\   \left|J \right| =2m-2  }} \left( \tilde{P}_{ J}^{(1)}\left(V_{11}\otimes Z\right)\right)_{1l}\cdot \right.\right.\\& \quad\quad\quad\quad\quad\quad\quad\quad\quad\quad\quad\quad\quad\quad\quad\quad\quad\quad\quad\left.\left. \left.\left(\overline{\displaystyle\sum_{k,u=1\atop{(k,u)\in\mathcal{S}\atop{k,u\neq 1}}}^{q}a_{ku}^{1l} \left( \left(V_{11}\otimes Z\right)_{1l}
\right)_{i}   \left(U_{11}\otimes W
\right)_{ku}} \right)  
 \left(U_{11}\otimes \tilde{W}\right)^{J}\right) \right\vert_{\left(Z,W\right)\in\mathcal{M}\atop{\left(V_{11}\otimes Z\right)_{11}=0}} \right\}  =0.\end{split}
 \end{equation*}

Extracting multiples of terms of bi-degree $(2,2)$ in $\left(V_{11}\otimes Z,\overline{V_{11}\otimes Z}\right)$ from (\ref{ecu56ror}),  we obtain
   
   \begin{equation}\begin{split}&   \left. 2\Re \left\{  \sqrt{-1} \displaystyle\sum_{l=1}^{N}\left(  \displaystyle\sum_{\substack{ J   \in\mathbb{N}^{q^{2}}    \\   \left|J  \right| =2m+1  }} n\left(J^{ '}_{0}\right)\left(\tilde{P}_{ J}^{(1)}\left(V_{11}\otimes Z\right)\right)_{1l}\overline{\left(V_{11}\otimes Z\right)}_{1l} \left(U_{11}\otimes W\right)_{0}^{J^{ '}_{0}} 
\hat{\tilde{W}}^{J^{ '}}\right)\right\vert_{\left(Z,W\right)\in\mathcal{M}\atop{\left(V_{11}\otimes Z\right)_{11}=0}}  \right\} +  \\&     \left. \displaystyle\sum_{l=1}^{N}\left(      \displaystyle\sum_{\substack{ J  \in\mathbb{N}^{q^{2}}    \\   \left|J \right| = m   }} \left(\tilde{Q}_{ J}^{(2)}\left(V_{11}\otimes Z\right)\right)_{1l}   
\hat{\tilde{W}}^{J}   \right)\cdot         \overline{ \left(     \displaystyle\sum_{\substack{ J  \in\mathbb{N}^{q^{2}}    \\   \left|J \right| = m   }}  \left(\tilde{Q}_{ J}^{(2)}\left(V_{11}\otimes Z\right)\right)_{1l}   
 \hat{\tilde{W}}^{J}  \right)}^{t}\right\vert_{\left(Z,W\right)\in\mathcal{M}\atop{\left(V_{11}\otimes Z\right)_{11}=0}}     
   +   2\Re \left\{  \sqrt{-1} \displaystyle\sum_{l=1}^{N}     \displaystyle\sum_{\substack{ J   \in\mathbb{N}^{q^{2}}    \\   \left|J  \right| =2m -2  }} \right. \\&   n\left(J_{0}\right)\left( \tilde{P}_{ J}^{(1)}\left(V_{11}\otimes Z\right)\right)_{1l}\cdot      \left.    \left.\left(\overline{\displaystyle\sum_{k,u=1\atop{(k,u)\in\mathcal{S}\atop{k,u\neq 1}}}^{q}a_{ku}^{1l} \left( \left(V_{11}\otimes Z\right)_{1l}
\right)_{i} \left(U_{11}\otimes W
\right)_{ku}} \right)    \left(U_{11}\otimes W\right)_{0}^{J^{ '}_{0}} 
\hat{\tilde{W}}^{J^{ '}}  \right\}\right\vert_{\left(Z,W\right)\in\mathcal{M}\atop{\left(V_{11}\otimes Z\right)_{11}=0}}  + \\& \quad\quad\quad\quad\quad\hspace{0.1 cm}   \left. \displaystyle\sum_{\alpha+\beta=2m\atop{l=1,\dots,N} }
 \left(      \displaystyle\sum_{\substack{ J  \in\mathbb{N}^{q^{2}}    \\   \left|J \right| =\alpha }} \left(\tilde{Q}_{ J}^{(2)}\left(V_{11}\otimes Z\right)\right)_{1l}  
\hat{\tilde{W}}^{J}   \right)\right. \cdot \left.         \overline{ \left(     \displaystyle\sum_{\substack{ J  \in\mathbb{N}^{q^{2}}    \\   \left|J \right| = \beta }}  \left(\tilde{Q}_{ J}^{(2)}\left(V_{11}\otimes Z\right)\right)_{1l}   
 \hat{\tilde{W}}^{J}  \right)}^{t}\right\vert_{\left(Z,W\right)\in\mathcal{M}\atop{\left(V_{11}\otimes Z\right)_{11}=0}}      
   =0. \end{split}  \label{sup1ror} 
 \end{equation}

Extracting multiples of terms of bi-degree $(4,4)$ in $\left(V_{11}\otimes Z,\overline{V_{11}\otimes Z}\right)$ from (\ref{ecu56ror}),  we obtain
 \begin{equation}\begin{split}&   \left. 2\Re \left\{  \displaystyle\sum_{l=1}^{N} \left(   \displaystyle\sum_{\substack{ J   \in\mathbb{N}^{q^{2}}    \\   \left|J  \right| =2m+1  }}\left(-n\left(J^{'''}_{0}\right)\right)\sqrt{-1}\left(\tilde{P}_{ J }^{(1)}\left(V_{11}\otimes Z\right)\right)_{1l}\overline{\left(V_{11}\otimes Z\right)}_{1l} \left(U_{11}\otimes W\right)_{0}^{J^{'''}_{0}}  \hat{\tilde{W}}^{J'''}  \right) 
  \right\}\right\vert_{\left(Z,W\right)\in\mathcal{M}\atop{\left(V_{11}\otimes Z\right)_{11}=0}} +\\&     \displaystyle\sum_{l=1}^{N}\left(        \displaystyle\sum_{\substack{ J  \in\mathbb{N}^{q^{2}}    \\   \left|J \right| = m   }}n\left(J^{'}_{0}\right)\left(\tilde{Q}_{ J}^{(2)}\left(V_{11}\otimes Z\right)\right)_{1l} \left(U_{11}\otimes W\right)_{0}^{J^{'}_{0}} \hat{\tilde{W}}^{J'}  \right)    \cdot \left.       \overline{   \left(    \displaystyle\sum_{\substack{ J ' \in\mathbb{N}^{q^{2}}    \\   \left|J \right| = m  }} \left(\tilde{Q}_{ J}^{(2)}\left(V_{11}\otimes Z\right)\right)_{1l} 
\hat{\tilde{W}}^{J}  \right) }^{t}\right\vert_{\left(Z,W\right)\in\mathcal{M}\atop{\left(V_{11}\otimes Z\right)_{11}=0}}  +\\&       \displaystyle\sum_{l=1}^{N}\left(     \displaystyle\sum_{\substack{ J  \in\mathbb{N}^{q^{2}}    \\   \left|J \right| = m -1 }}\left(\tilde{Q}_{ J}^{(4)}\left(V_{11}\otimes Z\right)\right)_{1l}  
     \hat{\tilde{W}}^{J}  \right) \cdot    \left.        
    \overline{ \left(      \displaystyle\sum_{\substack{ J  \in\mathbb{N}^{q^{2}}    \\   \left|J \right| = m -1 }}\left(\tilde{Q}_{ J}^{(4)}\left(V_{11}\otimes Z\right)\right)_{1l} \hat{\tilde{W}}^{J}   \right)}^{t}\right\vert_{\left(Z,W\right)\in\mathcal{M}\atop{\left(V_{11}\otimes Z\right)_{11}=0}}  + 2\Re \left\{\sqrt{-1}  \displaystyle\sum_{l=1}^{N}    \displaystyle\sum_{\substack{ J   \in\mathbb{N}^{q^{2}}    \\   \left|J  \right| =2m-2  }} \right. \\&       \left(-n\left(J^{ '''}_{0}\right)\right)\left(\tilde{P}_{ J }^{(1)}\left(V_{11}\otimes Z\right)\right)_{1l}  \cdot     \left.    \left.\left(\overline{\displaystyle\sum_{k,u=1\atop{(k,u)\in\mathcal{S}\atop{k,u\neq 1}}}^{q}a_{ku}^{1l}\left( \left(V_{11}\otimes Z\right)_{1l}
\right)_{i}\left(U_{11}\otimes W
\right)_{ku}} \right)  \left(U_{11}\otimes W\right)_{0}^{J^{'''}_{0}}   \hat{\tilde{W}}^{J'''}  
  \right\}\right\vert_{\left(Z,W\right)\in\mathcal{M}\atop{\left(V_{11}\otimes Z\right)_{11}=0}} +\\&          \displaystyle\sum_{\alpha+\beta=2m-4\atop{l=1,\dots,N}}
  \left(        \displaystyle\sum_{\substack{ J  \in\mathbb{N}^{q^{2}}    \\   \left|J \right| = \alpha  }}n\left(J^{'}_{0}\right)\left(\tilde{Q}_{ J}^{(2)}\left(V_{11}\otimes Z\right)\right)_{1l} \left(U_{11}\otimes W\right)_{0}^{J^{'}_{0}}    \hat{\tilde{W}}^{J'}  \right)    \cdot 
                     \left.       \overline{   \left(    \displaystyle\sum_{\substack{ J ' \in\mathbb{N}^{q^{2}}    \\   \left|J \right| = \beta }} \left(\tilde{Q}_{ J}^{(2)}\left(V_{11}\otimes Z\right)\right)_{1l}
                    \hat{\tilde{W}}^{J}  \right) }^{t}\right\vert_{\left(Z,W\right)\in\mathcal{M}\atop{\left(V_{11}\otimes Z\right)_{11}=0}}  +\\& \quad\quad\quad\quad\quad\quad\quad\quad\quad\quad   \displaystyle\sum_{\alpha+\beta=2m-2\atop{l=1,\dots,N}}
 \left(     \displaystyle\sum_{\substack{ J  \in\mathbb{N}^{q^{2}}    \\   \left|J \right| = \alpha }}\left(\tilde{Q}_{ J}^{(4)}\left(V_{11}\otimes Z\right)\right)_{1l}   
     \hat{\tilde{W}}^{J}  \right) \cdot      \left.        
    \overline{ \left(      \displaystyle\sum_{\substack{ J  \in\mathbb{N}^{q^{2}}    \\   \left|J \right| = \beta }}\left(\tilde{Q}_{ J}^{(4)}\left(V_{11}\otimes Z\right)\right)_{1l}  
\hat{\tilde{W}}^{J}   \right)}^{t}\right\vert_{\left(Z,W\right)\in\mathcal{M}\atop{\left(V_{11}\otimes Z\right)_{11}=0}}   =0. \end{split}\label{supror}
 \end{equation}

 We identify similar terms (\ref{sup1ror}) and (\ref{supror}). We compute
 
 \begin{equation}\begin{split}&  \left(  \left(\left(\left. \tilde{P}_{ J}^{( 1)}\left(V_{11}\otimes Z\right)\right\vert_{\left(V_{11}\otimes Z\right)_{11}=0} \right)_{1l} \right)_{ 1\leq l  \leq N }\right)_{J  \in\mathbb{N}^{q^{2}}\atop{\left|J \right| =2m+1}}  \quad\mbox{depending on}\\&\quad\quad\quad  \quad\quad\quad  \quad\quad\quad  \quad\quad\quad  \quad\quad\quad  \quad\quad\quad  \quad\quad\quad  \quad\quad    \left(  \left(\left(\left. \tilde{P}_{ J}^{( 1)}\left(V_{11}\otimes Z\right)\right\vert_{\left(V_{11}\otimes Z\right)_{11}=0} \right)_{1l} \right)_{ 1\leq l  \leq N }\right)_{J  \in\mathbb{N}^{q^{2}}\atop{\left|J \right| <2m+1}}. \end{split}
  \label{lugu4}
\end{equation}

Returning to the last first in order to identify similar terms, we compute

  \begin{equation}\begin{split}&  \left( \left(\left( \left.  \tilde{Q}_{ J}^{( 1)}\left(V_{11}\otimes Z\right)\right\vert_{\left(V_{11}\otimes Z\right)_{11}=0}   \right)_{1l}\right)_{ 1\leq l  \leq N}\right)_{J  \in\mathbb{N}^{q^{2}}\atop{\left|J \right| =m+1}}  \quad\mbox{depending on} \\&\quad\quad\quad  \quad\quad\quad  \quad\quad\quad  \quad\quad  \quad\quad\quad  \quad\quad\quad  \quad\quad\quad  \quad\quad  \left( \left(\left( \left.  \tilde{Q}_{ J}^{( 1)}\left(V_{11}\otimes Z\right)\right\vert_{\left(V_{11}\otimes Z\right)_{11}=0}   \right)_{1l}\right)_{ 1\leq l  \leq N}\right)_{J  \in\mathbb{N}^{q^{2}}\atop{\left|J \right| <m+1}}. \end{split}
 \label{lugu5}
\end{equation}

Extracting multiples of terms of bi-degree $(4,1)$ in $\left(V_{11}\otimes Z,\overline{V_{11}\otimes Z}\right)$ from (\ref{ecu56ror}), we obtain

\begin{equation*}\begin{split}&  2\Re \left\{\left. \displaystyle\sum_{l=1}^{N} \displaystyle\sum_{\substack{ J  \in\mathbb{N}^{q^{2}}    \\   \left|J \right| =2m  }}\left(\tilde{P}_{ J}^{(3)}\left(V_{11}\otimes Z\right)\right)_{kl}\overline{\left(V_{11}\otimes Z\right)}_{kl}   \hat{\tilde{W}}
^{J}\right\vert_{\left(Z,W\right)\in\mathcal{M}\atop{\left(V_{11}\otimes Z\right)_{11}=0}}   \right\} +    2\Re \left\{ \displaystyle\sum_{l=1}^{N} \displaystyle\sum_{\substack{ J  \in\mathbb{N}^{q^{2}}    \\   \left|J \right| =2m  }}     \left(\tilde{P}_{ J}^{(3)}\left(V_{11}\otimes Z\right)\right)_{kl}\cdot \right.\\&\quad\quad\quad\quad\quad\quad\quad\quad\quad\quad\quad\quad\quad\quad\quad\hspace{0.15 cm} \left.  \left. \left(\overline{\displaystyle\sum_{k,u=1\atop{(k,u)\in\mathcal{S}\atop{k,u\neq 1}}}^{q}a_{ku}^{1l} \left( \left(V_{11}\otimes Z\right)_{1l}
\right)_{i}  \cdot  \left(U_{11}\otimes W
\right)_{ku}} \right)_{kl}    \cdot \hat{\tilde{W}}
^{J}  \right\vert_{\left(Z,W\right)\in\mathcal{M}\atop{\left(V_{11}\otimes Z\right)_{11}=0}}\right\} +\\& \left.    \displaystyle\sum_{l=1}^{N}\left(\displaystyle\sum_{\substack{ J  \in\mathbb{N}^{q^{2}}    \\   \left|J \right| = m  }} \left(\tilde{Q}_{ J}^{(2)}\left(V_{11}\otimes Z\right)\right)_{kl}   \hat{\tilde{W}}
^{J}\right)\cdot        \overline{ \left(\displaystyle\sum_{\substack{ J  \in\mathbb{N}^{q^{2}}    \\   \left|J \right| = m   }} \left(\tilde{Q}_{ J}^{(2)}\left(V_{11}\otimes Z\right)\right)_{kl}     \hat{\tilde{W}}
^{J}\right)} \right\vert_{\left(Z,W\right)\in\mathcal{M}\atop{\left(V_{11}\otimes Z\right)_{11}=0}}  +\\&   \left.   \displaystyle\sum_{\alpha+\beta=2m\atop{l=1,\dots,N}}\left(\displaystyle\sum_{\substack{ J  \in\mathbb{N}^{q^{2}}    \\   \left|J \right| = \alpha  }} \left(\tilde{Q}_{ J}^{(2)}\left(V_{11}\otimes Z\right)\right)_{kl}   \hat{\tilde{W}}
^{J}\right)\cdot        \overline{ \left(\displaystyle\sum_{\substack{ J  \in\mathbb{N}^{q^{2}}    \\   \left|J \right| = \beta   }} \left(\tilde{Q}_{ J}^{(2)}\left(V_{11}\otimes Z\right)\right)_{kl}  \hat{\tilde{W}}
^{J}\right)}\right\vert_{\left(Z,W\right)\in\mathcal{M}\atop{\left(V_{11}\otimes Z\right)_{11}=0}} =0. \end{split}   \end{equation*}

Moving forward, it follows that
\begin{equation*} \begin{split}&  2\Re \left\{ \left.\displaystyle\sum_{l=1}^{N} \displaystyle\sum_{\substack{ J  \in\mathbb{N}^{q^{2}}    \\   \left|J \right| =2m -1 }} \left(\tilde{P}_{ J}^{(5)}\left(V_{11}\otimes Z\right)\right)_{kl}\overline{\left(V_{11}\otimes Z\right)}_{kl}    \hat{\tilde{W}}
^{J}  \right\vert_{\left(Z,W\right)\in\mathcal{M}\atop{\left(V_{11}\otimes Z\right)_{11}=0}} \right\} +  2\Re \left\{\displaystyle\sum_{l=1}^{N} \displaystyle\sum_{\substack{ J  \in\mathbb{N}^{q^{2}}    \\   \left|J \right| =2m -4 }}    \left(\tilde{P}_{ J}^{(5)}\left(V_{11}\otimes Z\right)\right)_{kl} \cdot   \right.\\&  \quad\quad\quad\quad\quad\quad\quad\quad\quad\quad\quad\quad\quad\quad\quad\quad\quad\quad\quad \left. \left.  \left(\overline{\displaystyle\sum_{k,u=1\atop{(k,u)\in\mathcal{S}\atop{k,u\neq 1}}}^{q}a_{ku}^{1l} \left( \left(V_{11}\otimes Z\right)_{1l}
\right)_{i}  \cdot  \left(U_{11}\otimes W
\right)_{ku}} \right)_{kl}     \hat{\tilde{W}}
^{J}\right\vert_{\left(Z,W\right)\in\mathcal{M}\atop{\left(V_{11}\otimes Z\right)_{11}=0}}  \right\}=0,\end{split}   \end{equation*}

  $$\vdots\quad\quad\quad\quad\quad\quad\quad\quad\quad\quad\quad\quad    \vdots\quad\quad\quad\quad \quad\quad\quad\quad \quad\quad\quad\quad \vdots\quad\quad\quad\quad \quad\quad  \quad\quad\quad \quad\quad\quad\quad\vdots $$

In consequence, we compute 
\begin{equation}\begin{split}&  \left(\left(  \left(\left.  \tilde{Q}_{ J}^{(2p)}\left(V_{11}\otimes Z\right) \right\vert_{\left(V_{11}\otimes Z\right)_{11}=0}  \right)_{1l}\right)_{ 1\leq l  \leq N} \right)_{1\leq p\leq  2m+1\atop{J  \in\mathbb{N}^{q^{2}}\atop{\left|J \right| =2m+1-p}}}  \mbox{depending on}\\& \quad\quad\quad\quad\quad\quad\quad\quad\quad\quad\quad\quad\quad
\quad\quad\quad\quad\quad\quad\quad\quad\quad\quad   \left(\left(  \left(\left.  \tilde{Q}_{ J}^{(2p)}\left(V_{11}\otimes Z\right) \right\vert_{\left(V_{11}\otimes Z\right)_{11}=0}  \right)_{1l}\right)_{ 1\leq l  \leq N} \right)_{1\leq p<  2m+1\atop{J  \in\mathbb{N}^{q^{2}}\atop{\left|J \right| <2m+1-p}}}, \\& \left(\left(\left( \left.\tilde{P}_{ J}^{(2p+1)}\left(V_{11}\otimes Z\right) \right\vert_{\left(V_{11}\otimes Z\right)_{11}=0}\right)_{1l} \right)_{  1\leq l  \leq N}\right)_{1 \leq p\leq  m\atop{J  \in\mathbb{N}^{q^{2}}\atop{\left|J \right| =m+1-p}}}  \mbox{depending on} \\& \quad\quad\quad\quad\quad\quad\quad\quad\quad\quad\quad\quad\quad
\quad\quad\quad\quad\quad\quad\quad\quad\quad\quad  \left(\left(\left( \left.\tilde{P}_{ J}^{(2p+1)}\left(V_{11}\otimes Z\right) \right\vert_{\left(V_{11}\otimes Z\right)_{11}=0}\right)_{1l} \right)_{  1\leq l  \leq N}\right)_{1 \leq p<  m\atop{J  \in\mathbb{N}^{q^{2}}\atop{\left|J \right| <m+1-p}}}.  \end{split}
 \label{lugu6}
\end{equation} 
 
 We summarize  (\ref{lugu1}),(\ref{lugu2}),(\ref{lugu3}),(\ref{lugu4}),(\ref{lugu5})  and  (\ref{lugu6}) :
\subsection{Formal Power Series} We consider the following formal expansions
$$    \left(\tilde{V}_{11}\otimes \left(f_{kl}\right)_{1\leq k\leq q\atop{1\leq l\leq 2N}}\right)_{12} \left(V_{11}\otimes Z,U_{11}\otimes W\right)+A_{2}\left(V_{11}\otimes Z,U_{11}\otimes W\right) =\left(V_{11}\otimes Z\right)
_{12}+\left(V_{11}\otimes Z\right)
_{11}\hat{f}_{12}\left(V_{11}\otimes Z,U_{11}\otimes W\right),$$  
 $$\vdots\quad\quad\quad\quad\quad\quad\quad\quad\quad\quad\quad\quad    \vdots\quad\quad\quad\quad \quad\quad\quad\quad \quad\quad\quad\quad \vdots\quad\quad\quad\quad \quad\quad\quad \quad \quad\quad\quad \quad\quad\quad\quad\vdots $$
 $$ \left(\tilde{V}_{11}\otimes \left(f_{kl}\right)_{1\leq k\leq q\atop{1\leq l\leq 2N}}\right)_{1N} \left(V_{11}\otimes Z,U_{11}\otimes W\right)+A_{N}\left(V_{11}\otimes Z,U_{11}\otimes W\right) =\left(V_{11}\otimes Z\right)
_{1N}+\left(V_{11}\otimes Z\right)
_{11}\hat{f}_{1N}\left(V_{11}\otimes Z,U_{11}\otimes W\right),$$
where $A_{2}\left(V_{11}\otimes Z,U_{11}\otimes W\right),\dots,A_{N}\left(V_{11}\otimes Z,U_{11}\otimes W\right)$ are formal power series determined by induction, and
$$    \left(\tilde{V}_{11}\otimes \left(\varphi_{kl}\right)_{1\leq k\leq q\atop{1\leq l\leq 2N}}\right)_{1\hspace{0.05 cm}N+1} \left(V_{11}\otimes Z,U_{11}\otimes W\right)+B_{1}\left(V_{11}\otimes Z,U_{11}\otimes W\right) =\left(V_{11}\otimes Z\right)
_{11}\hat{\varphi}_{11}\left(V_{11}\otimes Z,U_{11}\otimes W\right),$$  $$\left(\tilde{V}_{11}\otimes \left(\varphi_{kl}\right)_{1\leq k\leq q\atop{1\leq l\leq 2N}}\right)_{1\hspace{0.05 cm}N+2} \left(V_{11}\otimes Z,U_{11}\otimes W\right)+B_{2}\left(V_{11}\otimes Z,U_{11}\otimes W\right) =\left(V_{11}\otimes Z\right)
_{12}\hat{\varphi}_{12}\left(V_{11}\otimes Z,U_{11}\otimes W\right),$$
   $$\vdots\quad\quad\quad\quad\quad\quad\quad\quad\quad\quad\quad\quad    \vdots\quad\quad\quad\quad  \quad\quad\quad\quad \quad\quad\quad\quad\vdots\quad\quad\quad\quad\quad \quad\quad\quad \quad\quad\quad\quad\vdots $$ 
$$ \left(\tilde{V}_{11}\otimes \left(\varphi_{kl}\right)_{1\leq k\leq q\atop{1\leq l\leq 2N}}\right)_{1\hspace{0.05 cm}2N} \left(V_{11}\otimes Z,U_{11}\otimes W\right)+B_{N}\left(V_{11}\otimes Z,U_{11}\otimes W\right) = \left(V_{11}\otimes Z\right)
_{11}\hat{\varphi}_{1N}\left(V_{11}\otimes Z,U_{11}\otimes W\right),$$
where $B_{1}\left(V_{11}\otimes Z,U_{11}\otimes W\right),B_{2}\left(V_{11}\otimes Z,U_{11}\otimes W\right),\dots,B_{N}\left(V_{11}\otimes Z,U_{11}\otimes W\right)$ are formal power series determined by induction.

We focus on      the coefficients  of    the  terms  of the following type
\begin{equation*}\begin{split}\left(V_{11}\otimes Z\right)
_{11}\left(V_{11}\otimes Z\right)
^{I}&\left(\mbox{Re} \left(U_{11}\otimes W
\right)_{11}\right)^{j_{11}} \left(U_{11}\otimes W
\right)_{12}^{j_{12}}\dots 
 \left(U_{11}\otimes W
\right)_{1q}^{j_{1q}}\cdot\\&\overline{\left(U_{11}\otimes W
\right)}_{12}^{j_{21}}\left( \Re \left(U_{11}\otimes W
\right)_{22}\right)^{j_{22}}\dots  \left(U_{11}\otimes W
\right)_{2q}^{j_{2q}}\cdot\\&\quad\quad \vdots\quad\quad\quad\quad\quad\quad\quad \vdots\quad \quad\quad   \quad\hspace{0.1 cm}\quad\ddots\quad \quad\quad \vdots \\& \overline{ \left(U_{11}\otimes W
\right)}_{1q}^{j_{q1}} \overline{  \left(U_{11}\otimes W
\right)}_{2q}^{j_{q2}}\dots \left(  \Re \left(U_{11}\otimes W
\right)_{qq}\right)^{j_{qq}}\cdot\overline{\left(V_{11}\otimes Z\right)}_{1l},\quad\mbox{for all $I\in\mathbb{N}^{qN}$, $l=1,\dots,N$ and $J\in\mathbb{N}^{q^{2}}$.}\end{split}\end{equation*}

In particular, we obtain 
$\hat{f}_{12}\left(V_{11}\otimes Z,U_{11}\otimes W\right)=\dots=\hat{f}_{1N}\left(V_{11}\otimes Z,U_{11}\otimes W\right)=0.$ Then,  we write the formal weighted expansions
\begin{equation*}\begin{split}&     \left(\hat{\varphi}_{1l}\left(V_{11}\otimes Z,U_{11}\otimes W\right)\right)_{1\leq l\leq N}=\displaystyle\sum_{\alpha\in\mathbb{N}^{\star}}\left(\theta^{\alpha}\left(V_{11}\otimes Z,U_{11}\otimes W\right)\right)_{1\leq l\leq N},\quad\mbox{provided the polynomial}\\& \left(\theta^{\alpha}\left(V_{11}\otimes Z,U_{11}\otimes W\right)\right)_{1\leq l\leq N}=\left(\displaystyle\sum_{\substack{ J  \in\mathbb{N}^{q^{2}}, \hspace{0.05 cm} I  \in\mathbb{N}^{qN} \\ \left|I \right|+2\left|J \right| =\alpha  }} \theta_{ l}^{ I J  } \left(V_{11}\otimes Z\right)^{I} \left(U_{11}\otimes W\right)^{J}\right)_{1\leq l\leq N}. \end{split}
\end{equation*}

Provided $\alpha\in\mathbb{N}^{\star}$ such that $\alpha\neq 1,2$, we assume
\begin{equation*} \left(\varphi_{kl}^{\left(\alpha-1\right)}\left(Z,W\right)\right)_{1\leq k  \leq q\atop{1\leq l  \leq N}}=\mbox{O}_{q\times N}\hspace{0.1 cm}\mbox{and}   \hspace{0.1 cm} \left( f_{kl}^{\left(\alpha+1\right)}\left(Z,W\right)\right)_{1\leq k  \leq q\atop{1\leq l  \leq N}}=\mbox{O}_{q\times N}.
 \end{equation*}

In particular, we write
 $  f_{11}^{\left(\alpha+2\right)}\left(V_{11}\otimes Z, U_{11}\otimes W \right) =\left(V_{11}\otimes Z\right)
_{11}\cdot A\left(V_{11}\otimes Z, U_{11}\otimes W \right)+B\left(V_{11}\otimes Z, U_{11}\otimes W \right)$,  
where $A\left(V_{11}\otimes Z, U_{11}\otimes W \right)$ and $B\left(V_{11}\otimes Z, U_{11}\otimes W \right)$ are formal power series such that  $B\left(V_{11}\otimes Z, U_{11}\otimes W \right)$ is not a multiple of $\left(V_{11}\otimes Z\right)
_{11}$.

 We reiterate (\ref{eKu1rRR1}),(\ref{eKu1rRR2}) and (\ref{ecu3ror})   using   (\ref{lugu1}),(\ref{lugu2}),(\ref{lugu3}),(\ref{lugu4}),(\ref{lugu5})  and  (\ref{lugu6}). We obtain
 \begin{equation}\begin{split}& 2\Re \left\{\left. \left(\frac{i}{2}\left(U_{11}\otimes W \right)_{11}\right)\cdot \overline{  \left( A\left(V_{11}\otimes Z, U_{11}\otimes W \right)\right)}\right\vert_{\Im \left( U_{11}\otimes W\right)=\left( V_{11}\otimes Z\right)\cdot \left(\overline{ V_{11}\otimes Z}\right)^{t}}\right\}+\\& \left. 2 Re \left\{\displaystyle\sum_{l=1}^{N}\overline{\left(V_{11}\otimes Z\right)_{1l} }\cdot\left(  \displaystyle\sum_{\substack{ J  \in\mathbb{N}^{q^{2}}, \hspace{0.05 cm} I  \in\mathbb{N}^{qN} \\ \left|I \right|+2\left|J \right| =\alpha  }} \theta_{ l}^{ I J  } \left(V_{11}\otimes Z\right)^{I}\left(U_{11}\otimes W\right)^{J }\right)\right\vert_{\Im \left( U_{11}\otimes W\right)=\left( V_{11}\otimes Z\right)\cdot \left(\overline{ V_{11}\otimes Z}\right)^{t}}\right\}=0 .\end{split}\label{porrr}
\end{equation}
\begin{equation}2 \Re \left\{\left. \left(V_{11}\otimes Z\right)
_{11} \cdot\overline{\left(    B\left(V_{11}\otimes Z, U_{11}\otimes W \right)\right)}\right\vert_{\Im \left( U_{11}\otimes W\right)=\left( V_{11}\otimes Z\right)\cdot \left(\overline{ V_{11}\otimes Z}\right)^{t}}\right\}=0.\label{porrrt}
\end{equation} 
 
 \subsection{Proof of Theorem \ref{Teorem}}The vanishing of all terms is clear in (\ref{porrr}) and (\ref{porrrt}). More   generally, we obtain  the vanishing of the higher order terms   from the formal expansions  (\ref{dub1}) using induction, (\ref{tre}),(\ref{rezu111r}),(\ref{ecuA5se}),(\ref{ecuA5}),(\ref{ecuA3sec}),(\ref{ecuA8}), 
(\ref{ecuA8se}),(\ref{ecuA6se}) and    special  coordinates constructed according to Huang-Ji\cite{HJ}. 
 The  Classes   of Equivalence   (\ref{clase}) are  defined by the composition $\left(\tilde{G},\tilde{F}\right) 
\left(W,Z\right)= \left(\mathcal{C}'\circ \left( G,F\right) \circ   \mathcal{C}^{-1}\right)\left(W,Z\right)$, given the  Cayley type transformation  (\ref{Ca}),       the  Cayley type transformation of the target $\mathcal{BSD}$-Model denoted by $\mathcal{C}'$   
and   suitable compositions with   automorphisms of the Shilov Boundaries $S_{p,q}$ and $S_{p',q'}$. They extend to holomorphic automorphisms of the bounded symmetric domains $D_{p,q}$ and  $D_{p',q'}$
according to Kaup-Zaitsev\cite{KaZa1},\cite{KaZa2} and
 Kim-Zaitsev\cite{kz},\cite{kza}.

  \section{Acknowledgements}    Thanks to the Universities Ko\c{c}, Middle East Technical, Galatasaray and Hacettepe for   hospitality   prior to my deportation from Istanbul. Special Thanks to my (former) supervisor Prof. Dmitri Zaitsev for important conversations while I was reading \cite{kz} (and for his guidance  on  the main part\cite{V1}  (supported by  Science Foundation Ireland, Grant 06/RFP/MAT 018) of my doctoral thesis)).

\end{document}